\documentclass[review, sort&compress, 10pt]{elsarticle}

\usepackage[utf8]{inputenc}
\usepackage[T1]{fontenc}
\usepackage{amssymb}
\usepackage{amsthm}
\usepackage{graphicx}
\usepackage{subcaption}
\usepackage{amsmath}
\usepackage{amsfonts}
\usepackage{caption}
\usepackage{microtype}
\usepackage{csquotes}
\usepackage{xcolor}
\usepackage{hyperref}
\usepackage{tikz}
\usetikzlibrary{positioning,arrows.meta,calc,shapes.geometric,decorations.pathreplacing,patterns}
\usepackage{cleveref}
\usepackage{comment}
\usepackage{color}
\usepackage{mathrsfs} 
\usepackage{microtype}
\usepackage{algorithm}
\usepackage{algpseudocode}
\usepackage{pdflscape}
\usepackage[margin=2.25cm]{geometry}
\usepackage{nomencl}
\makenomenclature
\usepackage{framed}
\usepackage{multicol}
\usepackage{booktabs}

\renewcommand*\nompreamble{\begin{multicols}{2}}
	\renewcommand*\nompostamble{\end{multicols}}

\renewcommand\nomgroup[1]{%
	\item[\bfseries
	\ifstrequal{#1}{A}{Constants and parameters}{%
		\ifstrequal{#1}{B}{State variables}{%
			\ifstrequal{#1}{C}{Superscripts}{
				\ifstrequal{#1}{D}{Acronyms}{}}}}%
	]}

\theoremstyle{remark}
\newtheorem{remark}{Remark}
\journal{}

\begin{document}
	
	\begin{frontmatter}
		
		
		
		\title{Physics-informed post-processing of stabilized finite element solutions for transient convection-dominated problems}

		\author[Slymn1]{S\"uleyman Cengizci\corref{cor1}}
		\address[Slymn1]{Department of Computer Technologies, Antalya Bilim University, Antalya 07190, Turkey}
		\cortext[cor1]{Corresponding author}
		\ead{suleyman.cengizci@antalya.edu.tr}

		
		\author[omur]{Ömür Uğur}
		\address[omur]{Institute of Applied Mathematics, Middle East Technical University, Ankara 06800, Turkey}
		\ead{ougur@metu.edu.tr}
		
		\author[Natesan]{Srinivasan Natesan}
		\address[Natesan]{Department of Mathematics, Indian Institute of Technology Guwahati, Guwahati 781039, India}
		\ead{natesan@iitg.ac.in}

\begin{abstract}
	The numerical simulation of convection-dominated transient transport 
	phenomena involves significant computational challenges due to the 
	presence of sharp gradients and propagating fronts across the 
	spatiotemporal domain. Classical numerical discretization methods 
	often produce spurious oscillations, necessitating advanced 
	stabilization strategies. Moreover, even when stabilization 
	suppresses such numerical instabilities, additional regularization 
	may still be required to resolve localized sharp gradients. 
	Conversely, standalone physics-informed neural network (PINN) 
	methods struggle to capture sharp layers arising in 
	convection-dominated problems and typically require a prohibitively large 
	number of training epochs to learn the solution from scratch. 
	To this end, this work presents a hybrid computational 
	framework---extending the PASSC (PINN-Augmented SUPG with 
	Shock-Capturing) methodology, originally developed for steady 
	problems, to the unsteady regime---that combines stabilized finite 
	element methods (FEM) with PINNs for solving transient 
	convection-diffusion-reaction (CDR) equations. A semi-discrete 
	stabilized finite element method is adopted, where stabilization 
	is achieved through the streamline-upwind/Petrov--Galerkin (SUPG) 
	formulation augmented with a YZ$\beta$ shock-capturing operator to 
	resolve steep gradients and shock-like features. Rather than 
	training over the entire spatiotemporal domain, a PINN-based 
	post-processing correction is applied selectively near the 
	terminal time, where the neural network enhances the finite 
	element solution by assimilating spatial data from the last $K_s$ 
	snapshots while enforcing residual constraints derived from the 
	governing equations and boundary conditions. The network 
	architecture employs residual blocks with random Fourier features 
	to efficiently capture multiscale solution structures, and a 
	progressive training strategy with adaptive loss weighting 
	balances data fidelity and physical consistency. Comprehensive 
	numerical experiments on five benchmark problems---including 
	boundary and interior layers, traveling waves, and nonlinear 
	Burgers dynamics---demonstrate that the proposed framework 
	markedly improves solution accuracy at the terminal time compared 
	to standalone stabilized FEM solutions.
\end{abstract}

			
			
			
			
		\end{frontmatter}

		\section{Introduction}
      
		The mathematical description of physical phenomena in science and engineering 
heavily relies on (nonlinear) partial differential equations (PDEs). Convection-diffusion-reaction 
(CDR) equations constitute a particularly important class of PDEs, as they govern the 
transport and transformation of quantities such as mass, heat, and momentum. Since 
analytical solutions to such equations---especially when nonlinearities or coupled 
systems are involved---rarely possess closed-form expressions, computational techniques 
have become essential. To that end, one typically employs numerical discretization 
schemes, semi-analytical techniques, machine learning approaches, or any combination 
thereof to approximate solutions with acceptable accuracy and computational cost.

High-fidelity numerical simulation of CDR-type PDEs demands meticulous 
treatment of both spatial and temporal discretization. Classical numerical spatial 
discretization approaches, including finite difference, finite element, 
and finite volume schemes, tend to produce (node-to-node) spurious oscillations when 
convection dominates diffusion, thereby degrading the physical 
reliability of computed solutions. This long-standing difficulty has motivated the development of numerous stabilized formulations over the past several decades. Within the finite element community, the streamline-upwind/Petrov--Galerkin (SUPG)~\cite{Hughes79a, Brooks82a, Tezduyar82a, Tezduyar83a, Hughes84a}, Galerkin/least-squares (GLS)~\cite{Hughes89a}, and variational multiscale (VMS)~\cite{Hughes1995} methods have gained widespread acceptance as effective stabilization strategies. Despite their success, these techniques may still yield localized nonphysical oscillations near steep gradients, shocks, or thin boundary and internal layers. Shock-capturing (SC) (also called discontinuity-capturing) operators are therefore often appended to further suppress such oscillations in regions of rapid solution variation. A persistent challenge, however, lies in the proper determination of stabilization parameters, which are typically computed element-by-element. Although numerous empirically motivated parameter definitions exist, a universally applicable formulation remains elusive, and parameter tuning can impose a significant computational overhead, particularly for transient and high-dimensional problems. Comprehensive discussions of stabilized finite element techniques and their applications can be found in~\cite{Tezduyar91b, Franca1992, Tezduyar00a, Franca2000, Hughes2004, Burman2004, john2006, Franca2006, john2008, John2008b, Codina2008, Jiang2018, Liu2021, Jiang2023, Liu2024, Cengizci2026} and the references cited therein.

The rapid advancement of deep learning technology has spurred growing interest in its application to PDE solvers. Physics-informed neural networks (PINNs) represent a particularly promising development, wherein physical laws encoded through differential operators are embedded directly into the training loss. By minimizing residuals of the governing equations at selected collocation points, PINNs can approximate solutions without relying on labeled datasets. Nonetheless, PINNs and established numerical methods---including the standard Galerkin finite element method (GFEM)---have largely evolved along separate paths. Recent investigations, however, indicate that machine learning tools, and PINNs in particular, hold considerable promise for augmenting classical solvers by providing data-driven estimates of the artificial dissipation necessary for numerical stability. While various coupling strategies between PINNs and finite element formulations have been explored, the systematic incorporation of shock-capturing mechanisms into continuous FEM-based PINN frameworks has received very limited attention. Addressing this gap could yield improvements in both accuracy and efficiency, and extend the reach of stabilized methods to more demanding applications. 

A number of investigations have sought to leverage artificial neural networks (ANNs) for improving the stability and fidelity of conventional numerical schemes. Ray and Hesthaven~\cite{Ray2018, Ray2019}, for example, utilized multilayer perceptrons to detect oscillatory cells in high-order finite element computations and applied slope limiters locally to mitigate nonphysical behavior. Veiga et al.~\cite{Veiga2018} developed a parameter-free stabilization strategy for hyperbolic conservation laws, training ANNs to pinpoint regions where artificial dissipation is needed; their approach was validated within a discontinuous Galerkin (dGFEM) setting on scalar advection and the Euler system. Discacciati et al.~\cite{Discacciati2020} combined ANNs with Runge--Kutta dGFEM solvers to construct a universal artificial viscosity framework capable of maintaining optimal accuracy for smooth solutions while reliably capturing discontinuities. Schwander et al.~\cite{Schwander2021} trained networks to assess local regularity within a Fourier spectral context and benchmarked their method against entropy viscosity approaches on one- and two-dimensional conservation laws as well as the Euler equations.

The foundational work on PINNs by Raissi et al.~\cite{Raissi2019} demonstrated that differential equations can be solved by minimizing residuals at collocation points in an unsupervised or semi-supervised manner, without the need for labeled training data. Karniadakis et al.~\cite{Karniadakis2021} subsequently broadened this paradigm to encompass high-dimensional inverse problems, scenarios involving noisy observations, and the discovery of latent physics, underscoring the versatility of PINNs for both forward and inverse modeling. Mitusch et al.~\cite{Mitusch2021} presented a hybrid FEM--PINN implementation built on the open-source libraries \texttt{FEniCS}~\cite{logg2012automated, Alnaes2015} and \texttt{dolfin-adjoint}~\cite{Mitusch2019}, with applications to Poisson and CDR problems as well as cardiac electrophysiology. Huang et al.~\cite{Huang2025} provided a taxonomy of deep neural network (DNN) approaches for CDR-type equations along with an assessment of their applicability in scientific and engineering contexts. Jin et al.~\cite{Jin2021} introduced NSFnets, a PINN architecture tailored for the incompressible Navier--Stokes equations covering both laminar and turbulent regimes.

Yu and Hesthaven~\cite{Yu2022} constructed a data-driven artificial viscosity model for dGFEM that substantially outperformed classical modal decay indicators in smoothness detection. Yadav and Ganesan~\cite{Yadav2022} proposed SPDE-ConvNet, a convolutional neural network integrated with SUPG stabilization for singularly perturbed CDR equations, demonstrating effective parameter optimization and oscillation suppression. Huang~\cite{Huang2023} formulated a variationally consistent nodally integrated meshfree method aimed at advection-dominated problems, addressing the loss of Galerkin exactness inherent in nodal integration. Le-Duc et al.~\cite{LeDuc2023} developed the finite-element-informed neural network (FEI-NN) for parametric structural mechanics simulations. In subsequent work, Yadav and Ganesan~\cite{Yadav2024} introduced SPDE-NetII, an ANN-augmented SUPG scheme that achieved improved resolution of boundary and interior layers and surpassed the performance of standard PINNs and Variational Neural Networks (VarNet)~\cite{Khodayi19} in benchmark comparisons. Luong et al.~\cite{Luong2024} addressed the imposition of exact boundary conditions on complex geometries through a two-stage procedure: constructing DNN-based trial functions satisfying homogeneous essential conditions, followed by solution of the resulting boundary value problem.

Grossmann et al.~\cite{Grossmann2024} conducted a rigorous comparison between PINNs and FEM for linear and nonlinear PDEs, including the Poisson, Allen--Cahn, and semi-linear Schr\"{o}dinger equations, concluding that FEM generally attains superior accuracy while PINNs may offer computational advantages in certain settings. Yamazaki et al.~\cite{Yamazaki2024} integrated weak Galerkin formulations with implicit Euler time stepping to devise Finite Operator Learning (FOL), achieving high accuracy for time-dependent heat conduction in heterogeneous media. Efforts to accelerate training include the Fourier neural operators and Koopman neural networks explored by Buzaev et al.~\cite{Buzaev2023}, as well as the curriculum learning strategy of Wang et al.~\cite{YufengWang2024}, which prioritizes easier subdomains before tackling high-loss regions. Chang et al.~\cite{Chang2024a, Chang2024b} proposed singular-layer PINNs (sl-PINNs) that embed boundary layer structure directly into the network architecture, improving performance for singular perturbation problems. Frerichs-Mihov et al.~\cite{FrerichsMihov2024} designed specialized loss functions for convection-dominated flows, while Hsieh and Huang~\cite{Hsieh2024} developed a multiscale stabilized PINN formulation incorporating weak boundary enforcement, Sobolev training, and transfer learning to enhance stability near boundary layers. K\"{u}t\"{u}k and Y\"{u}cel~\cite{Ktk2025} embedded the energy dissipation structure of the Allen--Cahn equation into the PINN loss and demonstrated accurate modeling of phase separation and metastability across one to three dimensions. Most recently, Raina et al.~\cite{Raina2025} tackled two-parameter singular perturbation problems in fluid mechanics, identifying limitations of standard PINNs in resolving boundary layers and proposing finite basis PINNs (FB-PINNs) inspired by finite element methodology.

In the present work, we extend the hybrid framework proposed in~\cite{cengizci_pinn} by the authors of this current study 
---termed PASSC (\textbf{P}INN-\textbf{A}ugmented 
\textbf{S}UPG with \textbf{S}hock-\textbf{C}apturing)---originally 
developed for advection-dominated steady CDR-type PDEs, to the 
unsteady regime. The main contributions of this study can be summarized as follows. 
First, we present a hybrid training paradigm that leverages two 
sophisticated FEM ingredients---SUPG stabilization and the YZ$\beta$ 
shock-capturing mechanism~\cite{Tezduyar02a, Tezduyar04o}---to produce 
high-fidelity training data for unsteady convection-dominated problems 
exhibiting boundary and/or internal layers. Although SUPG stabilization 
performs well in many such scenarios, auxiliary shock-capturing is often 
necessary to adequately resolve sharp gradients and internal 
discontinuities~\cite{cengizci22, Cengizci2024_3d}. Accordingly, we 
employ SUPG-YZ$\beta$ solutions as training data while selectively 
imposing PDE residual constraints in stable interior regions. Second, 
the proposed architecture is dimension-agnostic, featuring deep residual 
networks with Fourier feature embeddings, layer normalization, and 
conservative hyperparameter choices. Third, the methodology is inherently 
scalable, accommodating problems in one, two, and three spatial dimensions.

    \section{Mathematical Framework} \label{sec:section2}
Let $n_\text{sd} \geq 1$ denote the spatial dimension, and let $\Omega \subset \mathbb{R}^{n_\text{sd}}$ be a bounded Lipschitz domain representing the physical region occupied by the flow, i.e., the spatial domain in which the problem is defined. To establish the mathematical framework addressed in this manuscript, we consider the following time-dependent CDR-type PDE as a representative model problem:
\begin{alignat}{2}
\frac{\partial u(t, \mathbf{x})}{\partial t} - \varepsilon \Delta u(t, \mathbf{x}) + \mathbf{b}(t, \mathbf{x})\cdot\nabla u(t, \mathbf{x}) + c(\mathbf{x}) u(t, \mathbf{x}) &= f(t, \mathbf{x}) &\; &\text{in } \; \Omega \times (t_0, t_\text{f}], \label{eq:CDR1}\\
u(t, \mathbf{x}) &= u_D(t, \mathbf{x}) &\; &\text{on } \; \partial\Omega \times (t_0, t_\text{f}], \label{eq:CDR2} \\
u(t_{0}, \mathbf{x}) &= u_0(\mathbf{x}) &\; &\text{in } \; \Omega. \label{eq:CDR3}
\end{alignat}
Here, $u$ denotes the unknown scalar field, where $t \in (t_0, t_\text{f}]$ is the temporal variable and $\mathbf{x} = (x_1, x_2, \ldots, x_{n_\text{sd}}) \in \Omega$ is the spatial position vector. The positive constant $\varepsilon$ represents the diffusion coefficient, which is assumed to satisfy $0 < \varepsilon \ll 1$ in convection-dominated regimes. The vector field
$\mathbf{b} \in L^{\infty}([t_0, t_{\mathrm{f}}];\,
\mathcal{W}^{1,\infty}(\Omega)^{n_{\mathrm{sd}}})$
denotes the convection velocity, $c \in L^{\infty}(\Omega)$ is the
reaction coefficient, and
$f \in L^2([t_0, t_{\mathrm{f}}];\, L^2(\Omega))$ is the source term. The function $u_0 \in L^2(\Omega)$ prescribes the initial 
condition, while $u_D \in C([t_0, t_\text{f}]; \mathcal{H}^{1/2}(\partial\Omega))$ 
specifies the time-dependent Dirichlet boundary data, where 
$\mathcal{H}^{1/2}(\partial\Omega)$ denotes the trace space of 
$\mathcal{H}^1(\Omega)$ on $\partial\Omega$, such that for each 
$t \in [t_0, t_\text{f}]$, the function $u_D(t, \cdot)$ admits an extension 
to $\mathcal{H}^1(\Omega)$. Throughout this work, we adopt standard Sobolev space notation: $L^2(\Omega)$ and $L^{\infty}(\Omega)$ denote the spaces of square-integrable and essentially bounded functions, respectively; $\mathcal{W}^{1,\infty}(\Omega)$ is the space of functions with essentially bounded first-order weak derivatives; and $\mathcal{H}^1(\Omega)$ denotes the space of functions  with square-integrable first-order derivatives. For a Banach 
space $X$, the notation $C([t_0, t_\text{f}]; X)$ and $L^2([t_0, t_\text{f}]; X)$ 
refers to the standard spaces of continuous and square-integrable 
$X$-valued functions, respectively.

\begin{remark} \label{rem:Neumann}
Although the governing equations~\eqref{eq:CDR1}--\eqref{eq:CDR3} are 
equipped with Dirichlet boundary conditions, the proposed framework 
readily accommodates Neumann-type conditions as well. For brevity and 
clarity of exposition, we focus exclusively on the Dirichlet case 
throughout the subsequent development.
\end{remark}

\begin{remark} \label{rem:notation}
For notational convenience in the finite element formulations of 
Sections~\ref{subsection:gfem}--\ref{subsection:yzb}, the explicit 
spatial dependence of field quantities is omitted whenever the 
meaning remains clear from context. Accordingly, the convection 
velocity $\mathbf{b}(t, \mathbf{x})$, reaction coefficient $c(\mathbf{x})$, 
and source term $f(t, \mathbf{x})$ are henceforth denoted by 
$\mathbf{b}$, $c$, and $f$, respectively.
\end{remark}

\subsection{Standard Galerkin finite element method formulation} \label{subsection:gfem}
The standard Galerkin finite element formulation corresponding to the model problem defined in Eqs.~\eqref{eq:CDR1}--\eqref{eq:CDR3} seeks a solution $u^h \in \mathcal{S}_u^h$ such that for all test functions $w^h \in \mathcal{V}_u^h$:
\begin{equation}
    \int_\Omega w^h \frac{\partial u^h}{\partial t} \, d\Omega + \int_\Omega \left( \varepsilon \nabla w^h \cdot \nabla u^h + w^h \, \mathbf{b} \cdot \nabla u^h + c \, w^h u^h \right) d\Omega = \int_\Omega w^h f \, d\Omega.
    \label{eq:gfem}
\end{equation}
In this formulation, the superscript $h$ denotes membership in a  finite-dimensional function space. The discrete trial and test  function spaces are defined as:
\begin{align}
\mathcal{S}^{h}_{u} &= \left\{ u^{h} \in \mathcal{H}^{1h}(\Omega) : u^{h} \doteq u_D^{h} 
\text{ on } \partial\Omega \right\}, \\
\mathcal{V}^{h}_{u} &= \left\{ w^{h} \in \mathcal{H}^{1h}(\Omega) : w^{h} \doteq 0 
\text{ on } \partial\Omega \right\},
\end{align}
where $u_D^{h}$ denotes the discrete representation of the boundary data $u_D$ 
on the finite element space, and the symbol $\doteq$ indicates that the boundary 
values are assigned in the sense of the finite element interpolation. The space $\mathcal{H}^{1h}(\Omega)$ denotes a finite-dimensional subspace of the Sobolev space $\mathcal{H}^1(\Omega)$ and is defined by:
\begin{equation}
\mathcal{H}^{1h}(\Omega) = \left\{ \phi^h \in \mathcal{C}^0(\bar{\Omega}) : \left. \phi^h \right|_{\Omega^e} \in \mathbb{P}_1(\Omega^e), \; \forall \Omega^e \in \mathcal{T}^h \right\}.
\end{equation}
Here, $\mathcal{C}^0(\bar{\Omega})$ represents the space of continuous functions defined over the closure of the computational domain $\Omega$. The set $\mathbb{P}_1(\Omega^e)$ denotes the space of piecewise linear polynomial functions defined on each finite element $\Omega^e$, and $\mathcal{T}^h$ represents the mesh partition of the domain into finite elements.
    
    This standard Galerkin formulation serves as a baseline for the numerical treatment of CDR-type PDEs. However, in convection-dominated regimes, it is well known that this approach typically suffers from numerical instabilities and nonphysical oscillations, motivating the need for stabilization techniques such as SUPG and shock-capturing methods discussed in subsequent sections.

    \subsection{SUPG-stabilized finite element formulation} \label{subsection:supg}
    The SUPG-stabilized finite element formulation corresponding to the model problem in Eqs.~\eqref{eq:CDR1}--\eqref{eq:CDR3} seeks a discrete solution $u^h \in \mathcal{S}_u^h$ such that, for all test functions $w^h \in \mathcal{V}_u^h$:
\begin{multline}
\int_\Omega w^h \frac{\partial u^h}{\partial t} \, d\Omega 
+ \int_\Omega \left( 
\varepsilon \nabla w^h \cdot \nabla u^h 
+ w^h \, \mathbf{b} \cdot \nabla u^h 
+ c \, w^h u^h 
\right) d\Omega 
- \int_\Omega w^h f \, d\Omega \\
+ \sum_{e=1}^{n_{\text{el}}} \int_{\Omega^e} 
\tau_{\text{SUPG}} \left( \mathbf{b} \cdot \nabla w^h \right)
\underbrace{\left( 
\frac{\partial u^h}{\partial t}
- \varepsilon \Delta u^h 
+ \mathbf{b} \cdot \nabla u^h 
+ c u^h - f 
\right)}_{\mathcal{R}(u^h)}
d\Omega^e = 0.
\label{eq:supg}
\end{multline}
The first integral in formulation~\eqref{eq:supg} corresponds to the standard GFEM formulation, while the second line represents the SUPG contribution. The solution and test function spaces $\mathcal{S}_u^h$ and $\mathcal{V}_u^h$ are identical to those defined in the GFEM formulation. The index $n_{\text{el}}$ denotes the number of finite elements in the mesh discretization, and $\tau_{\text{SUPG}}$ is the stabilization parameter associated with the SUPG formulation.

\begin{remark} \label{rem:polynomial_degree}
Throughout this study, all finite element computations employ the 
continuous piecewise-linear space $\mathcal{H}^{1h}(\Omega)$. Although 
higher-order polynomial bases (e.g., $\mathbb{P}_2$ or $\mathbb{P}_3$) 
may yield superior accuracy in smooth regions, they incur greater 
computational costs and can become susceptible to numerical instabilities 
when convection dominates, especially in the presence of sharp internal 
layers or steep gradients. Moreover, without adequate stabilization, 
elevated polynomial degrees tend to amplify spurious oscillations in 
such regimes. In contrast, linear elements offer reduced assembly and 
solution times while integrating seamlessly with SUPG stabilization and 
shock-capturing techniques. Consequently, we deliberately employ 
stabilized $\mathbb{P}_1$ elements to achieve an effective compromise 
between computational efficiency and robustness, thereby facilitating 
accurate simulations across a wide range of convection-dominated problems.
\end{remark}

In formulation~\eqref{eq:supg}, the parameter $\tau_{\text{SUPG}}$ dictates the level of artificial 
diffusion added to the discrete system for stabilization purposes. 
Selecting this parameter involves a delicate balance: excessively 
large values produce over-diffused solutions that compromise accuracy, 
while values that are too small permit oscillatory artifacts to persist. 
Thus, careful, problem-dependent tuning of $\tau_{\text{SUPG}}$ is 
indispensable. Although the literature offers a variety of robust 
definitions, no single universally optimal formula has emerged. In this present study, we prefer to adopt the stabilization parameter 
proposed by Shakib et al.~\cite{Shakib88}, with additional discussion 
available in~\cite{Tezduyar91c, Franca1992}:
\begin{equation} \label{eq:shakib_tau2}
\tau_{\text{SUPG}} = \left( \left(\frac{2}{\varDelta t} \right)^2 
+ \left(\frac{2\|\mathbf{b}\|}{h_K}\right)^2 
+ \left(\frac{4\varepsilon}{h_K^2}\right)^2 \right)^{-1/2}.
\end{equation}
Here, $\|\mathbf{b}\|$ denotes the magnitude of the convection velocity 
(with $\|\cdot\|$ representing the Euclidean norm used consistently 
throughout this paper), $\varDelta t$ represents the time-step size, and $h_K$ is the characteristic mesh size 
associated with element $K = \Omega^e$, defined as
\begin{equation}
h_K = \text{diam}(K), \quad \forall K \in \mathcal{T}^h.
\label{eq:dia}
\end{equation}

\begin{remark} \label{rem:consistency_supg}
A key property of the SUPG formulation in Eq.~\eqref{eq:supg} is that 
the stabilization term incorporates the residual, $\mathcal{R}(u^h)$, of the governing  PDE~\eqref{eq:CDR1} as a weighting factor. This construction ensures  consistency: whenever the exact solution is substituted into the  formulation, the stabilization contribution vanishes identically.
\end{remark} 

Extensive numerical evidence~\cite{john2006, john2008, John2008b, cengizci22, 
Cengizci2023_tr, Cengizci2023amc, Cengizci2024_3d, Cengizci2023zamm, 
Cengizci2024_mhd, Cengizci2024_heston, Cengizci2026} indicates that the 
SUPG method, while effective in many scenarios, can be inadequate for 
suppressing node-to-node oscillations near sharp gradients or shock-like 
features. To address this shortcoming, supplementary shock-capturing 
mechanisms are typically required. In the present work, we adopt the 
YZ$\beta$ shock-capturing technique, originally proposed by Tezduyar 
and collaborators~\cite{Tezduyar02a, Tezduyar02j, Tezduyar04o, Tezduyar05d, 
Tezduyar05f}, to enhance solution quality in regions exhibiting steep 
spatial variations.

\subsection{SUPG-YZ$\beta$ formulation} \label{subsection:yzb}
The GFEM finite element formulation, augmented with SUPG 
stabilization and YZ$\beta$ shock-capturing, reads: find 
$u^h \in \mathcal{S}_u^h$ such that for all 
$w^h \in \mathcal{V}_u^h$,
\begin{multline}
	\int_\Omega w^h \frac{\partial u^h}{\partial t} \, d\Omega
	+ \int_\Omega \Big(
	\varepsilon \nabla w^h \cdot \nabla u^h 
	+ \mathbf{b} \cdot \nabla u^h \, w^h 
	+ c\, u^h w^h 
	- f w^h 
	\Big) \, d\Omega \\
	+ \sum_{e=1}^{n_{\text{el}}} \int_{\Omega^e} 
	\tau_{\text{SUPG}} (\mathbf{b} \cdot \nabla w^h) 
	\Big( 
	\frac{\partial u^h}{\partial t}
	- \varepsilon \Delta u^h 
	+ \mathbf{b} \cdot \nabla u^h 
	+ c\, u^h - f 
	\Big) \, d\Omega \\
	+ \sum_{e=1}^{n_{\text{el}}} \int_{\Omega^e} 
	\nu_{\text{SHOC}} (\nabla w^h \cdot \nabla u^h) \, d\Omega 
	= 0.
	\label{eq:supg_shock}
\end{multline}
The first line in Eq.~\eqref{eq:supg_shock} represents the standard Galerkin formulation, the second introduces SUPG stabilization terms, and the third includes the shock-capturing term. The function spaces $\mathcal{S}_u^h$ and $\mathcal{V}_u^h$ are defined as previously. The YZ$\beta$ shock-capturing parameter, $\nu_{\text{SHOC}}$, can be defined as follows~\cite{Tezduyar02a, Tezduyar02j, Tezduyar04o, Tezduyar05d, Tezduyar05f}:
\begin{equation}
\nu_{\text{SHOC}} = \left| \text{Y}^{-1} \text{Z} \right|
\left( \sum_{i=1}^{n_\text{sd}} \left| \text{Y}^{-1} \frac{\partial u^h}{\partial x_i} \right|^2 \right)^{\frac{\beta}{2} - 1}
\left( \frac{h_{\text{SHOC}}}{2} \right)^\beta,
\label{eq:nu_shoc}
\end{equation}
where $\text{Y}$ is a scalar representing a reference value of the solution $u$.

\begin{remark} \label{rem:remarkY}
The parameter $\mathrm{Y}$ is a nonzero scaling constant selected to 
reflect a reference scale of the solution. A typical specification is
\begin{equation}
    \mathrm{Y} = \frac{u^{\max} + u^{\min}}{2},
\end{equation}
representing the mean of the expected solution range; more generally, 
any appropriate value between $u^{\min}$ and $u^{\max}$ may be chosen. 
For transient problems, the initial condition often provides a 
suitable basis for determining $\mathrm{Y}$. In the context of coupled 
systems of equations, $\mathrm{Y}$ is replaced by a diagonal matrix 
$\mathbf{Y}$ whose entries encode the reference scales for each 
governing equation; see~\cite{Tezduyar05f, Cengizci2023_tr} for 
further discussion.
\end{remark}

The quantity $\mathrm{Z}$ appearing in Eq.~\eqref{eq:nu_shoc} represents 
the pointwise residual of the time-dependent governing equation, expressed 
following~\cite{bazilevs_yzb_2007} as
\begin{equation}
\mathrm{Z} = \mathcal{R}(u^h) = \frac{\partial u^h}{\partial t} - \varepsilon \Delta u^h 
+ \mathbf{b} \cdot \nabla u^h + c u^h - f.
\label{eq:shock_residual}
\end{equation}

\begin{remark}

This residual definition extends the original formulation of 
Tezduyar~\cite{Tezduyar02a, Tezduyar02j, Tezduyar04o, Tezduyar05d, 
Tezduyar05f} by incorporating all terms from the governing equation. 
As a consequence, the shock-capturing contribution vanishes identically 
when the exact solution is attained, thereby avoiding spurious numerical 
diffusion and ensuring consistency of the formulation. For systems of 
coupled equations, the scalar $\mathrm{Z}$ generalizes to a residual 
vector $\mathbf{Z}$, and in view of Remark~\ref{rem:remarkY}, the 
absolute value $\left| \mathrm{Y}^{-1} \mathrm{Z} \right|$ is replaced 
by the norm $\Vert \mathbf{Y}^{-1} \mathbf{Z} \Vert$.
\end{remark}

\begin{remark} \label{rem:sharpness_yzb}
The exponent $\beta$ in Eq.~\eqref{eq:nu_shoc} controls the sensitivity 
of the shock-capturing operator to solution gradients. In particular, 
$\beta = 1$ is appropriate for moderately sharp transitions, while 
$\beta = 2$ provides stronger regularization for severely steep gradients. 
Intermediate values may also be employed, as discussed in~\cite{Tezduyar05f}. 
Given that the present work targets strongly convection-dominated problems 
exhibiting pronounced solution fronts, we adopt $\beta = 2$ throughout.
\end{remark}

The shock-capturing length scale, $h_{\text{SHOC}}$, is defined as:
\begin{equation}
h_{\text{SHOC}} = \left[ 2 \left( \sum_{a=1}^{n_{\text{en}}} \left| \mathbf{j} \cdot \nabla N_a \right| \right) \right]^{-1}, 
\quad \text{with} \quad \mathbf{j} = \frac{\nabla u^h}{\|\nabla u^h\|},
\label{eq:h_shoc}
\end{equation}
where $n_{\text{en}}$ is the number of element nodes, and $N_a$ is the shape function corresponding to node $a$.

\begin{remark} \label{rem:yzb_lenght}
Many shock-capturing formulations define the length scale 
$h_{\text{SHOC}}$ using directionally weighted measures aligned 
with the local convection field. Despite their potential for improved 
accuracy, such approaches often require additional mesh-related 
preprocessing and can complicate implementation. In the present work, 
we adopt a simpler strategy by setting $h_{\text{SHOC}} = h_K$ 
(see~ Eq.~\eqref{eq:dia}), which avoids directional dependence and 
reduces computational cost. While this approximation may be less 
precise on strongly anisotropic element configurations, the numerical 
results presented herein demonstrate that it provides stable and 
sufficiently accurate predictions across all test cases considered.
\end{remark}

\begin{remark} \label{rem:Z_yzb}
Based on the considerations in Remarks~\ref{rem:sharpness_yzb} 
and~\ref{rem:yzb_lenght}, the shock-capturing viscosity 
$\nu_{\text{SHOC}}$ in Eq.~\eqref{eq:nu_shoc} simplifies to
\begin{equation} \label{eq:nu_shocb2}
    \nu_{\text{SHOC}} = 
    \frac{\left| \mathrm{Y}^{-1} \mathrm{Z} \right| h_K^2}{4}.
\end{equation}
\end{remark}

As with the SUPG stabilization parameter, selecting an appropriate 
value for $\nu_{\text{SHOC}}$ remains a challenging and inherently 
problem-specific task. In practice, this coefficient is often 
determined through empirical calibration or heuristic tuning. This 
difficulty highlights the appeal of data-driven methodologies such 
as PINNs, which can learn and adjust stabilization parameters 
directly from physical constraints and available data. Beyond parameter 
estimation, PINNs can serve as post-processing correctors that refine 
stabilized finite element solutions by minimizing residuals derived 
from the governing equations and boundary/initial conditions. This 
capability is examined in detail in the following section.

\begin{remark} \label{whysuphandyzb}
The finite element component of the proposed hybrid methodology is not limited to the particular stabilization and shock-capturing schemes employed here. Other  stabilized formulations, including GLS or VMS, can be seamlessly  integrated into the framework. Likewise, alternative shock-capturing  operators or spurious oscillations at layers diminishing (SOLD)  methods may replace the YZ$\beta$ technique adopted in this work. 
The choice of SUPG stabilization combined with YZ$\beta$ shock-capturing 
is motivated primarily by the authors' prior experience with these 
approaches in the context of convection-dominated transport problems.
\end{remark}

\section{Integration of PINNs with SUPG-YZ$\beta$ Formulation}
Following the framework introduced for steady-state problems in~\cite{cengizci_pinn}, this 
section presents the hybrid computational strategy that couples 
stabilized finite element approximations with physics-informed neural 
networks. The network architecture, 
together with the multi-phase training procedure, is depicted in 
Figures~\ref{fig:selective_physics_flowchart} 
and~\ref{fig:nn-architecture}.

\begin{remark} \label{rem:dimension}
Although the proposed framework applies to CDR-type problems in one, 
two, or three spatial dimensions, all formulations and algorithmic 
descriptions are presented for the two-dimensional setting 
($n_{\text{sd}}=2$) unless stated otherwise. This convention is 
adopted for notational simplicity and readability; the generalization 
to other dimensions is straightforward.
\end{remark}

\subsection{Motivation for Hybrid PINN-FEM Integration}

Coupling PINNs with the SUPG-YZ$\beta$ formulation mitigates key 
shortcomings that arise when either method is applied standalone. 
The SUPG-YZ$\beta$ approach, though numerically robust, depends on 
user-specified parameters whose optimal values are generally determined 
through empirical rules or case-by-case calibration, potentially 
limiting performance across varying flow scenarios. Conversely, 
conventional PINNs, despite their conceptual elegance, frequently 
encounter difficulties in resolving thin boundary layers and may 
suffer from prolonged training times in strongly convection-dominated 
settings. The proposed hybrid strategy overcomes these obstacles 
through several innovations. Numerically stable SUPG-YZ$\beta$ 
solutions provide informative priors, allowing the neural network 
to concentrate on resolving fine-scale features rather than learning 
fundamental PDE behavior from the outset. Simultaneously, the 
framework exploits the stability of SUPG-YZ$\beta$ discretizations 
to accelerate convergence while retaining the representational 
flexibility of PINNs for complex physical phenomena. This constitutes 
the principal contribution of the present study relative to existing 
FEM-PINN coupling strategies (see, e.g.,~\cite{Mitusch2021, LeDuc2023, 
Yadav2024}): by anchoring network training in stabilized finite element 
solutions, the methodology inherits the robustness of classical 
numerical schemes while benefiting from the adaptability of 
data-driven learning.

\subsection{Network Architecture and Fourier Feature Mapping} 
\label{sec:network_arch}

The PINN architecture employs random Fourier feature 
mapping to enhance the network's ability to 
capture high-frequency variations characteristic of boundary layers 
and internal shocks. Given spatiotemporal input coordinates 
$\mathbf{z} = (t, \mathbf{x}) \in \mathbb{R}^{n_\text{sd}+1}$, the Fourier feature embedding is defined as:
\begin{equation}
\boldsymbol{\phi}(\mathbf{z}) = \begin{bmatrix} \mathbf{z} \\ 
\sin(\mathbf{B}^\top \mathbf{z}) \\ 
\cos(\mathbf{B}^\top \mathbf{z}) \end{bmatrix} 
\in \mathbb{R}^{n_\text{sd}+1 + 2 n_\text{F}},
\label{eq:fourier_features}
\end{equation}
where $\mathbf{B} \in \mathbb{R}^{(n_\text{sd}+1) \times n_\text{F}}$ is a random 
frequency matrix whose entries are independently sampled from a 
normal distribution $\mathcal{N}(0, \sigma^2)$ and held fixed 
(i.e., not trained) throughout optimization. Here, $n_\text{F}$ denotes 
the number of Fourier features and $\sigma$ controls the frequency 
bandwidth. The trigonometric functions are applied element-wise to 
the projection $\mathbf{B}^\top \mathbf{z} \in \mathbb{R}^{n_\text{F}}$, 
and the resulting feature vector $\boldsymbol{\phi}(\mathbf{z})$ 
concatenates the raw coordinates with the sine and cosine 
projections.

The Fourier-enhanced features are first mapped to a hidden 
dimension $n_h$ through an input layer:
\begin{equation}
\mathbf{h}_0 = \texttt{SiLU}\!\left(\mathbf{W}^{(0)} 
\boldsymbol{\phi}(\mathbf{z}) + \tilde{\mathbf{b}}^{(0)}\right),
\label{eq:input_mapping}
\end{equation}
where $\mathbf{W}^{(0)} \in \mathbb{R}^{n_h \times (n_\text{sd}+1+2n_\text{F})}$ and 
$\tilde{\mathbf{b}}^{(0)} \in \mathbb{R}^{n_h}$. The network then 
processes this representation through $n_r$ residual blocks, each 
of which transforms the hidden state according to:
\begin{align}
\mathbf{p}_i &= \texttt{SiLU}\!\left(\texttt{Linear}_1
(\mathbf{h}_i)\right), \label{eq:residual_step1} \\
\mathbf{q}_i &= \texttt{Linear}_2(\mathbf{p}_i), 
\label{eq:residual_step2} \\
\mathbf{h}_{i+1} &= \texttt{SiLU}\!\left(\texttt{LayerNorm}
(\mathbf{q}_i + \mathbf{h}_i)\right), 
\label{eq:residual_block}
\end{align}
where $\mathbf{h}_i \in \mathbb{R}^{n_h}$ denotes the hidden 
state at the $i$-th block for $i = 0, 1, \ldots, n_r - 1$. Here, 
$\texttt{Linear}_1$ and $\texttt{Linear}_2$ denote standard fully 
connected (affine) transformations of the form
\begin{equation}
\texttt{Linear}(\mathbf{v}) = \mathbf{W}\mathbf{v} + \tilde{\mathbf{b}},
\end{equation}
with trainable weight $\mathbf{W} \in \mathbb{R}^{n_h \times n_h}$ 
and bias $\tilde{\mathbf{b}} \in \mathbb{R}^{n_h}$; each residual block 
uses two such layers with distinct parameter sets. The skip 
connection $\mathbf{q}_i + \mathbf{h}_i$ followed by layer 
normalization mitigates gradient degradation in deep architectures. 
The specific hyperparameter values ($n_h$, $n_r$, $n_\text{F}$, $\sigma$) 
adopted for each numerical experiment are reported in 
Section~\ref{sec:section5}.

The \texttt{SiLU} (sigmoid-weighted linear unit, Swish) activation function is defined as:
\begin{equation}
\texttt{SiLU}(x) = x \cdot \texttt{sigmoid}(x) = \frac{x}{1 + e^{-x}}.
\label{eq:silu}
\end{equation}
This smooth, non-monotonic activation promotes better gradient flow compared to \texttt{ReLU} (rectified linear unit)-based alternatives and has demonstrated superior performance in physics-informed learning tasks, particularly for problems involving smooth solution variations. \texttt{LayerNorm} provides training stability through normalization:
\begin{equation}
\texttt{LayerNorm}(\mathbf{x}) = \frac{\mathbf{x} - \mu_\text{LN}}{\sqrt{\sigma_\text{LN}^2 + \epsilon_\text{LN}}} \odot \boldsymbol{\gamma_\text{LN}} + \boldsymbol{\beta_\text{LN}},
\label{eq:layernorm}
\end{equation}
where $\mu_\text{LN}$ and $\sigma^2_\text{LN}$ are the mean and variance computed over the feature dimension, $\boldsymbol{\gamma_\text{LN}}$ and $\boldsymbol{\beta_\text{LN}}$ are learnable scale and shift parameters, respectively, $\epsilon_\text{LN}$ is a small constant for numerical stability, and $\odot$ denotes the element-wise (Hadamard) product.

\subsection{Data-driven Training Approach}
Using the numerically stable FEM approximations obtained by employing 
the SUPG-YZ$\beta$ formulation, we train the PINN model with the 
mean squared error data loss function:
\begin{equation}
	\mathcal{L}_{\text{data}}(\theta) = \frac{1}{K_s} \sum_{k=1}^{K_s} 
	\left[ \frac{1}{N_h} \sum_{i=1}^{N_h} \left\vert 
	u_{\text{NN}}(t_k, \mathbf{x}_i; \theta) 
	- u^{h}_{\text{SUPG-YZ}\beta}(t_k, \mathbf{x}_i) 
	\right\vert^2 \right],
	\label{eq:data_loss}
\end{equation}
where $K_s$ denotes the number of temporal snapshots selected near the 
terminal time and $N_h$ is the number of nodal vertices in the finite 
element mesh $\mathcal{T}^h$, yielding a total of $N = K_s \times N_h$ 
space--time training pairs. For each snapshot $t_k$, the inner sum 
computes the spatial mean squared error over all mesh nodes 
$\{\mathbf{x}_i\}_{i=1}^{N_h}$, and the outer sum averages these 
snapshot-wise errors over the $K_s$ selected time levels. The 
spatiotemporal coordinates $(t_k, \mathbf{x}_i)$ correspond to the 
FEM solution evaluated at discrete time levels $t_k$ and mesh nodes 
$\mathbf{x}_i$. The term 
$u_{\text{NN}}(t_k, \mathbf{x}_i; \theta)$ represents the neural 
network prediction, $u^{h}_{\text{SUPG-YZ}\beta}(t_k, \mathbf{x}_i)$ 
is the corresponding stabilized FEM solution, and $\theta$ denotes 
the trainable network parameters. This loss function enforces the 
neural network to learn the spatiotemporal evolution and boundary 
layer structures captured by the stabilized finite element method. 
The data-driven approach offers several advantages: (i)~it leverages 
the computational expertise embedded in mature FEM solvers, (ii)~it 
avoids the challenges associated with balancing multiple loss terms 
in purely physics-informed training, and (iii)~it ensures that the 
learned solution inherits the stability properties of the underlying 
SUPG-YZ$\beta$ formulation.
\begin{remark} \label{rem:selective}
	The selective enforcement strategy detailed in 
	Section~\ref{sec:Selective} applies exclusively to the PDE residual 
	loss $\mathcal{L}_{\mathrm{pde}}(\theta)$, which is computed only at 
	interior collocation points. Both the data loss 
	$\mathcal{L}_{\mathrm{data}}(\theta)$ and the boundary loss 
	$\mathcal{L}_{\mathrm{bc}}(\theta)$ are computed over their complete 
	respective point sets without any distance-based or selective 
	filtering: $\mathcal{L}_{\mathrm{data}}$ utilizes all 
	$N = K_s \times N_h$ space--time points to maintain global solution 
	fidelity, while $\mathcal{L}_{\mathrm{bc}}$ employs all boundary 
	nodes on $\partial \Omega$.
\end{remark}

\subsection{Hybrid Physics-informed Training with Selective Enforcement}

To incorporate physical constraints while maintaining data fidelity, we propose a multiphase hybrid training strategy with selective physics enforcement. This approach recognizes that not all regions of the computational domain are equally suitable for physics-based learning. The hybrid loss function combines weighted components, for example:
\begin{equation}
\mathcal{L}_{\text{hybrid}}(\theta) = w_{\text{data}} \mathcal{L}_{\text{data}}(\theta) + w_{\text{pde}} \mathcal{L}_{\text{pde}}(\theta) + w_{\text{bc}} \mathcal{L}_{\text{bc}}(\theta),
\label{eq:hybrid_loss}
\end{equation}
where each component is to be defined below, and the nonnegative scalar weights $w_{\text{data}}, w_{\text{pde}}, w_{\text{bc}}$ control the relative importance of data consistency, PDE residual minimization, and boundary condition enforcement, respectively.

\begin{remark} \label{rem:curriculum_learn}
    These weights serve as balancing factors to tune the contribution of each term to the overall training objective. The dynamic evolution of these weights throughout training enables a curriculum learning approach, where the network first learns to approximate the reference solution before gradually incorporating physical constraints. This prevents early training instabilities that can arise when physics terms dominate before a basic solution approximation is established.
\end{remark}

The PDE residual loss enforces the governing CDR-type equation 
over randomly sampled interior collocation points:
\begin{equation}
\mathcal{L}_{\mathrm{pde}}(\theta) = \frac{1}{N_{\mathrm{int}}} 
\sum_{j=1}^{N_{\mathrm{int}}} \left\vert 
\mathcal{R^\text{NN}}(t_j, \mathbf{x}_j; \theta) \right\vert^2,
\label{eq:pde_loss}
\end{equation}
where the interior collocation set is defined as
\begin{equation}
\mathcal{X}_{\mathrm{int}} = \left\{ (t_j, \mathbf{x}_j) \in 
[t_{\mathrm{start}}, t_{\mathrm{end}}] \times \Omega \;:\; 
d_{\partial\Omega}(\mathbf{x}_j) > d_{\min} \right\}.
\label{eq:interior_set}
\end{equation}
Here, the spatiotemporal collocation points $(t_j, \mathbf{x}_j)$ 
are drawn uniformly at random within the training time window 
$[t_{\mathrm{start}}, t_{\mathrm{end}}]$ and the spatial domain 
$\Omega$, and $d_{\partial\Omega}(\mathbf{x})$ denotes the 
Euclidean distance from point $\mathbf{x}$ to the nearest boundary 
$\partial\Omega$. The threshold parameter $d_{\min}$ controls the 
exclusion zone near boundaries where PDE residuals are not enforced; 
its value is problem-dependent and should be chosen based on the 
expected boundary layer thickness and mesh resolution 
(see Section~\ref{sec:distance-based} for details). The term 
$N_{\mathrm{int}}$ denotes the number of interior collocation points 
that satisfy the distance criterion. The strong-form PDE residual 
$\mathcal{R^\text{NN}}$ for the PINNs is defined as:
\begin{equation}
\mathcal{R^\text{NN}}(t, \mathbf{x}; \theta) = 
\frac{\partial u_{\mathrm{NN}}}{\partial t} 
- \varepsilon \, \Delta u_{\mathrm{NN}} 
+ \mathbf{b}(t, \mathbf{x}) \cdot \nabla u_{\mathrm{NN}} 
+ c(\mathbf{x}) \, u_{\mathrm{NN}} 
- f(t, \mathbf{x}),
\label{eq:pde_residual}
\end{equation}
where $u_{\mathrm{NN}} = u_{\mathrm{NN}}(t, \mathbf{x}; \theta)$ 
and all partial derivatives are computed via automatic 
differentiation with respect to the spatiotemporal input coordinates. 
The temporal derivative $\partial u_{\mathrm{NN}} / \partial t$ is 
obtained entirely through automatic differentiation, eliminating the 
need for finite-difference approximations or access to the FEM 
solution at adjacent time steps.

The boundary condition loss enforces the Dirichlet condition 
$u = u_D$ on $\partial\Omega$ by penalizing deviations of the 
neural network prediction from the prescribed boundary values:
\begin{equation}
\mathcal{L}_{\mathrm{bc}}(\theta) = \frac{1}{N_{\mathrm{bc}}} 
\sum_{k=1}^{N_{\mathrm{bc}}} \left\vert 
u_{\mathrm{NN}}(t_k, \mathbf{x}_k; \theta) 
- u_D(t_k, \mathbf{x}_k) \right\vert^2,
\label{eq:bc_loss}
\end{equation}
where $\mathbf{x}_k \in \partial\Omega$ are collocation points 
distributed uniformly along each edge of the domain boundary, 
$t_k \in [t_{\mathrm{start}}, t_{\mathrm{end}}]$ are sampled 
at discrete time levels within the training window, and 
$u_D(t_k, \mathbf{x}_k)$ denotes the prescribed Dirichlet data, 
which may in general depend on both space and time. The total 
number of boundary collocation points is 
$N_{\mathrm{bc}} = N_{\mathrm{sides}} \times 
N_{\mathrm{pts}} \times N_{\mathrm{times}}$, where 
$N_{\mathrm{sides}}$ is the number of boundary edges, 
$N_{\mathrm{pts}}$ is the number of points per edge, and 
$N_{\mathrm{times}}$ is the number of sampled time levels.

\begin{remark}\label{rem:Ntimes_vs_Ks}
	The number of temporal levels used for the boundary condition loss,
	$N_{\mathrm{times}}$, need not coincide with the number of snapshots
	$K_s$ employed in the data loss~$\mathcal{L}_{\mathrm{data}}(\theta)$.
	In the data loss (Eq.~\eqref{eq:data_loss}), $K_s$ denotes the number
	of FEM solution snapshots selected near the terminal time, each
	contributing $N_h$ spatial training points. In the boundary loss
	(Eq.~\eqref{eq:bc_loss}), $N_{\mathrm{times}}$ controls the temporal
	sampling density along the boundary and may be chosen independently
	of~$K_s$ to balance boundary fidelity against computational cost.
	In particular, since boundary data are typically smooth or piecewise
	constant in time, a smaller value of $N_{\mathrm{times}}$ often
	suffices.
\end{remark}

\begin{remark} \label{remark:lift}
When the network architecture incorporates a lift function 
$u_{\mathrm{lift}}$ and a distance function 
$d(\mathbf{x}) = x_1(1-x_1)\,x_2(1-x_2)$ such that the 
output takes the form
\begin{equation}
u_{\mathrm{NN}}(t, \mathbf{x}; \theta) = u_{\mathrm{lift}}
(t, \mathbf{x}) + d(\mathbf{x})\, \mathcal{N}(t, \mathbf{x}; 
\theta),
\label{eq:lift_architecture}
\end{equation}
the Dirichlet conditions are satisfied exactly by construction, 
since $d(\mathbf{x}) = 0$ on $\partial\Omega$ implies 
$u_{\mathrm{NN}} = u_{\mathrm{lift}} = u_D$ on the boundary. 
In this case, $\mathcal{L}_{\mathrm{bc}} \equiv 0$ and the 
boundary loss term can be omitted from the training objective. 
This lift-based enforcement is applicable when the boundary 
data $u_D$ is sufficiently smooth; for problems with 
discontinuous or piecewise-defined boundary conditions, the 
explicit penalty formulation~\eqref{eq:bc_loss} is employed 
instead.
\end{remark}

\begin{remark} \label{rem:soft_bc}
The present framework employs soft boundary condition enforcement via $\mathcal{L}_{\mathrm{bc}}(\theta)$ rather than hard enforcement techniques, such as the unified PINN approach of Luong et al.~\cite{Luong2024}, where the network output is algebraically constrained to satisfy essential boundary conditions exactly through trial functions. While hard enforcement eliminates the need for tuning $w_{\mathrm{bc}}$ and guarantees exact Dirichlet satisfaction, soft enforcement is preferred here for several reasons: (i) it provides flexibility in balancing competing objectives during multi-phase training transitions; (ii) the selective physics enforcement strategy excludes boundary-adjacent regions from PDE residual computation, and a separate boundary loss allows explicit control over boundary condition fidelity in these zones; (iii) since SUPG-YZ$\beta$ solutions already satisfy boundary conditions by construction, boundary data is implicitly learned through $\mathcal{L}_{\text{data}}(\theta)$, making $\mathcal{L}_{\text{bc}}(\theta)$ serve primarily as a regularization mechanism during physics-dominant training phases. 
\end{remark}

\subsection{Enhanced Selective Physics Enforcement Strategy}
\label{sec:Selective}
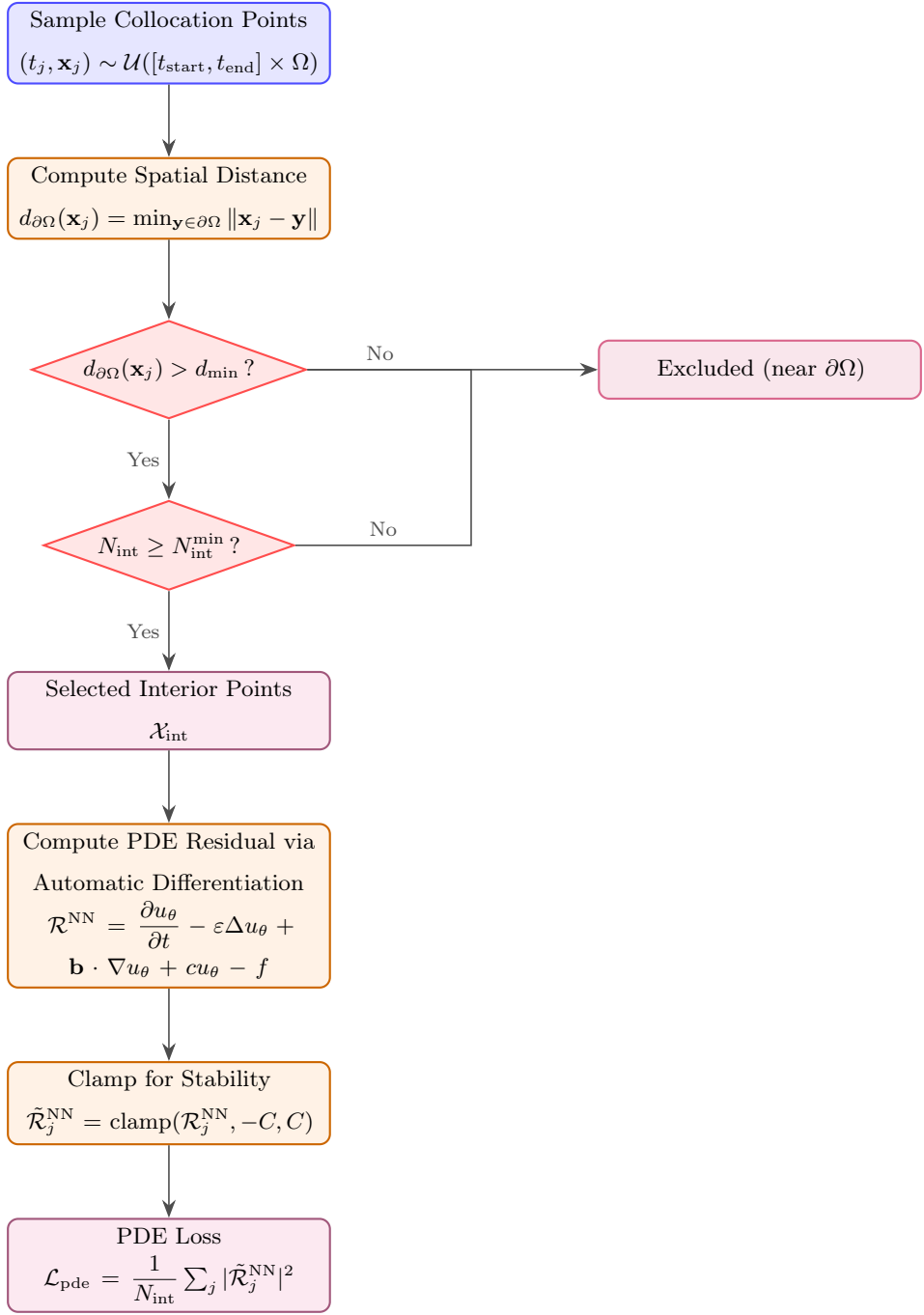
\begin{figure}[!htbp]
    \centering
    \begin{tikzpicture}[
        node distance=1.0cm,
        >={Stealth[length=2.5mm]},
        input/.style={
            rectangle, rounded corners=4pt, 
            draw=blue!70, line width=0.8pt, fill=blue!10, 
            text width=4.2cm, minimum height=0.9cm, 
            align=center, font=\small
        },
        process/.style={
            rectangle, rounded corners=4pt, 
            draw=orange!80!black, line width=0.8pt, 
            fill=orange!10, 
            text width=4.2cm, minimum height=0.9cm, 
            align=center, font=\small
        },
        decision/.style={
            diamond, aspect=2.8,
            draw=red!70, line width=0.8pt, fill=red!10, 
            align=center, font=\small,
            inner sep=2pt
        },
        output/.style={
            rectangle, rounded corners=4pt, 
            draw=magenta!60!black, line width=0.8pt, 
            fill=magenta!10, 
            text width=4.2cm, minimum height=0.9cm, 
            align=center, font=\small
        },
        exclude/.style={
            rectangle, rounded corners=4pt, 
            draw=purple!60, line width=0.8pt, fill=purple!10, 
            text width=4.2cm, minimum height=0.8cm, 
            align=center, font=\small
        },
        arrow/.style={draw=black!70, line width=0.6pt, ->},
        line/.style={draw=black!70, line width=0.6pt},
        yesno/.style={font=\footnotesize, text=black!70}
    ]
    
    \node[input] (input) {Sample Collocation Points\\
    $(t_j, \mathbf{x}_j) \sim \mathcal{U}(
    [t_\text{start}, t_\text{end}] \times \Omega)$};
    
    \node[process, below=of input] (distance) {
    Compute Spatial Distance\\
    $d_{\partial\Omega}(\mathbf{x}_j) = 
    \min_{\mathbf{y} \in \partial\Omega} 
    \|\mathbf{x}_j - \mathbf{y}\|$};
    
    \node[decision, below=1.1cm of distance] (dec1) {
    $d_{\partial\Omega}(\mathbf{x}_j) 
    > d_\text{min}$\,?};
    
    \node[decision, below=1.1cm of dec1] (dec2) {
    $N_{\mathrm{int}} \geq 
    N_{\text{int}}^{\text{min}}$\,?};
    
    \node[output, below=1.1cm of dec2] (selected) {
    Selected Interior Points\\
    $\mathcal{X}_{\mathrm{int}}$};
    
    \node[process, below=of selected] (residual) {
    Compute PDE Residual via Automatic Differentiation\\
    $\mathcal{R^\text{NN}} = \dfrac{\partial u_\theta}{\partial t} 
    - \varepsilon\Delta u_\theta 
    + \mathbf{b}\cdot\nabla u_\theta 
    + cu_\theta - f$};
    
    \node[process, below=of residual] (clamp) {
    Clamp for Stability\\
    $\tilde{\mathcal{R}}_j^\text{NN} = 
    \mathrm{clamp}(\mathcal{R}_j^\text{NN}, -C, C)$};
    
    \node[output, below=of clamp] (pdeloss) {
    PDE Loss\\
    $\mathcal{L}_{\mathrm{pde}} = 
    \dfrac{1}{N_\text{int}}
    \sum_{j}|\tilde{\mathcal{R}}_j^\text{NN}|^2$};
    
    \node[exclude, right=4.0cm of dec1] (excluded) {
    Excluded (near $\partial\Omega$)};
    
    \draw[arrow] (input) -- (distance);
    \draw[arrow] (distance) -- (dec1);
    \draw[arrow] (dec1) -- node[left, yesno] {Yes} (dec2);
    \draw[arrow] (dec2) -- node[left, yesno] {Yes} 
    (selected);
    \draw[arrow] (selected) -- (residual);
    \draw[arrow] (residual) -- (clamp);
    \draw[arrow] (clamp) -- (pdeloss);
    
    \draw[arrow] (dec1.east) -- 
    node[above, yesno, xshift=-1cm] {No} 
    (excluded.west);
    
    \draw[line] (dec2.east) -- 
    node[above, yesno] {No} ++(2.42,0) 
    |- (dec1.east);
    
    \end{tikzpicture}
    \caption{Flowchart illustrating the selective physics 
    enforcement strategy. At each training iteration, 
    spatiotemporal collocation points 
    $(t_j, \mathbf{x}_j)$ are sampled uniformly at random 
    within the training time window and spatial domain. 
    Only points whose spatial distance to the boundary 
    satisfies $d_{\partial\Omega}(\mathbf{x}_j) > 
    d_\text{min}$ and for which the minimum count criterion 
    $N_\text{int} \geq N_\text{int}^\text{min}$ holds are 
    retained for PDE residual computation. Boundary-adjacent 
    points are excluded to preserve numerical stability near 
    sharp gradients. The strong-form residual, including the 
    temporal derivative $\partial u_\theta / \partial t$ 
    obtained via automatic differentiation, is clamped to 
    $[-C, C]$ before computing the loss to prevent gradient 
    explosion.}
    \label{fig:selective_physics_flowchart}
\end{figure}

The key innovation lies in the selective application of physics constraints, which addresses a fundamental challenge in hybrid PINN training---determining where and when to apply physics-based learning. Instead of enforcing the PDE residuals uniformly over the entire domain, we adopt a multi-criteria selection strategy that limits physics enforcement to regions where it is most effective. The overall procedure is summarized in Figure~\ref{fig:selective_physics_flowchart}.

\subsubsection{Distance-based Selection Criterion} \label{sec:distance-based}

We introduce a distance function that measures the proximity of a point
$\mathbf{x} \in \Omega$ to the boundary $\partial\Omega$ for a
rectangular domain $\Omega = [a,b]^2$:
\begin{equation}
d_{\text{boundary}}(\mathbf{x})= \min_{i} \left\{ \frac{x_i - a}{b - a}, \frac{b - x_i}{b - a} \right\},
\label{eq:distance_function}
\end{equation}
where the minimum is taken over all spatial coordinates. This computes the minimum normalized distance to any boundary face, where each term represents the distance to the nearest boundary in the $i$-th coordinate direction, scaled by the domain width $(b-a)$. We apply the following distance criterion:
\begin{equation}
d_{\text{boundary}}(\mathbf{x}) > d_{\text{min}},
\label{eq:distance_criterion}
\end{equation}
where $d_{\text{min}}$ is chosen such that an optimal balance between physics enforcement and boundary layer preservation is achieved. This selective enforcement ensures that only interior points satisfying the distance criterion contribute to the PDE residual loss $\mathcal{L}_{\text{pde}}(\theta)$, thereby enhancing numerical stability and preserving boundary integrity.

\begin{remark}
\label{rem:distance_normalization}
The threshold parameter $d_{\min}$ is defined as a normalized quantity relative to a characteristic length scale $L$ of the domain. For a general bounded domain $\Omega \subset \mathbb{R}^n$, let
\begin{equation}
d(\mathbf{x}, \partial\Omega) = \inf_{\mathbf{y} \in \partial\Omega} \|\mathbf{x} - \mathbf{y}\|_2
\end{equation}
denote the Euclidean distance from a point $\mathbf{x} \in \Omega$ to the boundary $\partial\Omega$. The normalized distance is then defined as
\begin{equation}
\tilde{d}(\mathbf{x}) = \frac{d(\mathbf{x}, \partial\Omega)}{L}.
\end{equation}
For rectangular domains $\Omega = [0, L_1] \times [0, L_2] \times \cdots \times [0, L_n]$, the characteristic length is $L = \max_{1 \leq i \leq n} L_i$, and the distance simplifies to
\begin{equation}
d(\mathbf{x}, \partial\Omega) = \min_{1 \leq i \leq n} \min(x_i, L_i - x_i).
\end{equation}
The selective enforcement condition $\tilde{d}(\mathbf{x}) > d_{\min}$ ensures scale-invariant behavior across problems with different domain sizes. For the unit domains employed in our numerical examples ($L=1$), the normalized and absolute distances coincide.
\end{remark}

\subsubsection{Implementation Details and Computational Considerations}

The selective physics enforcement is implemented through an efficient masking operation based on the actual code implementation:
\begin{equation}
\mathcal{S}_{\text{physics}}(\mathbf{x}) = \mathbb{I}_{d_{\text{boundary}}(\mathbf{x}) > d_{\min}}, \quad \text{provided that } N_{\text{int}} \geq N_{\text{int}}^{\text{min}},
\label{eq:selection_function_impl}
\end{equation}
where the indicator function $\mathbb{I}_{d_{\text{boundary}}(\mathbf{x}) > d_{\min}}$ 
selects interior points satisfying the distance criterion, and the global condition 
$N_{\text{int}} \geq N_{\text{int}}^{\text{min}}$ ensures a minimum number of interior points for stable  PDE loss computation.

The identification of boundary points is performed using a vectorized 
distance computation. For a unit square domain $[0,1]^2$:
\begin{align}
d_{\text{left}} &= x_1, \quad d_{\text{right}} = 1 - x_1, \label{eq:37}\\
d_{\text{bottom}} &= x_2, \quad d_{\text{top}} = 1 - x_2, \label{eq:38} \\
d_{\text{boundary}} &= \min(d_{\text{left}}, d_{\text{right}}, 
                         d_{\text{bottom}}, d_{\text{top}}). \label{eq:39}
\end{align}

\begin{remark}
\label{rem:unit_domain_distance}
Equations~\eqref{eq:37}--\eqref{eq:39} are specific to the unit square domain $\Omega = [0,1]^2$ employed in the numerical examples. For a general rectangular domain $\Omega = [0, L_1] \times [0, L_2]$, these expressions generalize to $d_{\mathrm{left}} = x_1$, $d_{\mathrm{right}} = L_1 - x_1$, $d_{\mathrm{bottom}} = x_2$, $d_{\mathrm{top}} = L_2 - x_2$, and the normalized distance is $\tilde{d}(\mathbf{x}) = d_{\mathrm{boundary}} / L$ with $L = \max(L_1, L_2)$, as detailed in Remark~\ref{rem:distance_normalization}.
\end{remark}

The number of PDE collocation points, $N_{\mathrm{pde}}$, is a
problem-dependent hyperparameter whose value is selected to balance
computational cost with adequate physics sampling density. In general,
$N_{\mathrm{pde}}$ should be chosen large enough to ensure sufficient
coverage of the interior domain---particularly in problems featuring thin
layers or steep gradients---while remaining small enough to avoid excessive
computational overhead from second-order automatic differentiation. The
specific values adopted for each test case are reported alongside the
corresponding numerical examples in~Section~\ref{sec:section5}. Similarly, the boundary exclusion threshold $d_{\min}$ and the minimum
interior point count $N_{\mathrm{int}}^{\min}$ (see
Eq.~\eqref{eq:selection_function_impl}) are problem-dependent
parameters; their values for each test case are reported in
Section~\ref{sec:section5}.

\subsection{Multi-phase Weight Evolution}

To balance data fidelity with physics enforcement throughout the training process, we implement an adaptive weighting strategy in which the relative importance of each loss term evolves based on training progress. The relative weights assigned to the data loss $\mathcal{L}_{\text{data}}$, PDE residual loss $\mathcal{L}_{\text{pde}}$, and boundary loss $\mathcal{L}_{\text{bc}}$ evolve dynamically through distinct phases based on the training progress, e.g.,:
\begin{align}
\text{Phase I (Data-dominant):} \quad &w_{\text{data}} = w_{\text{data}}^\text{I}, \quad w_{\text{pde}} = w_{\text{pde}}^\text{I}, \quad w_{\text{bc}} = w_{\text{bc}}^\text{I}, \label{eq:weights_phase1}\\
\text{Phase II (Transition):} \quad & w_{\text{data}} = w_{\text{data}}^\text{II}, \quad w_{\text{pde}} = w_{\text{pde}}^\text{II}, \quad w_{\text{bc}} = w_{\text{bc}}^\text{II}, \label{eq:weights_phase2}\\
\text{Phase III (Physics-dominant):} \quad & w_{\text{data}} = w_{\text{data}}^\text{III}, \quad w_{\text{pde}} = w_{\text{pde}}^\text{III}, \quad w_{\text{bc}} = w_{\text{bc}}^\text{III}, \label{eq:weights_phase3}
\end{align}
where the superscripts of weights indicate the phase associated with them.

The phase boundaries and loss-weight selections (see Eqs.~\eqref{eq:weights_phase1}--\eqref{eq:weights_phase3}) to be used are selected based on numerical experiments. For example, in one-dimensional CDR problems, approximately $1000-2000$ total epochs can be sufficient for stable and accurate convergence, whereas two-dimensional cases typically require around $2500-5000$ epochs. These values should, therefore, be interpreted as empirically informed guidelines rather than universally optimal hyperparameters, and may need adjustment for other PDEs or geometries. While more principled dynamic weighting strategies---such as gradient-norm balancing or uncertainty-based scaling---have recently been explored in the broader PINN literature, our empirically-tuned three-phase approach proved robust across all benchmark problems considered. Integrating such adaptive mechanisms into hybrid FEM--PINN frameworks represents a promising direction for future research that could potentially enhance the method's generalizability.

The rationale behind this progressive strategy is threefold: (i) initializing  with data-dominant learning prevents early training instabilities that arise  when physics terms dominate before a reasonable solution approximation is  established; (ii) the gradual transition allows the network to smoothly adapt  to the changing optimization landscape; and (iii) physics-dominant refinement 
in the final phase enables the hybrid approach to potentially improve upon  the FEM reference solution by enforcing the governing PDE directly. The  specific weight values for each phase are problem-dependent and will be  reported in the numerical examples.

\subsection{FEM-PINN Data Transfer and Computational Interface}

The seamless integration between \texttt{FEniCS} and \texttt{PyTorch} is achieved through a carefully designed data transfer protocol that maintains numerical accuracy while enabling efficient GPU-accelerated training. Spatial coordinates $\{\mathbf{x}_i\}_{i=1}^{N}$ are extracted from the FEM mesh using \texttt{FEniCS} function space mappings via \texttt{V.tabulate\_dof\_coordinates()}, and paired with snapshot times $\{t_k\}_{k=1}^{K_s}$ to construct space-time training tuples $\{(t_k, \mathbf{x}_i, u_{\text{SUPG-YZ}\beta}^k(\mathbf{x}_i))\}$. This ensures exact correspondence between finite element nodes and the points used for loss computation across multiple time levels. High-accuracy evaluation of $u_{\text{SUPG-YZ}\beta}^k(\mathbf{x})$ at each snapshot is performed using \texttt{FEniCS}'s built-in point evaluation capabilities, maintaining the accuracy of the finite element approximation. Efficient detection of boundary points using spatial tolerance $\epsilon_{\text{boundary}} = 10^{-8}$ and geometric proximity tests are employed.

The hybrid framework leverages \texttt{PyTorch}'s automatic differentiation capabilities for computing space-time PDE residuals. The network accepts three-dimensional input $(t, x_1, x_2)$ and produces a scalar output $u_{\text{NN}}(t, x_1, x_2)$. The temporal derivative $\partial u_{\text{NN}} / \partial t$ and spatial gradients $\nabla u_{\text{NN}}$ are computed simultaneously using \texttt{torch.autograd.grad} with \texttt{create\_graph=True}, enabling higher-order differentiation through the same computational graph. The Laplacian $\Delta u_{\text{NN}}$ is obtained through successive differentiation of spatial first derivatives. The strong-form PDE residual, $\mathcal{R^\text{NN}}$, is evaluated at randomly sampled interior collocation points within the training time window $[t_{\text{start}}, t_{\text{end}}]$, with boundary-proximal points excluded via a distance mask ($d > 0.02$) to avoid interference with the explicit boundary condition loss. Efficient batched computation of all derivatives across multiple collocation points leverages GPU parallelization.

\subsection{Optimization Strategy and Training Configuration}
\label{sec:optimization}

Training is performed using the \texttt{AdamW} stochastic 
optimization algorithm~\cite{loshchilov2017decoupled}, adopted for 
its improved generalization capabilities and decoupled weight decay 
mechanism compared to standard 
\texttt{Adam}~\cite{kingma2015adam}. The optimizer is configured 
with momentum parameters $\beta_1 = 0.9$, $\beta_2 = 0.999$, and 
a weight decay coefficient of $\lambda = 10^{-7}$, which serves as 
implicit $L^2$ regularization to prevent overfitting and promote 
smoother solution profiles---particularly beneficial for interior 
layer problems where solution regularity is crucial.

The initial learning rate is scaled according to the batch size 
following a square-root scaling rule:
\begin{equation}
	\alpha = 8 \times 10^{-5} \times \sqrt{\frac{B}{256}},
	\label{eq:lr_scaling}
\end{equation}
where $B$ denotes the mini-batch size. This scaling heuristic 
maintains stable gradient variance across different batch 
configurations and prevents oscillations when different loss 
components have conflicting gradient directions---a phenomenon 
amplified in the time-dependent setting where the temporal 
derivative introduces an additional source of gradient competition.

A \texttt{ReduceLROnPlateau} scheduler monitors the epoch-averaged 
total loss and reduces the learning rate by a multiplicative factor 
$\gamma = 0.9$ after $150$ epochs of stagnation, with a minimum 
learning rate floor of $\alpha_{\min} = 10^{-6}$. The extended 
patience period accommodates the inherent variability in 
physics-informed loss landscapes, where temporary plateaus often 
precede significant improvements as the network learns to balance 
multiple objectives. The multiplicative decay 
$\alpha \leftarrow \gamma \alpha$ provides a gradual annealing that 
preserves training momentum while enabling fine-grained convergence 
in later epochs.

\subsubsection{Gradient Stabilization and Numerical Safeguards}

To ensure gradient stability and suppress numerical divergence 
during training, we implement comprehensive numerical safeguards 
at three levels: gradient, loss, and residual.

Gradient clipping is applied with a maximum norm constraint 
$\|\nabla_\theta \mathcal{L}\|_2 \leq 1.0$, preventing catastrophic 
parameter updates in regions of sharp solution gradients where 
second-order derivatives may become large. A vectorized gradient 
validation mechanism inspects the concatenated parameter gradient 
tensor at each optimization step:
\begin{equation}
	\text{if } \neg\,\texttt{isfinite}\!\left(\bigoplus_{l=1}^{L} 
	\nabla_{\theta_l} \mathcal{L}\right) \text{ then discard update 
		and increment NaN counter},
\end{equation}
where $\bigoplus$ denotes concatenation across all $L$ parameter 
groups. Training is terminated if the cumulative NaN count exceeds 
$10$, indicating persistent numerical instability. At the loss 
level, individual components are similarly guarded:
\begin{equation}
	\text{if } \texttt{isnan}(\mathcal{L}_{\text{hybrid}}) \text{ or } 
	\texttt{isinf}(\mathcal{L}_{\text{hybrid}}) \text{ then } 
	\mathcal{L}_{\text{hybrid}} \leftarrow \mathcal{L}_{\text{data}},
\end{equation}
ensuring that training can recover from transient PDE residual 
spikes by falling back to the data fidelity loss alone.

Additionally, value clipping via the $\texttt{clamp}(x, a, b)$ 
function is employed to truncate extreme values at multiple stages 
of the computation:
\begin{align}
	\text{Residual clipping:} \quad 
	&\mathcal{R}^{\text{NN}}(\mathbf{x}) \leftarrow 
	\texttt{clamp}\bigl(\mathcal{R}^{\text{NN}}(\mathbf{x}),\, 
	-10.0,\, 10.0\bigr), \label{eq:residual_clip}\\
	\text{Loss clipping:} \quad 
	&\mathcal{L}_{\text{hybrid}} \leftarrow 
	\texttt{clamp}\bigl(\mathcal{L}_{\text{hybrid}},\, 0.0,\, 
	100.0\bigr), \label{eq:loss_clip}
\end{align}
where the $\texttt{clamp}$ function is defined as:
\begin{equation}
	\texttt{clamp}(x, a, b) = \begin{cases} 
		a, & \text{if } x < a, \\
		x, & \text{if } a \leq x \leq b, \\
		b, & \text{if } x > b.
	\end{cases}
\end{equation}
The residual clipping bounds of $\pm 10.0$ are chosen to suppress 
outlier contributions from collocation points near sharp layers 
while retaining physically meaningful gradient information.

\subsubsection{Batch Processing and Mixed-Precision Training}

All training data---comprising $K_s \times N_h$ space--time 
tuples---are preloaded entirely onto the GPU at initialization to 
eliminate host-to-device transfer overhead. Epoch-level shuffling is 
performed via in-place random permutation on the device using 
\texttt{torch.randperm}, yielding $\lceil N_{\text{data}} / B \rceil$ 
mini-batches per epoch, each receiving a full optimizer update. PDE residual evaluation is performed periodically during each epoch
rather than at every mini-batch, thereby amortizing the computationally
expensive second-order automatic differentiation over multiple data
batches. In the default configuration, the residual is computed on the
first mini-batch of every epoch; however, this frequency is treated as a
tunable parameter and may be adjusted on a per-problem basis---for
instance, evaluating the residual every $k$-th mini-batch when a higher
physics sampling rate is beneficial or when the per-batch cost must be
further reduced. The specific evaluation frequency adopted for each test
case is reported in Section~\ref{sec:section5}. In all cases,
$N_{\mathrm{pde}}$ collocation points are randomly sampled within the
training time window
$[t_{\mathrm{start}}, t_{\mathrm{end}}] \times \Omega$ at each
evaluation instance.

On CUDA-enabled devices, mixed-precision training via Automatic 
Mixed Precision (AMP) is selectively applied: the data fidelity and 
boundary condition forward passes execute in \texttt{float16} 
arithmetic for computational throughput, while PDE residual 
computations involving higher-order automatic differentiation are 
maintained in full \texttt{float32} precision to preserve the 
numerical accuracy of second-derivative calculations. A dynamic 
\texttt{GradScaler} adjusts the loss scaling factor to prevent 
gradient underflow in reduced precision. The best model state, 
determined by the minimum epoch-averaged total loss, is maintained 
in CPU memory and checkpointed to disk every $100$ epochs for fault 
tolerance. Training proceeds for $N_{\text{epochs}} = 5000$ epochs 
across all test cases considered in this work.

\subsection{Hyperparameter Selection}
\label{sec:hyperparams}

The network architecture employs $n_r$ residual blocks with hidden 
dimension $n_h$ and $n_{\text{F}}$ random Fourier features at scale 
$\sigma$; the specific values adopted for each test problem are 
reported in Section~\ref{sec:section5}. Across all examples, the 
architectural capacity is chosen to balance representational power 
against overfitting risk: too few parameters would struggle to 
resolve sharp interfaces, while substantially larger networks would 
risk memorizing FEM data without learning the underlying physics. 
The Fourier feature embedding provides the spectral bias correction 
necessary for representing high-frequency layer 
structures~\cite{tancik2020fourier}.

Training utilizes the last $K_s$ temporal snapshots from the 
SUPG-YZ$\beta$ solution, covering a time window 
$[t_{\text{start}}, t_{\text{end}}] \subset [0, t_f]$ near the 
terminal time. This windowed approach focuses the PINN's corrective 
capacity on the near-terminal dynamics where advected fronts have 
developed their most challenging features, while providing 
sufficient temporal context for the network to learn meaningful 
$\partial u / \partial t$ behavior through automatic differentiation. 
The total training dataset comprises $K_s \times N_h$ space--time 
points; the values of $K_s$ and the batch size $B$ are 
problem-dependent and are specified alongside each numerical example.

For problems with homogeneous Dirichlet data, the lift-based 
distance function enforcement described in 
Remark~\ref{remark:lift} renders the boundary loss unnecessary 
($w_{\text{bc}} = 0$). When the boundary data is 
piecewise-discontinuous or non-homogeneous, an explicit boundary 
loss with $n_{\text{bc}}$ collocation points per edge is employed 
instead. Similarly, the number of PDE collocation points 
$N_{\text{pde}}$ and the boundary exclusion threshold $d_{\min}$ are 
selected on a per-problem basis to balance physics sampling density 
against computational cost; representative values are given in 
Section~\ref{sec:section5}.

	\def\nh{128}
	\def\nr{8}
	\def\nF{24}
	\def\inputdim{3}
	\def\fourierdim{51}
	\def\narrowdim{64}
	
	\begin{figure}[!htbp]
		\centering
		\resizebox{0.95\textwidth}{!}{
		\begin{tikzpicture}[
			>=Stealth,
			neuron/.style={circle, draw, minimum size=11pt, inner sep=0pt, line width=0.5pt},
			input neuron/.style={neuron, fill=blue!20, draw=blue!50},
			fourier neuron/.style={neuron, fill=violet!18, draw=violet!50},
			hidden neuron/.style={neuron, fill=cyan!15, draw=cyan!50},
			res neuron/.style={neuron, fill=green!15, draw=green!45!black},
			narrow neuron/.style={neuron, fill=red!12, draw=red!45},
			output neuron/.style={neuron, fill=red!22, draw=red!65, minimum size=16pt},
			conn/.style={line width=0.12pt, draw=black!12},
			skip conn/.style={line width=0.8pt, draw=orange!65!black, dashed, 
				-{Stealth[length=3.5pt]}},
			layer label/.style={font=\small\bfseries, align=center},
			sublabel/.style={font=\scriptsize, text=gray!65!black, align=center},
			]
			
			\def\neuronsp{0.50}
			\def\layersep{2.0}
			
			\def\xI{0}
			\node[input neuron] (I-1) at (\xI, 1.2)  {\tiny$t$};
			\node[input neuron] (I-2) at (\xI, 0)    {\tiny$x_1$};
			\node[input neuron] (I-3) at (\xI, -1.2) {\tiny$x_2$};
			
			\def\xF{2.0}
			\foreach \i in {1,...,4} {
				\pgfmathsetmacro{\yy}{2.1 - (\i-1)*0.55}
				\node[fourier neuron] (F-t\i) at (\xF, \yy) {};
			}
			\node[font=\normalsize, text=gray!60] (F-dots) at (\xF, 0) {$\vdots$};
			\foreach \i in {1,...,4} {
				\pgfmathsetmacro{\yy}{-0.55 - (\i-1)*0.55}
				\node[fourier neuron] (F-b\i) at (\xF, \yy) {};
			}
			
			\def\xH{4.0}
			\foreach \i in {1,...,5} {
				\pgfmathsetmacro{\yy}{2.65 - (\i-1)*0.55}
				\node[hidden neuron] (H-t\i) at (\xH, \yy) {};
			}
			\node[font=\normalsize, text=gray!60] (H-dots) at (\xH, 0) {$\vdots$};
			\foreach \i in {1,...,5} {
				\pgfmathsetmacro{\yy}{-0.55 - (\i-1)*0.55}
				\node[hidden neuron] (H-b\i) at (\xH, \yy) {};
			}
			
			\foreach \blk in {1,2,3} {
				\pgfmathsetmacro{\xR}{4.0 + \blk * 2.0}
				\foreach \i in {1,...,5} {
					\pgfmathsetmacro{\yy}{2.65 - (\i-1)*0.55}
					\node[res neuron] (R\blk-t\i) at (\xR, \yy) {};
				}
				\node[font=\normalsize, text=gray!60] (R\blk-dots) at (\xR, 0) {$\vdots$};
				\foreach \i in {1,...,5} {
					\pgfmathsetmacro{\yy}{-0.55 - (\i-1)*0.55}
					\node[res neuron] (R\blk-b\i) at (\xR, \yy) {};
				}
			}
			
			\def\xEll{9.0}
			\node[font=\Large\bfseries, text=green!40!black] at (\xEll, -0.0) {$\cdots$};
			
			\def\xN{12.6}
			\foreach \i in {1,...,4} {
				\pgfmathsetmacro{\yy}{1.85 - (\i-1)*0.55}
				\node[narrow neuron] (N-t\i) at (\xN, \yy) {};
			}
			\node[font=\normalsize, text=gray!60] (N-dots) at (\xN, 0) {$\vdots$};
			\foreach \i in {1,...,4} {
				\pgfmathsetmacro{\yy}{-0.55 - (\i-1)*0.55}
				\node[narrow neuron] (N-b\i) at (\xN, \yy) {};
			}
			
			\def\xO{14.6}
			\node[output neuron] (O-1) at (\xO, 0) {\small$u$};
			
			
			\foreach \i in {1,2,3} {
				\foreach \j in {1,...,4} {
					\draw[conn] (I-\i) -- (F-t\j);
					\draw[conn] (I-\i) -- (F-b\j);
				}
			}
			
			\foreach \src in {t1,t2,t3,t4,b1,b2,b3,b4} {
				\foreach \dst in {t1,t2,t3,t4,t5,b1,b2,b3,b4,b5} {
					\draw[conn] (F-\src) -- (H-\dst);
				}
			}
			
			\foreach \src in {t1,t2,t3,t4,t5,b1,b2,b3,b4,b5} {
				\foreach \dst in {t1,t2,t3,t4,t5,b1,b2,b3,b4,b5} {
					\draw[conn] (H-\src) -- (R1-\dst);
				}
			}
			
			\foreach \src in {t1,t2,t3,t4,t5,b1,b2,b3,b4,b5} {
				\foreach \dst in {t1,t2,t3,t4,t5,b1,b2,b3,b4,b5} {
					\draw[conn] (R1-\src) -- (R2-\dst);
				}
			}
			
			\foreach \src in {t1,t2,t3,t4,t5,b1,b2,b3,b4,b5} {
				\foreach \dst in {t1,t2,t3,t4,t5,b1,b2,b3,b4,b5} {
					\draw[conn] (R2-\src) -- (R3-\dst);
				}
			}
			
			\foreach \src in {t1,t2,t3,t4,t5,b1,b2,b3,b4,b5} {
				\foreach \dst in {t1,t2,t3,t4,b1,b2,b3,b4} {
					\draw[conn] (R3-\src) -- (N-\dst);
				}
			}
			
			\foreach \src in {t1,t2,t3,t4,b1,b2,b3,b4} {
				\draw[conn, line width=0.25pt, draw=black!20] (N-\src) -- (O-1);
			}
			
			\draw[skip conn] 
			([yshift=3pt]H-t1.north) to[out=70, in=110] ([yshift=3pt]R1-t1.north);
			\draw[skip conn] 
			([yshift=3pt]R1-t1.north) to[out=70, in=110] ([yshift=3pt]R2-t1.north);
			\draw[skip conn] 
			([yshift=3pt]R2-t1.north) to[out=70, in=110] ([yshift=3pt]R3-t1.north);
			
			\draw[decorate, decoration={brace, amplitude=7pt, raise=4pt}, 
			green!50!black, line width=0.7pt]
			(5.5, 3.6) -- (10.5, 3.6)
			node[midway, above=13pt, font=\small\bfseries, text=green!40!black] 
			{$\times\;\nr$ Residual Blocks};
			
			\def\labelY{-3.6}
			
			\node[layer label, text=blue!60!black] at (\xI, \labelY) {Input};
			\node[sublabel] at (\xI, \labelY - 0.4) {$\mathbb{R}^{\inputdim}$};
			\node[sublabel] at (\xI, \labelY - 0.75) {$(t, x_1, x_2)$};
			
			\node[layer label, text=violet!60!black] at (\xF, \labelY) {Fourier};
			\node[sublabel] at (\xF, \labelY - 0.4) {$\mathbb{R}^{\fourierdim}$};
			\node[sublabel] at (\xF, \labelY - 0.75) {$[\mathbf{z};\sin;\cos]$};
			
			\node[layer label, text=cyan!45!black] at (\xH, \labelY) {Input Layer};
			\node[sublabel] at (\xH, \labelY - 0.4) {$\nh$ neurons};
			\node[sublabel] at (\xH, \labelY - 0.75) {SiLU};
			
			\node[layer label, text=green!40!black] at (8, \labelY) {Hidden Layers};
			\node[sublabel] at (8, \labelY - 0.4) {$\nh$ neurons/block};
			\node[sublabel] at (8, \labelY - 0.75) {Skip $+$ LayerNorm $+$ SiLU};
			
			\node[layer label, text=red!55!black] at (\xN, \labelY) {Narrowing};
			\node[sublabel] at (\xN, \labelY - 0.4) {$\narrowdim$ neurons};
			\node[sublabel] at (\xN, \labelY - 0.75) {SiLU};
			
			\node[layer label, text=red!65!black] at (\xO, \labelY) {Output};
			\node[sublabel] at (\xO, \labelY - 0.4) {$1$ neuron};
			\node[sublabel] at (\xO, \labelY - 0.75) {$u_\mathrm{NN}(t,\mathbf{x})$};
		\end{tikzpicture}}
		\caption{Schematic of the PASSC neural network architecture 
			(shown for $n_h=\nh$, $n_r=\nr$, $n_\text{F}=\nF$, $n_\text{sd}=2$).
			The spatiotemporal input 
			$\mathbf{z}=(t,x_1,x_2)\in\mathbb{R}^{\inputdim}$ 
			is mapped through a random Fourier feature embedding to 
			$\boldsymbol{\varphi}(\mathbf{z})\in\mathbb{R}^{\fourierdim}$, 
			projected to a hidden dimension of $\nh$ via the input layer with SiLU 
			activation, and then processed by $\nr$ residual blocks---each comprising 
			two fully connected layers, a skip (residual) connection, and layer 
			normalization. A narrowing layer ($\nh\!\to\!\narrowdim$) followed by a 
			single output neuron produces the predicted solution 
			$u_\mathrm{NN}(t,\mathbf{x})$. }
		\label{fig:nn-architecture}
	\end{figure}
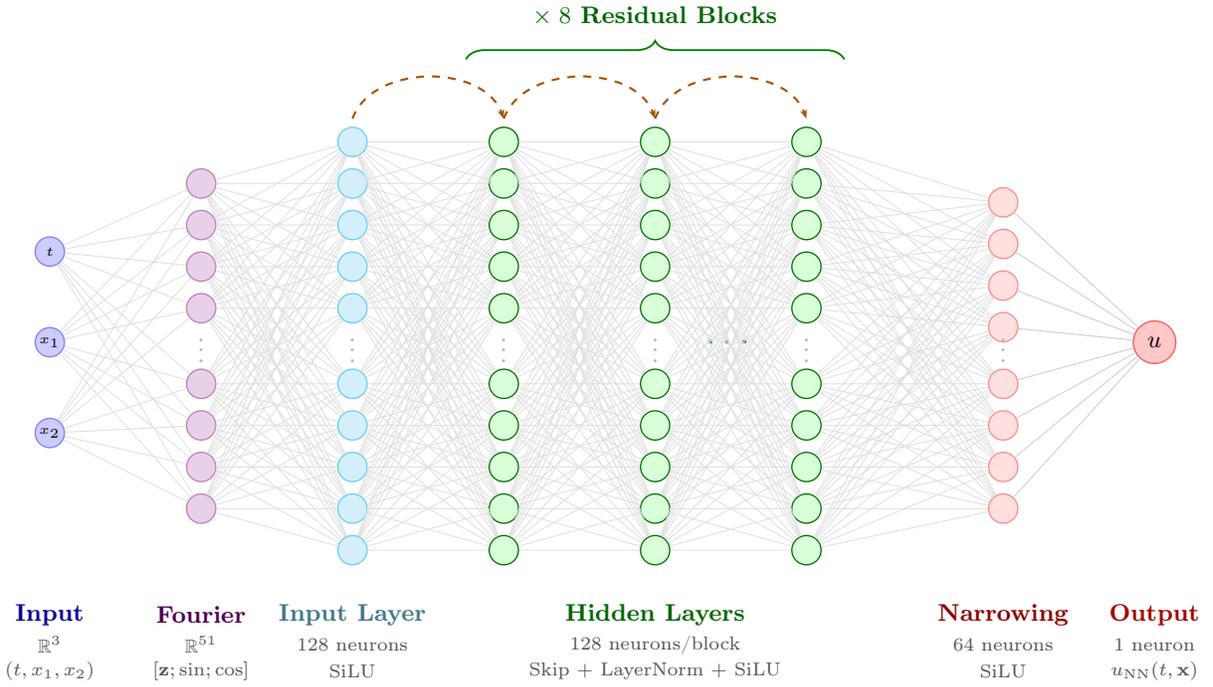

\section{Further Computational Details} \label{sec:section4}

This section describes the computational framework developed to implement and validate the proposed hybrid methodology. It covers the scientific computing environment, spatial discretization schemes, mesh generation strategies, and solution techniques for the resulting linear and nonlinear algebraic systems, as well as the \texttt{PyTorch} environment and \texttt{\texttt{CUDA}} technology. 
    
\subsection{Scientific Computing Environment: \texttt{FEniCS}}

\texttt{FEniCS}~\cite{fenics16, logg2012automated, Alnaes2015} is an open-source scientific computing platform specifically designed for the numerical solution of partial differential equations (PDEs). It provides a robust and flexible infrastructure for developing FEM-based solvers and has been selected as the primary computational environment in this study. Owing to its versatility, \texttt{FEniCS} offers an ideal framework for implementing both stabilized finite element methods and PINNs. Its \texttt{C++} and/or \texttt{Python}-based interface allows for concise expression of complex mathematical models and seamless integration with a wide array of machine learning libraries---an essential capability for optimizing stabilization and shock-capturing parameters within a unified environment. Moreover, \texttt{FEniCS} is compatible with high-performance computing (HPC) systems, supporting parallel computations that are essential for solving large-scale and computationally intensive problems. This scalability is particularly advantageous when addressing time-dependent and three-dimensional PDEs, contributing to the robustness and flexibility of the numerical investigations in this study. Here, \texttt{FEniCS} (2019.2.0.64.dev0) serves as the core computational platform for implementing and validating the standard and stabilized finite element formulations and their integration with PINNs. The SUPG stabilization and YZ$\beta$ shock-capturing techniques are developed within the \texttt{FEniCS} environment and combined with machine learning strategies.

\texttt{FEniCS} also has a newer version called \texttt{FEniCSx}~\cite{dolfinx2023}, which offers enhanced performance through modern \texttt{C++} architecture, improved scalability, and better integration with contemporary scientific computing ecosystems. For more details on both the legacy and new versions of the \texttt{FEniCS} project, interested readers can refer to the project webpage \url{https://fenicsproject.org/}.

\subsection{\texttt{PyTorch} and \texttt{\texttt{CUDA}}}
\texttt{PyTorch}~\cite{Paszke19a} is utilized in this study as the core deep learning framework for implementing PINNs. Its dynamic computation graph, user-friendly interface, and strong support for GPU acceleration via \texttt{\texttt{CUDA}} (Compute Unified Device Architecture) make it a suitable choice for training neural networks in complex, physics-driven settings. \texttt{PyTorch} enables the automatic differentiation of loss functions that include not only standard PDE residuals but also contributions from SUPG stabilization and YZ$\beta$ shock-capturing terms, as introduced in Section 3. This capability is essential for efficiently computing gradients during backpropagation and updating model parameters in a stable manner. In order to enhance computational performance, particularly for high-resolution spatiotemporal problems and deep network architectures, all \texttt{PyTorch}-based PINN models are trained on NVIDIA GPUs using \texttt{\texttt{CUDA}}. The \texttt{\texttt{CUDA}} backend allows substantial speedup in both forward passes and gradient computations, which significantly reduces the total training time. 

The integration between \texttt{PyTorch} and \texttt{FEniCS} is facilitated by the \texttt{Python} ecosystem, allowing seamless data transfer between the finite element solver and the neural network. In particular, quantities such as element-wise parameters, geometric characteristics, and PDE residuals computed in \texttt{FEniCS} are passed as inputs to the \texttt{PyTorch} network to enable localized adaptivity during training. 

For more details on \texttt{PyTorch}, the reader is referred to the official webpage of the project: \url{https://pytorch.org/}.

\subsection{Finite Element Mesh Construction}
In this work, 1D and 2D geometries are considered, and finite element meshes are generated to support the numerical simulations. The built-in \texttt{mshr} module of \texttt{FEniCS} is employed to generate finite element meshes directly within the \texttt{Python} environment. This enables rapid mesh generation and seamless integration into the overall simulation workflow.

\subsection{Time Integration} \label{sec:time_integration}
The spatially semi-discrete system arising from the stabilized finite
element formulation is advanced in time using the backward Euler
method. Given the discrete solution $u^{h,n} \in \mathcal{S}^h_u$
at time level $t^n$, the solution $u^{h,n+1}$ at
$t^{n+1} = t^n + \varDelta t$ is obtained by solving the following
fully discrete variational problem: find
$u^{h,n+1} \in \mathcal{S}^h_u$ such that
\begin{equation}
	\begin{split}
		& \int_\Omega \frac{u^{h,n+1} - u^{h,n}}{\varDelta t}\,w^h\,d\Omega
		+ \int_\Omega \left( \varepsilon\,\nabla u^{h,n+1} \cdot \nabla w^h
		+ (\mathbf{b}^{n+1} \cdot \nabla u^{h,n+1})\,w^h
		+ c\,u^{h,n+1}\,w^h - f^{n+1}\,w^h \right) d\Omega \\[4pt]
		& + \sum_{e=1}^{n_\text{el}} \int_{\Omega^e}
		\tau_\text{SUPG}\,(\mathbf{b}^{n+1} \cdot \nabla w^h)\,
		\mathcal{R}^h(u^{h,n+1})\,d\Omega
		+ \sum_{e=1}^{n_\text{el}} \int_{\Omega^e}
		\nu_\text{SHOC}\,\nabla w^h \cdot \nabla u^{h,n+1}\,d\Omega = 0,
	\end{split}
	\label{eq:fully_discrete}
\end{equation}
for all $w^h \in \mathcal{V}^h_u$, where
$\mathbf{b}^{n+1} = \mathbf{b}(t^{n+1}, \mathbf{x})$ and
$f^{n+1} = f(t^{n+1}, \mathbf{x})$ denote the convection velocity and
source term evaluated at the implicit time level, respectively. For
nonlinear problems in which the convection velocity depends on the
solution itself (e.g., Burgers' equation with
$\mathbf{b} = (u, u)^{\top}$), the implicit dependence
$\mathbf{b}^{n+1} = \mathbf{b}(u^{h,n+1})$ is resolved through the
Newton--Raphson linearization described in the 
next section. The strong-form residual at time
level $n+1$ is defined as
\begin{equation}
	\mathcal{R}^h(u^{h,n+1}) = \frac{u^{h,n+1} - u^{h,n}}{\varDelta t}
	- \varepsilon\,\Delta u^{h,n+1}
	+ \mathbf{b}^{n+1} \cdot \nabla u^{h,n+1}
	+ c\,u^{h,n+1} - f^{n+1}.
	\label{eq:strong_residual_discrete}
\end{equation}

\subsection{Numerical Solution of Linear and Nonlinear Systems of Algebraic Equations} 
The stabilized formulation of the governing equations, as described in Eq.~\eqref{eq:supg_shock}, leads to a system of nonlinear equations. These nonlinear systems are solved using the Newton--Raphson (N--R) method, which provides a robust framework for linearizing and iteratively solving nonlinear PDEs. At each nonlinear iteration, the resulting linear systems are addressed using different strategies depending on the problem size and computational setting. For small-scale problems, direct solvers based on LU (lower-upper) decomposition are employed due to their efficiency and ease of implementation. In serial computations, the built-in \texttt{PETSc} (Portable, Extensible Toolkit for Scientific Computation)~\cite{petsc-user-ref} support in \texttt{FEniCS} is used to access a wide variety of linear algebra routines.  In cases where direct solvers become infeasible due to memory constraints or extreme problem size, iterative solvers such as GMRES (Generalized Minimal Residual Method)~\cite{Saad86a} are adopted. GMRES is supported by ILU (incomplete lower-upper) preconditioning to accelerate convergence in sparse and poorly conditioned systems.

\section{Test Computations} \label{sec:section5}
This section presents a comprehensive set of numerical examples to evaluate the performance and robustness of the proposed hybrid framework. The test problems are carefully selected to cover a wide range of transient scenarios, assessing the accuracy, stability, and adaptive capabilities of the method, particularly in the presence of steep gradients, discontinuities, and shock-like features. Where applicable, comparisons with classical methods are provided to highlight the advantages of the developed approach. All computations are performed on a system equipped with an Intel Core Ultra 9 275HX processor (24 cores) and an NVIDIA GeForce RTX 5070 GPU (8 GB VRAM), running Ubuntu 24.04.3 LTS with CUDA Toolkit 13.0.

\subsection{Example 1.} 
As the first test example, consider the following 1D parabolic initial-boundary value problem~\cite{Gowrisankar2014}:
\begin{equation}
\begin{cases}
u_t - \varepsilon u_{xx} + (1 + x(1-x))u_x = f(t, x), & (t, x) \in (0,1] \times (0,1), \\[6pt]
u(0, x) = u_0(x), & 0 < x < 1, \\[6pt]
u(t,0) = 0, \quad u(t, 1) = 0, & 0 \leq t \leq 1,
\end{cases}
\end{equation}
where the initial data $u_0(x)$ and the source function $f(t, x)$ are chosen to fit the exact solution
\begin{equation}
u(t, x) = \exp(-t)\left(C_1 + C_2 x - \exp\left(\frac{x-1}{\varepsilon}\right)\right),
\end{equation}
with $C_1 = \exp(-1/\varepsilon)$ and $C_2 = 1 - \exp(-1/\varepsilon)$.

\begin{remark}
	In \textit{Example 1}, since the problem is one-dimensional, we simplify the notation by using $x$ instead of $x_1$ to denote the spatial coordinate.
\end{remark}

For the finite element setting, we choose 
$\varepsilon = 10^{-4}$, $n_{\text{el}} = 200$, 
$\varDelta t = 0.0025$ ($N_t = 400$), and Y $= 0.1$. 
The PINN architecture employs $n_h = 128$, $n_r = 8$, 
$n_\text{F} = 24$, and $\sigma = 4.0$. Training uses the last 
$K_s = 10$ temporal snapshots near the terminal time 
$t_\text{f} = 1.0$, with a batch size $256$, an initial learning rate 
$\alpha = 10^{-4}$, and a gradient clipping threshold 
$g_{\max} = 0.3$. The PDE residual is evaluated at 
$N_{\text{pde}} = 512$ randomly sampled interior 
collocation points every third mini-batch to reduce 
computational cost. The phase boundaries and adaptive 
weights are presented in Table~\ref{tab:phase_example_one}.

\begin{table}[htb]
\centering
\caption{Phase boundaries and adaptive weights used for solving \textit{Example 1}.}
\label{tab:phase_example_one}
\begin{tabular}{llccc}
\toprule
\text{Phase} & \text{Epochs} & $w_{\mathrm{data}}$ & $w_{\mathrm{pde}}$ & $w_{\mathrm{bc}}$  \\
\midrule
I   & 0--1499   & 1.0  & 0.5  & 0.1 \\
II  & 1500--2999 & 0.8  & 0.95 & 0.1 \\
III & 3000--4999 & 0.35 & 5.0  & 0.1 \\
\bottomrule
\end{tabular}
\end{table}

\begin{figure}[h!]
    \centering
    \begin{subfigure}[b]{0.49\linewidth}
        \includegraphics[width=\linewidth]{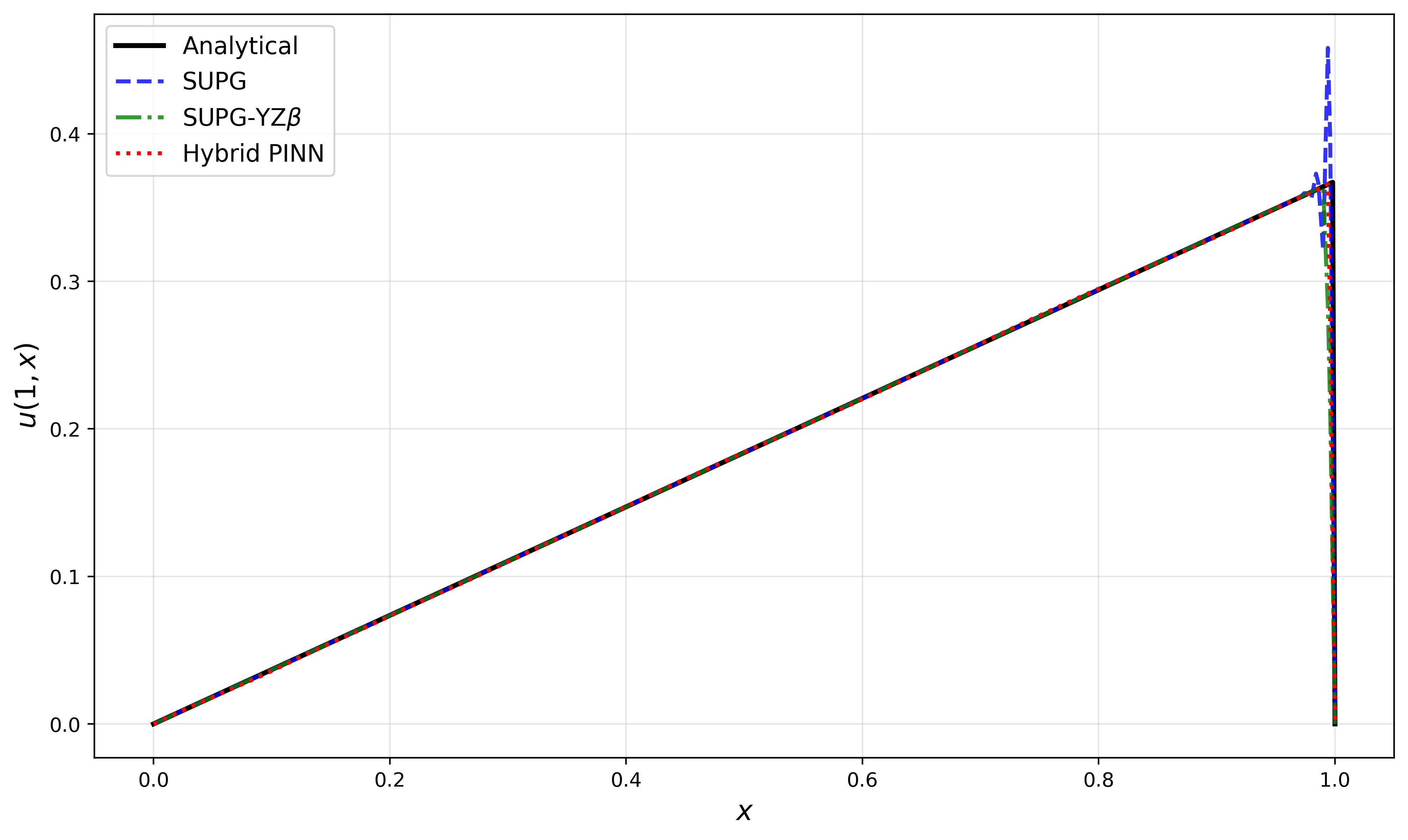}
        \caption{}
        \label{fig:fig1a}
    \end{subfigure}
    \hfill
    \begin{subfigure}[b]{0.49\linewidth}
        \includegraphics[width=\linewidth]{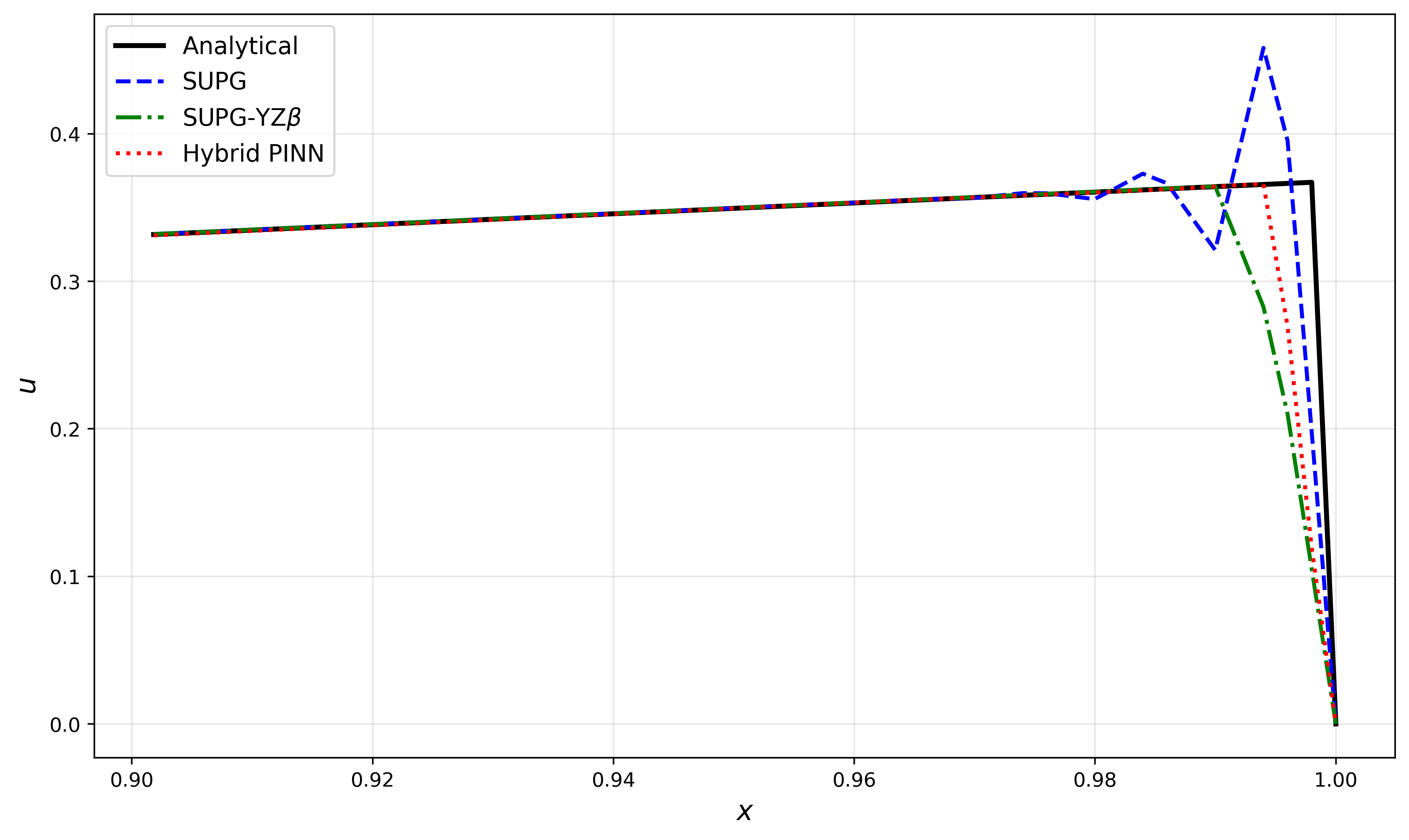}
        \caption{}
        \label{fig:fig1b}
    \end{subfigure}
    \caption{Comparison of the numerical solutions obtained by SUPG, SUPG-YZ$\beta$, and the hybrid PINN with the analytical solution for \textit{Example 1} with $\varepsilon = 10^{-4}$: (\subref{fig:fig1a}) entire domain and (\subref{fig:fig1b}) zoom in the boundary layer region near $x = 1$.}
    \label{fig:1_comparison_fem}
\end{figure}

Figure~\ref{fig:1_comparison_fem} presents the numerical 
solutions at the terminal time $t_\text{f} = 1$ for 
$\varepsilon = 10^{-4}$. On the global scale 
(Figure~\ref{fig:fig1a}), all three numerical 
approximations---SUPG, SUPG-YZ$\beta$, and the hybrid 
PINN---are visually indistinguishable from the analytical 
solution throughout the interior of the domain, where the 
solution varies smoothly and nearly linearly. The critical 
differences emerge in the boundary layer region near 
$x = 1$, which is magnified in Figure~\ref{fig:fig1b}. 
There, the SUPG solution exhibits pronounced node-to-node 
oscillations that extend into the interior, a well-known 
artifact of insufficient crosswind dissipation in the 
presence of exponentially steep gradients of width 
$\mathcal{O}(\varepsilon)$. The addition of the 
YZ$\beta$ shock-capturing operator suppresses the 
oscillatory overshoot, yielding a monotone profile; however, 
the SUPG-YZ$\beta$ solution introduces a slight smearing 
of the layer, as the isotropic artificial diffusion broadens 
the transition region beyond its analytical extent. The 
hybrid PINN correction recovers a sharper and more accurate 
representation of the boundary layer, closely tracking the 
analytical profile. This improvement is attributable to the 
physics-informed refinement: the network, having learned the 
global solution structure from the SUPG-YZ$\beta$ reference 
data, enforces the governing PDE residual via automatic 
differentiation, which acts as a regularizer that sharpens 
the layer without reintroducing spurious oscillations.

\begin{figure}[h!]
    \centering
    \begin{subfigure}[b]{0.49\linewidth}
        \includegraphics[width=\linewidth]{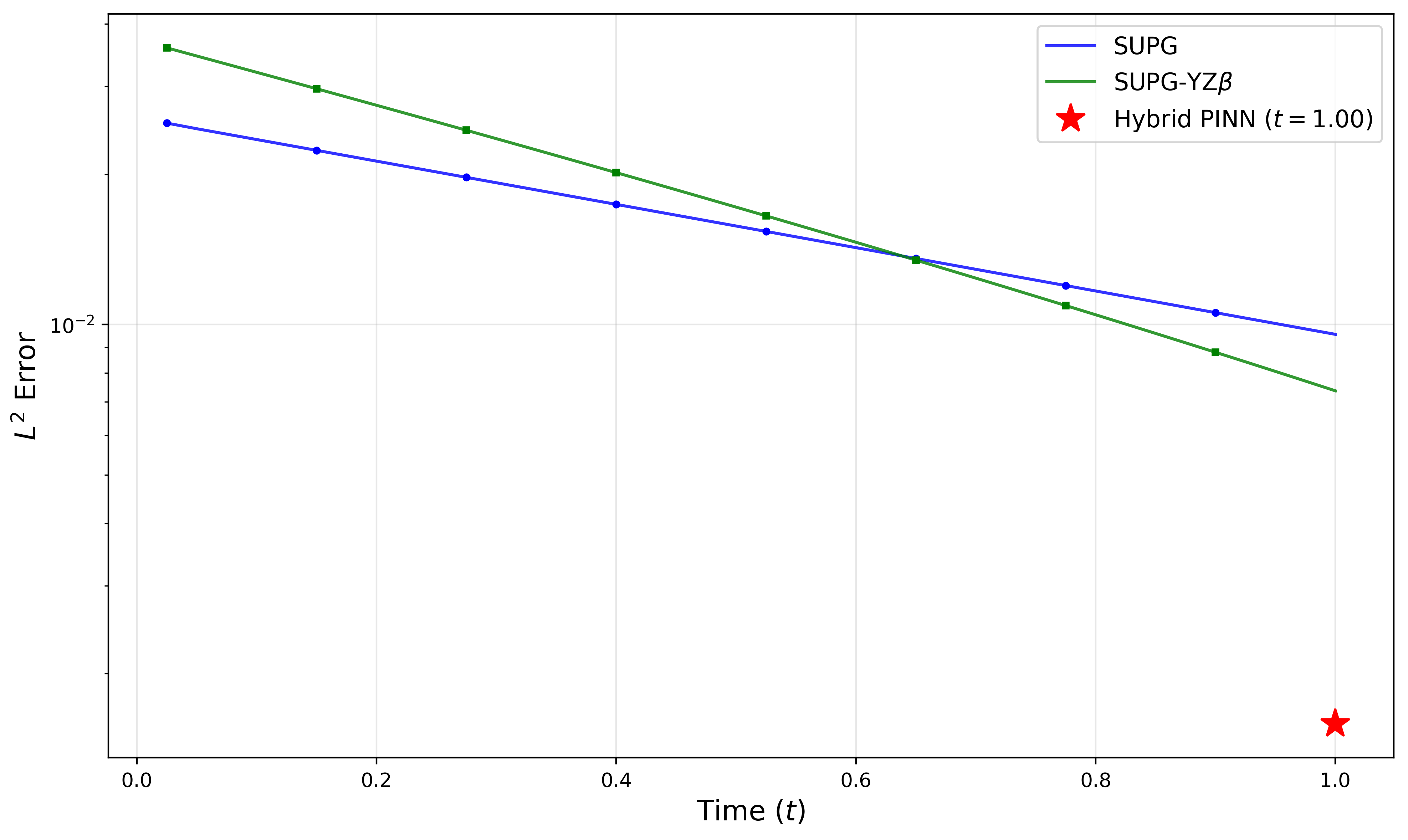}
        \caption{}
        \label{fig:fig2a}
    \end{subfigure}
    \hfill
    \begin{subfigure}[b]{0.49\linewidth}
        \includegraphics[width=\linewidth]{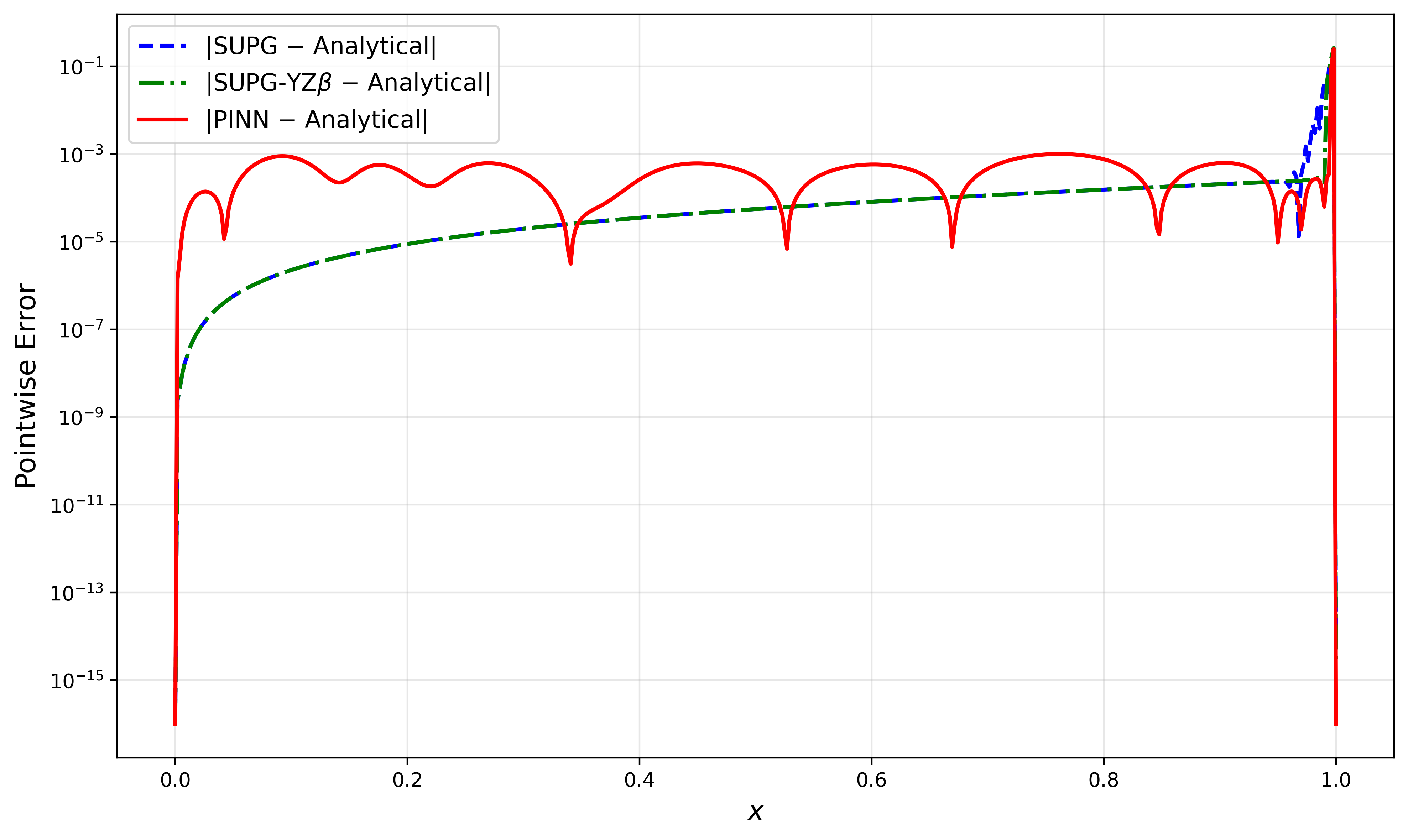}
        \caption{}
        \label{fig:fig2b}
    \end{subfigure}
    \caption{Error analysis for \textit{Example~1} with 
    	$\varepsilon = 10^{-4}$: 
    	(\subref{fig:fig2a})~evolution of the $L^2$ error 
    	over the time interval $[0, 1]$ for SUPG and 
    	SUPG-YZ$\beta$, with the hybrid PINN error at the 
    	terminal time $t_\text{f} = 1$ indicated by the star marker; 
    	(\subref{fig:fig2b})~pointwise absolute error at $t_\text{f} = 1$ 
    	across the spatial domain.}
    \label{fig:2_comparison_fem}
\end{figure}

The error analysis in Figure~\ref{fig:2_comparison_fem} 
provides a quantitative assessment of the three methods. 
Figure~\ref{fig:fig2a} displays the temporal evolution of 
the $L^2$ error norm. Both FEM solvers exhibit a 
monotonically decreasing error over the simulation interval, 
with the SUPG-YZ$\beta$ formulation consistently 
outperforming the SUPG-only discretization owing to the 
additional shock-capturing dissipation. Notably, the hybrid 
PINN correction at the terminal time $t_\text{f} = 1$ (red star) 
achieves an $L^2$ error approximately two orders of 
magnitude below the SUPG-YZ$\beta$ solution, demonstrating 
the substantial accuracy gain afforded by the 
physics-informed post-processing step. The spatial distribution of the pointwise error at $t_\text{f} = 1$, 
shown in Figure~\ref{fig:fig2b}, reveals further structural 
differences among the methods. The SUPG solution maintains an 
error of $\mathcal{O}(10^{-5})$ in the smooth interior but 
exhibits a sharp spike near $x = 1$, where the boundary 
layer oscillations produce errors exceeding 
$\mathcal{O}(10^{-1})$. The SUPG-YZ$\beta$ solution 
distributes the error more uniformly at 
$\mathcal{O}(10^{-5})$ across the domain, effectively 
eliminating the boundary layer spike at the cost of slightly 
elevated error in the smooth region due to the additional 
artificial diffusion. The hybrid PINN error profile, by 
contrast, displays a qualitatively different character: the 
pointwise error oscillates between $\mathcal{O}(10^{-7})$ 
and $\mathcal{O}(10^{-3})$, with several localized dips 
approaching machine precision. This non-uniform error 
distribution reflects the spectral approximation properties 
of the Fourier feature embedding, which allocates 
representational capacity across multiple frequency scales 
rather than distributing the error uniformly.

\begin{figure}[h!]
	\centering
	\begin{subfigure}[b]{0.49\linewidth}
		\includegraphics[width=\linewidth]{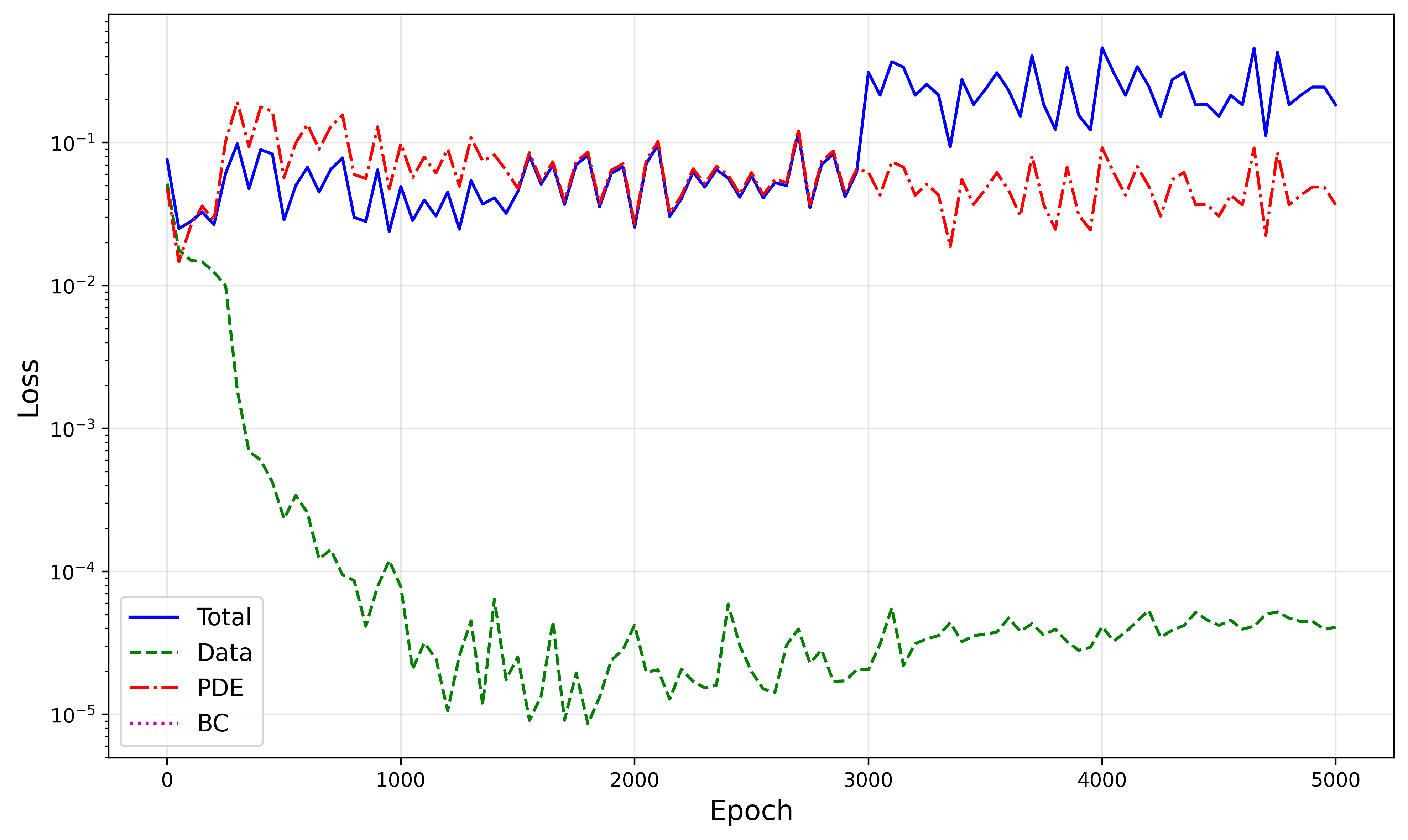}
		\caption{}
		\label{fig:fig3a}
	\end{subfigure}
	\hfill
	\begin{subfigure}[b]{0.49\linewidth}
		\includegraphics[width=\linewidth]{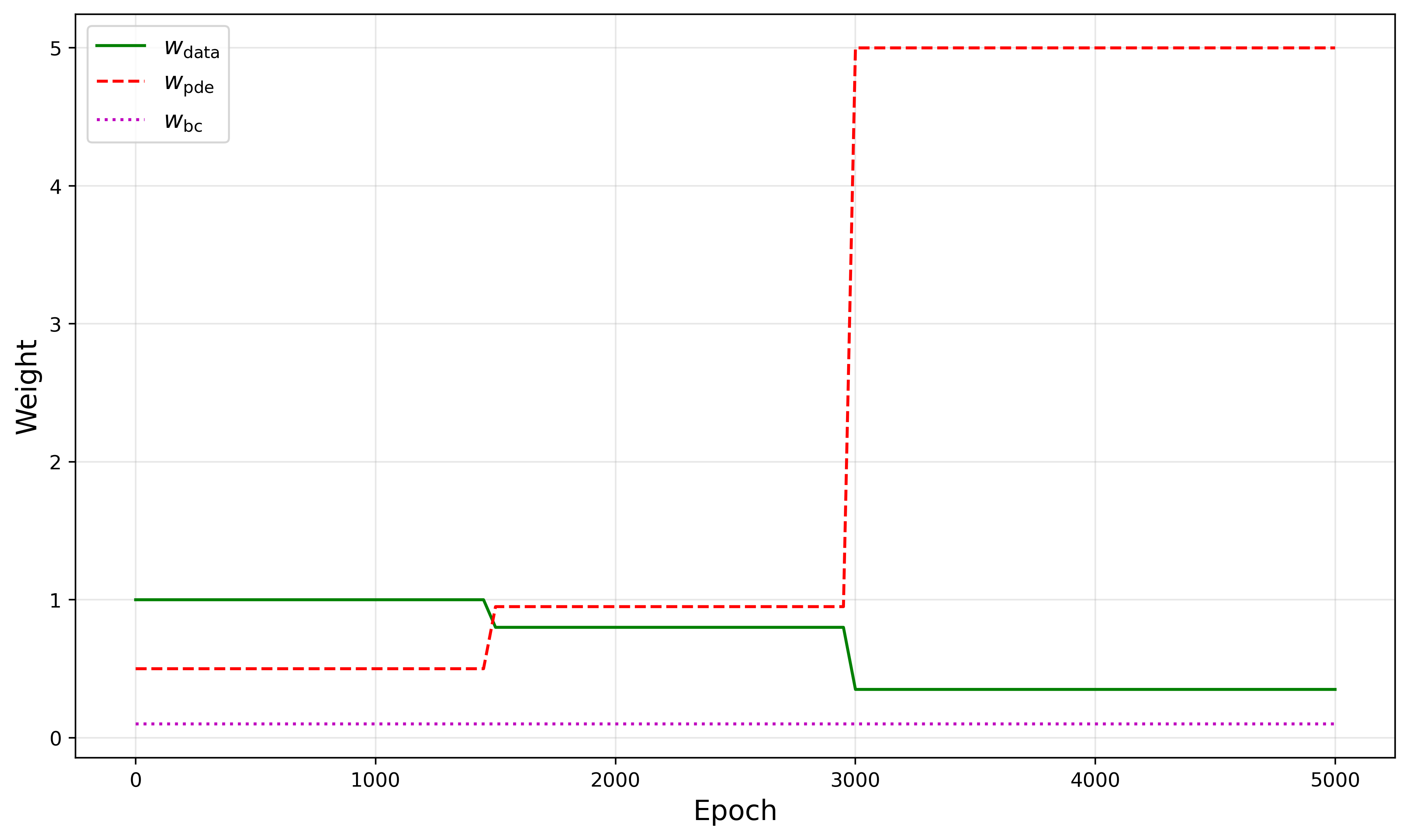}
		\caption{}
		\label{fig:fig3b}
	\end{subfigure}
	\caption{Training diagnostics for \textit{Example~1} with 
		$\varepsilon = 10^{-4}$: 
		(\subref{fig:fig3a})~evolution of the individual loss 
		components (total, data, PDE, and BC) over $5000$ 
		training epochs on a logarithmic scale; 
		(\subref{fig:fig3b})~adaptive weight schedule 
		$w_{\text{data}}$, $w_{\text{pde}}$, and 
		$w_{\text{bc}}$ governing the three-phase training 
		strategy.}
	\label{fig:3_training}
\end{figure}

The training dynamics of the hybrid PINN for \textit{Example~1} are 
illustrated in Figure~\ref{fig:3_training}. The loss 
evolution in Figure~\ref{fig:fig3a} reveals the distinct 
behavior of each component across the three training phases. 
During Phase~I (epochs $0$--$1499$), where data fidelity 
dominates with $w_{\text{data}} = 1.0$, the data loss 
decreases rapidly from $\mathcal{O}(10^{-2})$ to 
$\mathcal{O}(10^{-4})$, indicating that the network 
successfully learns the global solution structure from the 
SUPG-YZ$\beta$ reference snapshots. The PDE loss remains 
at $\mathcal{O}(10^{-1})$ during this phase, as the 
physics constraint receives only moderate weight 
($w_{\text{pde}} = 0.5$). The transition to Phase~II 
(epochs $1500$--$2999$) slightly increases the PDE 
influence to $w_{\text{pde}} = 0.95$ while reducing the 
data weight to $w_{\text{data}} = 0.8$, producing a 
gradual shift toward physics-driven learning. In Phase~III 
(epochs $3000$--$4999$), the PDE weight rises sharply to 
$w_{\text{pde}} = 5.0$ while the data weight drops to 
$w_{\text{data}} = 0.35$, as shown in the weight schedule 
of Figure~\ref{fig:fig3b}. The total loss increases 
accordingly, reflecting the now-dominant PDE residual 
contribution; however, the data loss continues to decrease 
steadily, confirming that physics-informed refinement 
does not compromise the learned solution fidelity. The 
boundary condition loss remains identically zero throughout 
training, as the lift-based distance function 
$d(x) = x(1-x)$ enforces homogeneous Dirichlet conditions 
exactly by construction (see~Remark~15), rendering the 
weight $w_{\text{bc}} = 0.1$ inconsequential.

\subsection{Example 2.}
Consider the following 2D convection-diffusion-reaction equation (hump changing its height)~\cite{Giere2015, John2008b}:
\begin{equation}
\begin{cases}
\displaystyle \frac{\partial u}{\partial t} - \varepsilon \Delta u + \mathbf{b} \cdot \nabla u + cu = f, & (x_1,x_2) \in \Omega, \quad t \in (0,t_{\text{f}}], \\[8pt]
u(0,x_1,x_2) = u_0(x_1,x_2), & (x_1,x_2) \in \Omega, \\[6pt]
u = 0, & (x_1,x_2) \in \partial\Omega, \quad t \in [0,t_{\text{f}}],
\end{cases}
\end{equation}
where $\Omega = (0,1)^2$, $t_{\text{f}} = 0.5$, $\varepsilon = 10^{-6}$, $\mathbf{b} = (2,3)^\top$, and $c = 1$. The forcing term $f$, the initial condition $u_0$, and the boundary conditions are chosen such that the exact solution
\begin{equation}
u(t,x_1, x_2) = 16 \sin(\pi t) \, x_{1}(1-x_{1}) \, x_2(1-x_2) \left[ \frac{1}{2} + \frac{\arctan\left( 2\varepsilon^{-1/2} \left( 0.25^2 - (x_1-0.5)^2 - (x_2-0.5)^2 \right) \right)}{\pi} \right]
\end{equation}
is satisfied. 

This problem was also considered in~\cite{Giere2015, John2008b}. In~\cite{John2008b}, the  authors compare several shock-capturing mechanisms and show  that they all fail for this particular type of problem by 
producing spurious stripes around the hump.

For the finite element setting, we choose 
$\varepsilon = 10^{-6}$, $n_{\text{el}} = 8,192$ 
(structured triangulation of a $64 \times 64$ grid), 
$\varDelta t = 0.0025$ ($N_t = 200$), and Y$= 0.7$. The PINN 
architecture employs $n_h = 128$, $n_r = 8$, $n_\text{F} = 24$, 
and $\sigma = 4.0$. Training uses the last $K_s = 5$ 
temporal snapshots near $t_\text{f} = 0.5$, with batch size 
$16000$, initial learning rate 
$\alpha = 8 \times 10^{-5} \sqrt{N_{\text{batch}}/256} 
\approx 6.3 \times 10^{-4}$ (square-root scaling), and 
gradient clipping threshold $g_{\max} = 1.0$. The PDE 
residual is evaluated at $N_{\text{pde}} = 256$ randomly 
sampled interior collocation points per mini-batch. The 
phase boundaries and adaptive weights are reported in 
Table~\ref{tab:phase_example_three}. 

\begin{table}[htb]
\centering
\caption{Phase boundaries and adaptive weights used for solving \textit{Example 2}.}
\label{tab:phase_example_three}
\begin{tabular}{llccc}
\toprule
\text{Phase} & \text{Epochs} & $w_{\mathrm{data}}$ & $w_{\mathrm{pde}}$ & $w_{\mathrm{bc}}$  \\
\midrule
I   & 0--999    & 1.0  & 0.05 & 0.01 \\
II  & 1000--1499 & 0.5  & 0.1  & 0.05 \\
III & 1500--4999 & 0.1  & 0.5  & 0.1  \\
\bottomrule
\end{tabular}
\end{table}

\begin{figure}[h!]
	\centering
	\begin{subfigure}[b]{0.325\linewidth}
		\includegraphics[width=\linewidth]{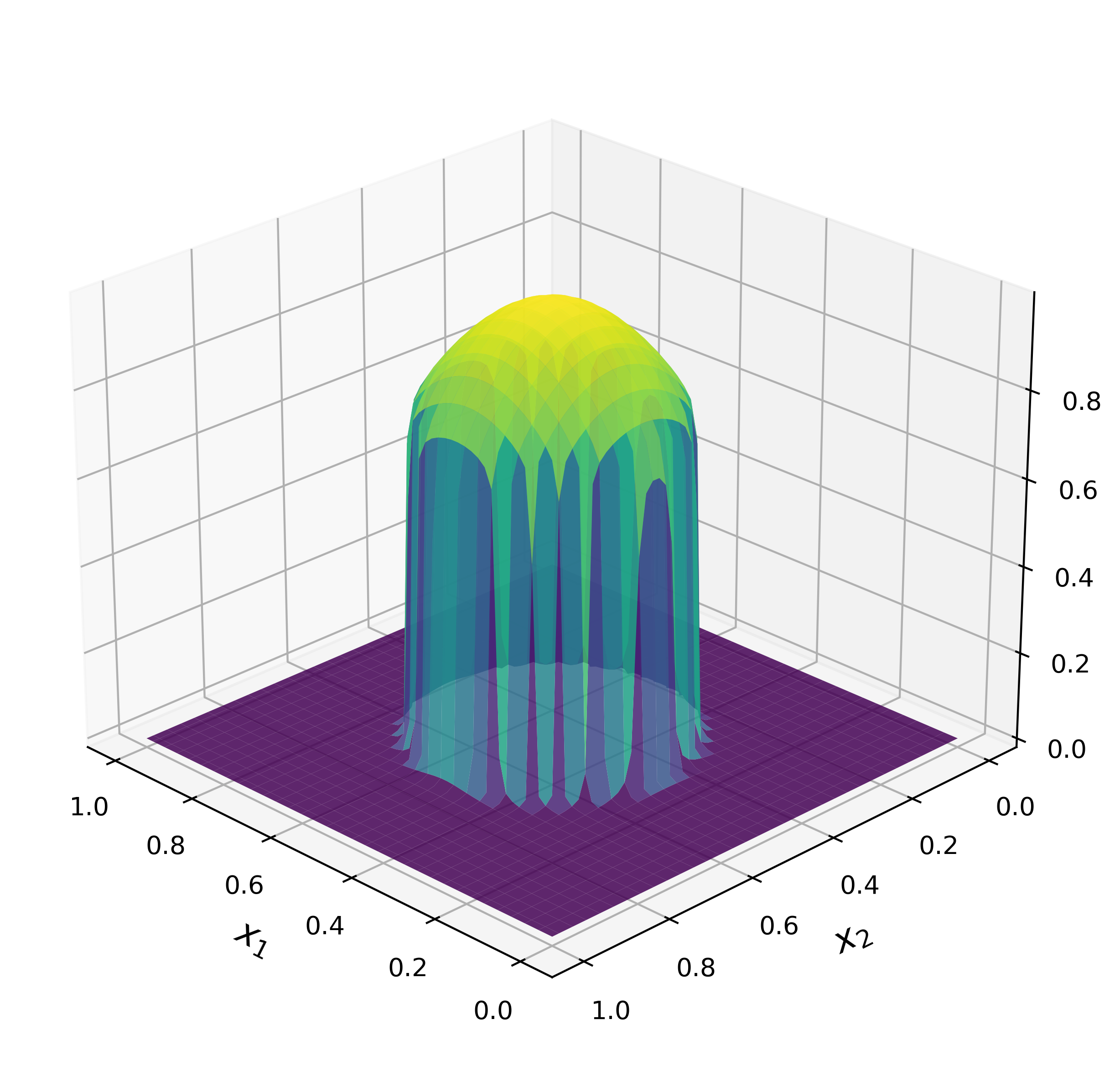}
		\caption{}
		\label{fig:fig4a}
	\end{subfigure}
	\hfill
	\begin{subfigure}[b]{0.325\linewidth}
		\includegraphics[width=\linewidth]{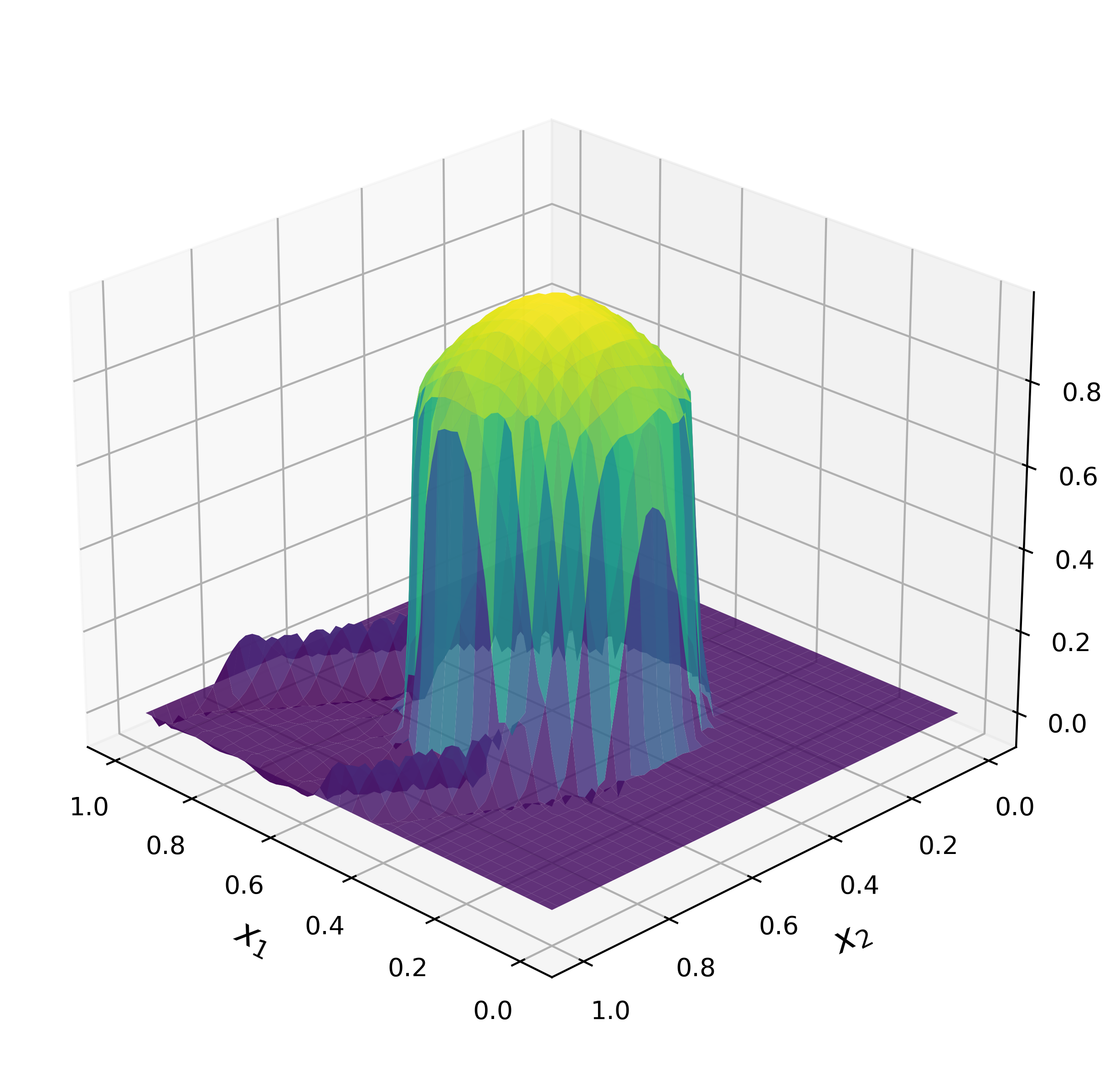}
		\caption{}
		\label{fig:fig4b}
	\end{subfigure}
	\hfill
	\begin{subfigure}[b]{0.325\linewidth}
		\includegraphics[width=\linewidth]{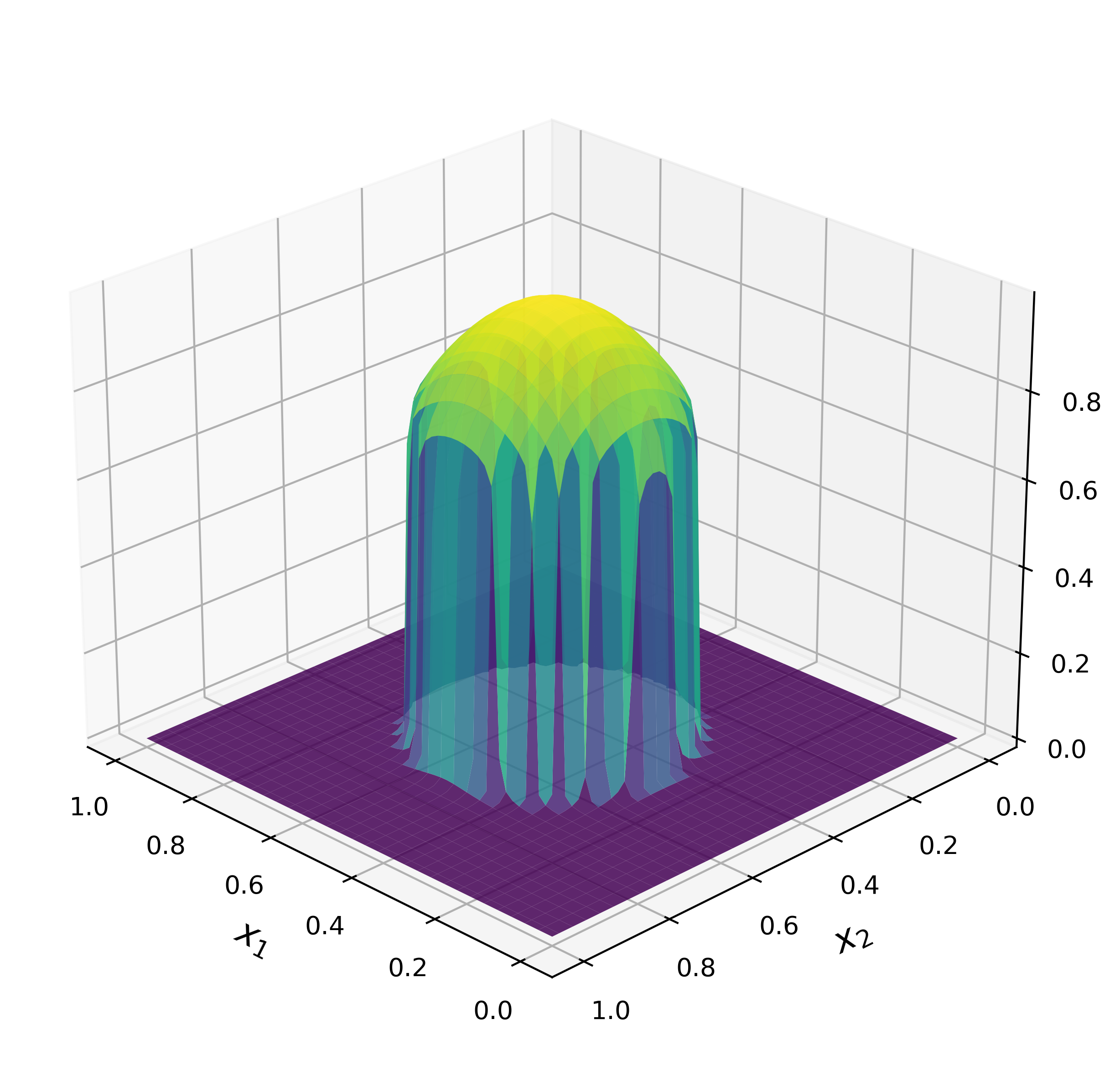}
		\caption{}
		\label{fig:fig4c}
	\end{subfigure}
	\caption{Three-dimensional surface plots of the solution 
		to \textit{Example~2} at $t = 0.5$ with 
		$\varepsilon = 10^{-6}$: 
		(\subref{fig:fig4a})~analytical solution; 
		(\subref{fig:fig4b})~SUPG-YZ$\beta$ finite element 
		approximation; 
		(\subref{fig:fig4c})~hybrid PINN correction.}
	\label{fig:4_surface}
\end{figure}

Figure~\ref{fig:4_surface} displays the three-dimensional 
solution profiles at $t = 0.5$ for \textit{Example~2}. The 
analytical solution (Figure~\ref{fig:fig4a}) features a 
smooth, compactly supported hump centered at 
$(x_1, x_2) = (0.5, 0.5)$ whose amplitude is modulated 
by the factor $16\sin(\pi t)$, reaching its peak value 
at the selected time level. The hump possesses extremely 
steep flanks governed by the arctangent layer of width 
$\mathcal{O}(\sqrt{\varepsilon}) = \mathcal{O}(10^{-3})$, 
making this a particularly challenging test case. As noted 
in~\cite{John2008b}, all classical shock-capturing 
mechanisms fail for this problem type, and the 
SUPG-YZ$\beta$ approximation in Figure~\ref{fig:fig4b} 
confirms this observation: pronounced spurious stripes 
appear around the base and flanks of the hump, together 
with localized oscillations in the flat region where the 
solution should vanish identically. These artifacts arise 
because the isotropic shock-capturing diffusion interacts 
adversely with the circular symmetry of the internal layer, 
generating crosswind distortions that propagate along the 
convection direction $\mathbf{b} = (2, 3)^{\top}$. The 
hybrid PINN correction (Figure~\ref{fig:fig4c}) 
substantially eliminates these spurious features, 
recovering a surface profile that closely matches the 
analytical solution. The hump shape, peak amplitude, and 
steep flanks are faithfully reproduced, while the 
oscillatory stripes visible in the FEM approximation are 
effectively suppressed. This result demonstrates the 
capacity of the physics-informed post-processing to 
correct structural artifacts that persist even after 
stabilization and shock-capturing have been applied.

\begin{figure}[h!]
	\centering
	\begin{subfigure}[b]{0.495\linewidth}
		\includegraphics[width=\linewidth]{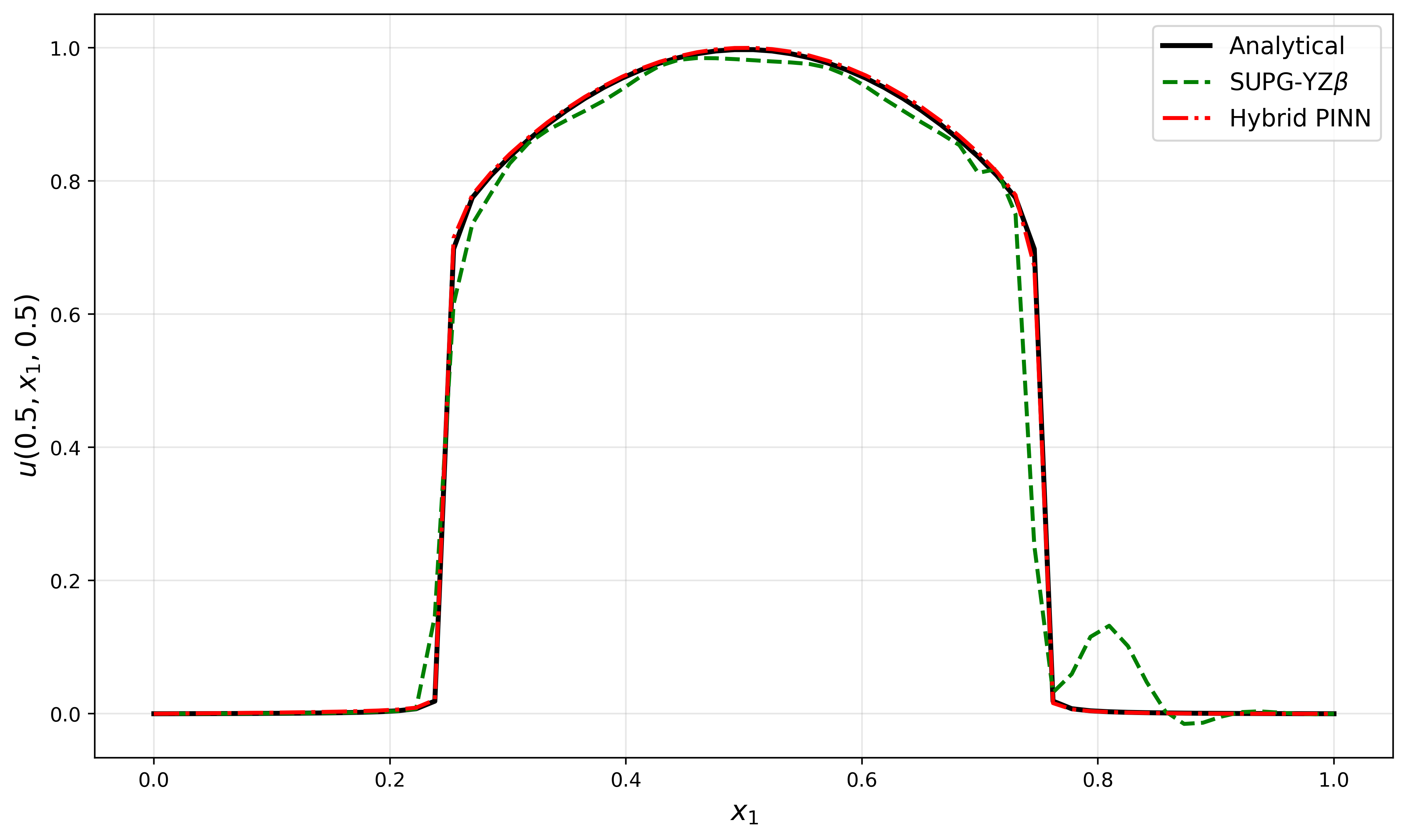}
		\caption{}
		\label{fig:fig5a}
	\end{subfigure}
	\hfill
	\begin{subfigure}[b]{0.495\linewidth}
		\includegraphics[width=\linewidth]{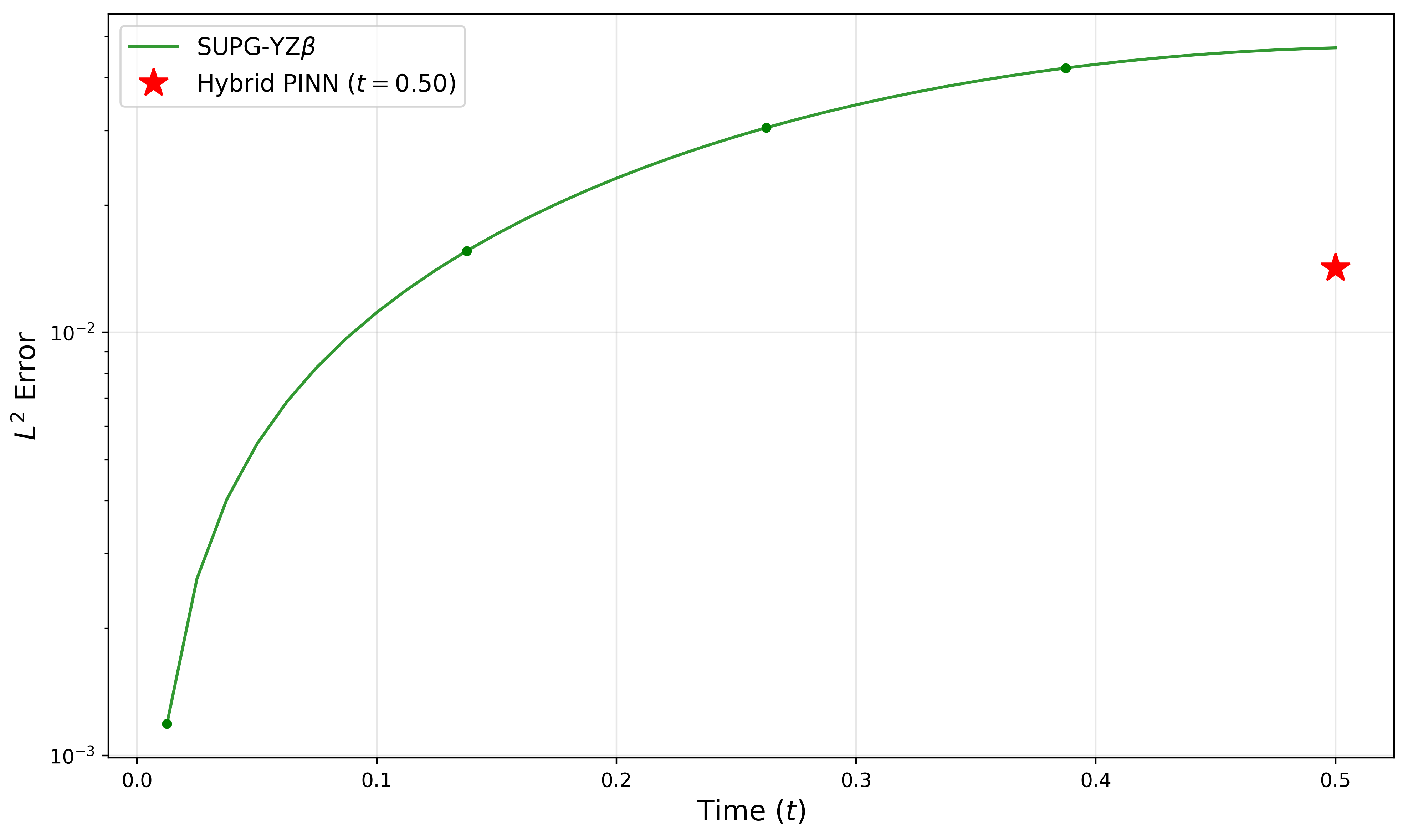}
		\caption{}
		\label{fig:fig5b}
	\end{subfigure}
	\caption{Quantitative comparison for \textit{Example~2} 
		with $\varepsilon = 10^{-6}$: 
		(\subref{fig:fig5a})~cross-sectional profiles along 
		$x_2 = 0.5$ at $t = 0.5$, comparing the analytical 
		solution, SUPG-YZ$\beta$, and hybrid PINN; 
		(\subref{fig:fig5b})~evolution of the $L^2$ error over 
		the time interval $[0, 0.5]$ for SUPG-YZ$\beta$, with 
		the hybrid PINN error at the terminal time indicated 
		by the star marker.}
	\label{fig:5_cross_l2}
\end{figure}

A more quantitative assessment is provided in 
Figure~\ref{fig:5_cross_l2}. The cross-sectional profile 
along $x_2 = 0.5$ at $t = 0.5$ 
(Figure~\ref{fig:fig5a}) reveals the detailed structure of 
the internal layer. The analytical solution displays a 
sharp, symmetric hump with a flat plateau near unity and 
steep flanks transitioning to zero over a distance of 
$\mathcal{O}(\sqrt{\varepsilon})$. The SUPG-YZ$\beta$ 
approximation captures the overall shape but exhibits two 
characteristic deficiencies: the steep flanks are smeared 
relative to the exact profile, and a spurious secondary 
bump of amplitude ${\sim}0.13$ appears near $x_1 \approx 
0.8$, corresponding to the stripe artifact visible in the 
three-dimensional view of Figure~\ref{fig:fig4b}. The 
hybrid PINN correction eliminates this spurious feature 
entirely and recovers steeper flanks that track the 
analytical profile more faithfully, although a slight 
undershoot is observed at the downstream transition near 
$x_1 \approx 0.75$. The temporal evolution of the $L^2$ error norm in 
Figure~\ref{fig:fig5b} shows a monotonically increasing 
trend for the SUPG-YZ$\beta$ solution, rising from 
$\mathcal{O}(10^{-3})$ at early times to 
$\mathcal{O}(10^{-2})$ at $t = 0.5$. This growth reflects 
the progressive accumulation of the stripe artifacts as the 
hump amplitude increases with $\sin(\pi t)$. The hybrid 
PINN error at the terminal time (red star) lies below the 
SUPG-YZ$\beta$ error curve, confirming that the 
physics-informed correction yields a measurable improvement 
in global accuracy despite operating only on the final 
temporal snapshots.

\begin{figure}[h!]
	\centering
	\begin{subfigure}[b]{0.495\linewidth}
		\includegraphics[width=\linewidth]{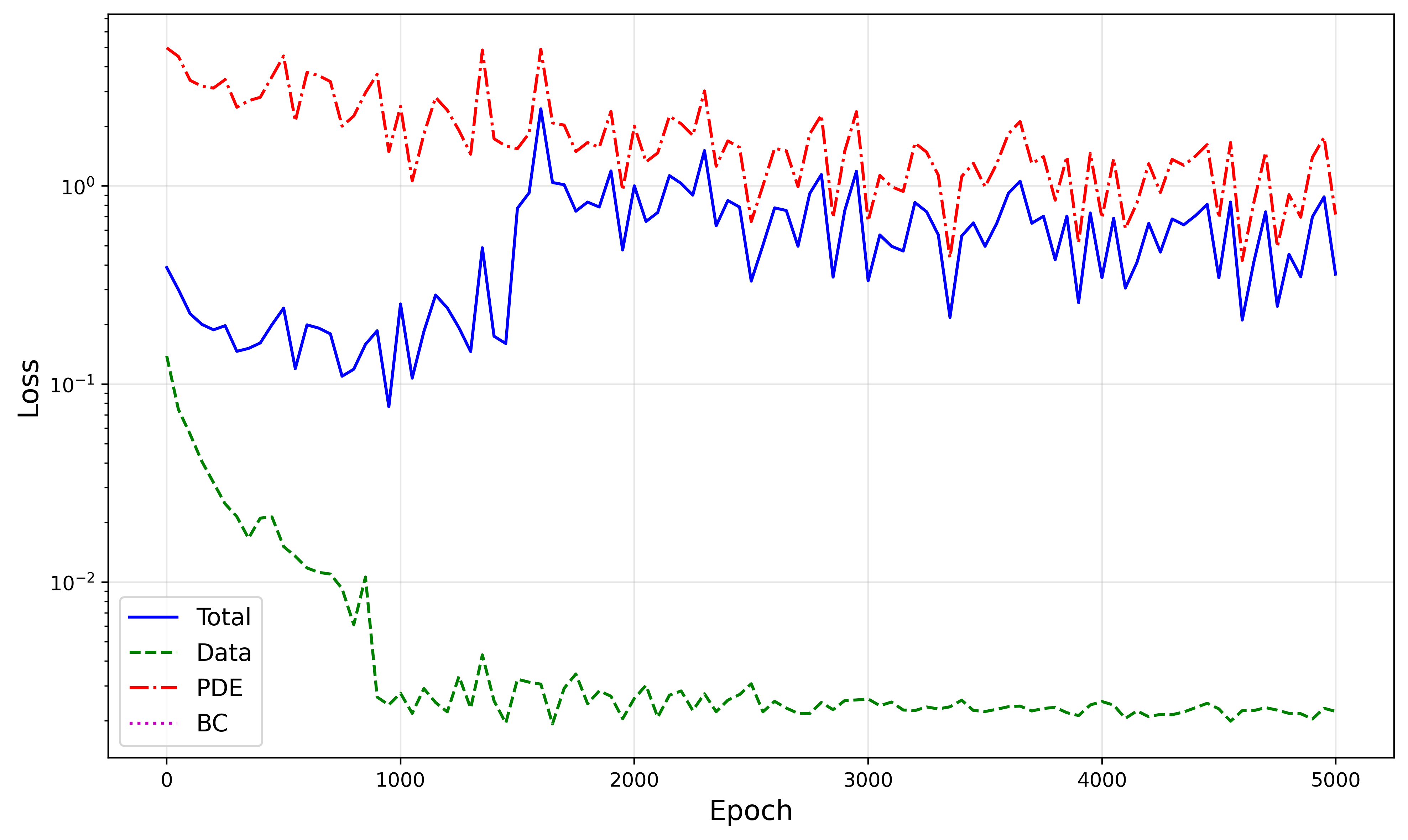}
		\caption{}
		\label{fig:fig6a}
	\end{subfigure}
	\hfill
	\begin{subfigure}[b]{0.495\linewidth}
		\includegraphics[width=\linewidth]{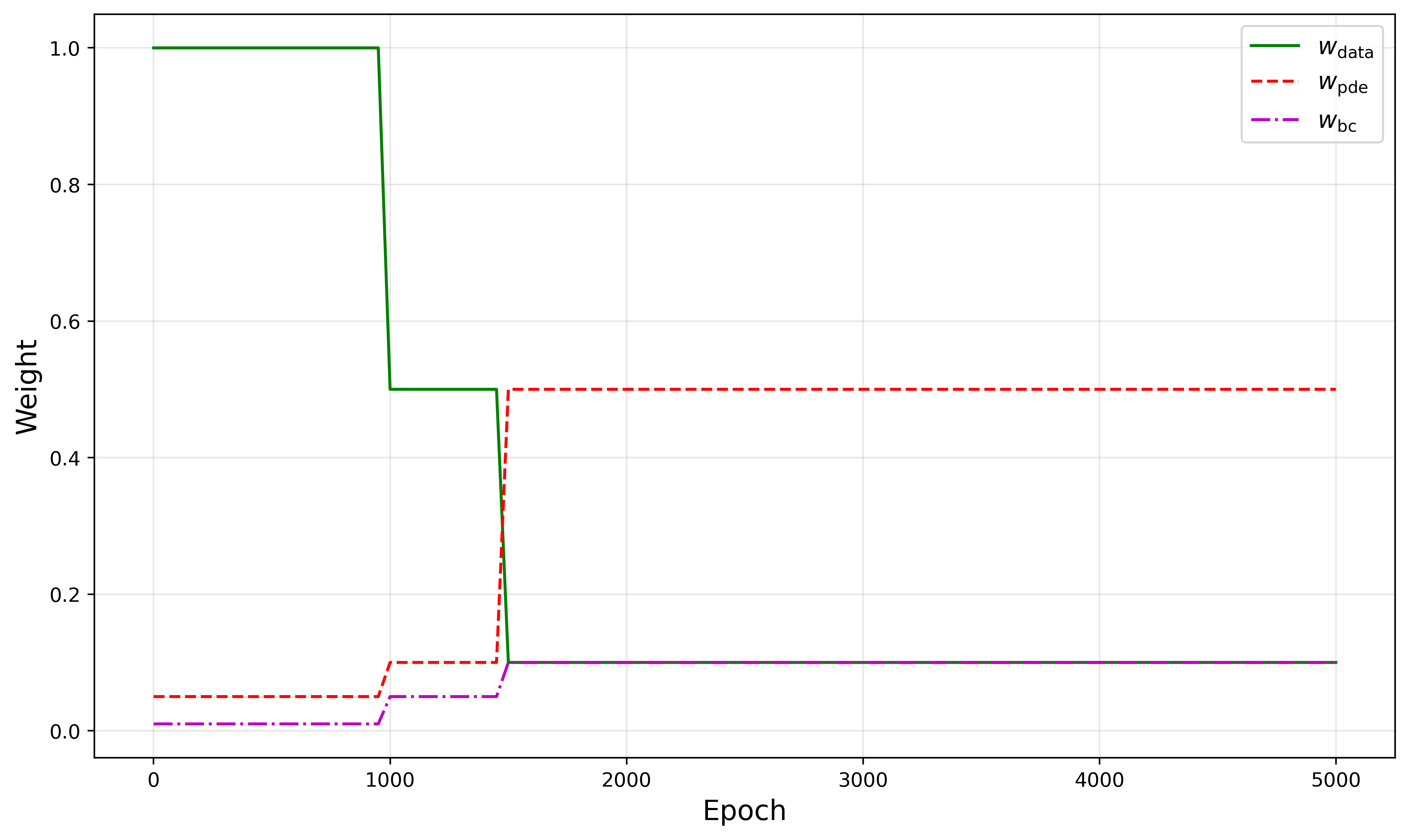}
		\caption{}
		\label{fig:fig6b}
	\end{subfigure}
	\caption{Training diagnostics for \textit{Example~2} 
		with $\varepsilon = 10^{-6}$: 
		(\subref{fig:fig6a})~evolution of the individual loss 
		components (total, data, PDE, and BC) over $5000$ 
		training epochs on a logarithmic scale; 
		(\subref{fig:fig6b})~adaptive weight schedule 
		$w_{\text{data}}$, $w_{\text{pde}}$, and 
		$w_{\text{bc}}$ governing the three-phase training 
		strategy.}
	\label{fig:6_training}
\end{figure}

The training diagnostics for \textit{Example~2} are presented in 
Figure~\ref{fig:6_training}. Compared with the 1D case 
of \textit{Example~1} (see\ Figure~\ref{fig:3_training}), 
the loss landscape exhibits markedly different behavior, 
reflecting the increased complexity of the two-dimensional 
hump problem. As shown in Figure~\ref{fig:fig6a}, the 
data loss decreases steadily from 
$\mathcal{O}(10^{-1})$ to $\mathcal{O}(10^{-2})$ 
during Phase~I (epochs $0$--$999$), where the data weight 
dominates with $w_{\text{data}} = 1.0$. However, the PDE 
loss remains persistently elevated at 
$\mathcal{O}(1)$--$\mathcal{O}(10)$ throughout training, 
substantially higher than in \textit{Example~1}. This elevated 
residual is expected: the steep internal layer of width 
$\mathcal{O}(\sqrt{\varepsilon}) = \mathcal{O}(10^{-3})$ 
generates large pointwise PDE residuals that resist 
minimization, particularly in the transition regions 
around the hump flanks.

The weight schedule in Figure~\ref{fig:fig6b} reveals the 
deliberately conservative physics enforcement adopted for 
this problem. Upon transitioning from Phase~I to Phase~II 
at epoch $1000$, the data weight drops from $1.0$ to 
$0.5$ while the PDE weight increases modestly from $0.05$ 
to $0.1$. In Phase~III (epochs $1500$--$4999$), the 
weights stabilize at $w_{\text{data}} = 0.1$ and 
$w_{\text{pde}} = 0.5$---notably milder than the 
aggressive physics-dominant configuration used in 
\textit{Example~1} ($w_{\text{pde}} = 5.0$). This restraint is 
deliberate: excessive PDE weight enforcement for the hump 
problem risks reintroducing the spurious stripe artifacts 
that the PINN correction is designed to eliminate. The 
boundary condition loss remains at low levels throughout 
training, consistent with the lift-based enforcement via 
the distance function 
$d(\mathbf{x}) = x_1(1-x_1)\,x_2(1-x_2)$.

\subsection{Example 3.}
Consider the following 2D convection-diffusion-reaction equation (traveling wave)~\cite{Giere2015}:
\begin{equation}
\begin{cases}
\displaystyle \frac{\partial u}{\partial t} - \varepsilon \Delta u + \mathbf{b} \cdot \nabla u + cu = f, & (x_1,x_2) \in \Omega, \quad t \in (0,t_{\text{f}}], \\[8pt]
u = 0, & (x_1,x_2) \in \partial\Omega, \quad t \in [0,t_{\text{f}}], \\[6pt]
u(0,x_1,x_2) = u_0(x_1,x_2), & (x_1,x_2) \in \Omega,
\end{cases}
\end{equation}
where $\Omega = (0,1)^2$, $t_{\text{f}} = 1$, $\varepsilon = 10^{-8}$, $\mathbf{b} = (\cos(\pi/3), \sin(\pi/3))^\top$, and $c = 1$. The forcing term $f$, the initial condition $u_0$, and the boundary conditions are chosen such that the exact solution
\begin{equation}
u(t,x_1,x_2) = 0.5 \sin(\pi x_1) \sin(\pi x_2) \left[ \tanh\left( \frac{x_1 + x_2 - t - 0.5}{\sqrt{\varepsilon}} \right) + 1 \right]
\end{equation}
is satisfied. The solution possesses a moving internal layer of width $\mathcal{O}(\sqrt{\varepsilon})$.

For the finite element setting, we set 
$n_{\text{el}} = 8192$ (structured triangulation of a 
$64 \times 64$ grid), $\varDelta t = 0.001$ ($N_t = 1000$), 
and Y$= 0.25$. The moving internal layer of width 
$\mathcal{O}(\sqrt{\varepsilon}) = \mathcal{O}(10^{-4})$ 
necessitates a finer temporal resolution compared to 
\textit{Example~2} to accurately track the propagating front. The 
PINN architecture employs a smaller network with 
$n_h = 96$, $n_r = 6$, $n_\text{F} = 16$, and $\sigma = 4.0$, 
reflecting the lower-dimensional structure of the traveling 
wave solution compared to the hump problem. Training uses 
the last $K_s = 5$ temporal snapshots near $t_\text{f} = 1.0$, 
with batch size $16000$, initial learning rate 
$\alpha \approx 6.3 \times 10^{-4}$ (square-root scaling), 
and gradient clipping threshold $g_{\max} = 1.0$. The PDE 
residual is evaluated at $N_{\text{pde}} = 384$ randomly 
sampled interior collocation points per mini-batch. Since 
the network architecture incorporates the lift-based 
distance function 
$d(\mathbf{x}) = x_1(1-x_1)\,x_2(1-x_2)$, the boundary 
loss is omitted entirely ($w_{\text{bc}} = 0$) throughout 
all training phases, as Dirichlet conditions are satisfied 
exactly by construction (see~Remark~15). The phase 
boundaries and adaptive weights are reported in 
Table~\ref{tab:phase_example_four}.

\begin{table}[htb]
\centering
\caption{Phase boundaries and adaptive weights used for solving \textit{Example 3}.}
\label{tab:phase_example_four}
\begin{tabular}{llccc}
\toprule
\text{Phase} & \text{Epochs} & $w_{\mathrm{data}}$ & $w_{\mathrm{pde}}$ & $w_{\mathrm{bc}}$ \\
\midrule
I   & 0--999     & 1.0 & 0.05 & 0.0 \\
II  & 1000--1499 & 0.5 & 0.1  & 0.0 \\
III & 1500--4999 & 0.1 & 0.5  & 0.0 \\
\bottomrule
\end{tabular}
\end{table}

\begin{figure}[h!]
	\centering
	\begin{subfigure}[b]{0.245\linewidth}
		\includegraphics[width=\linewidth]{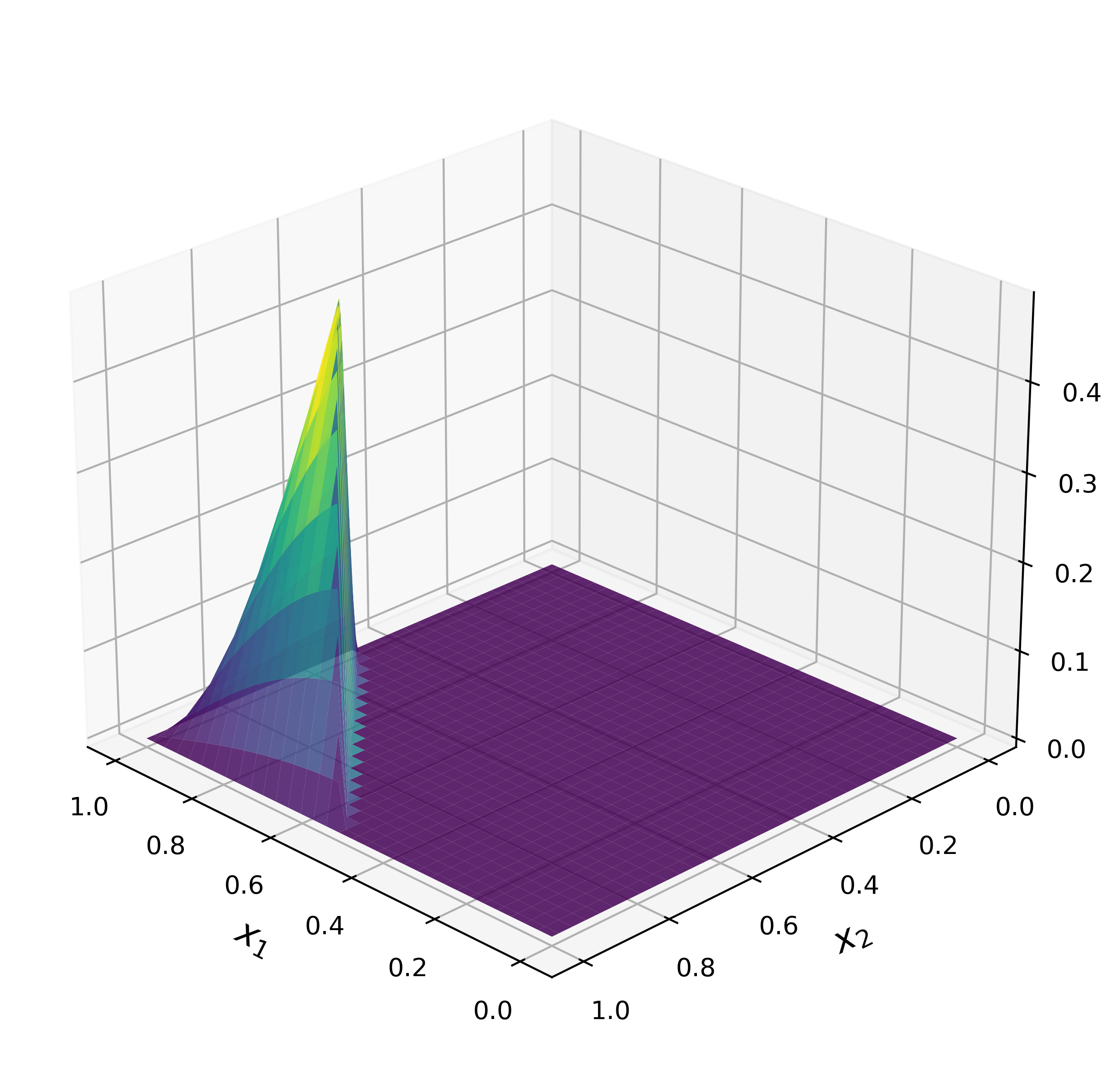}
		\caption{}
		\label{fig:fig7a}
	\end{subfigure}
	\hfill
	\begin{subfigure}[b]{0.245\linewidth}
		\includegraphics[width=\linewidth]{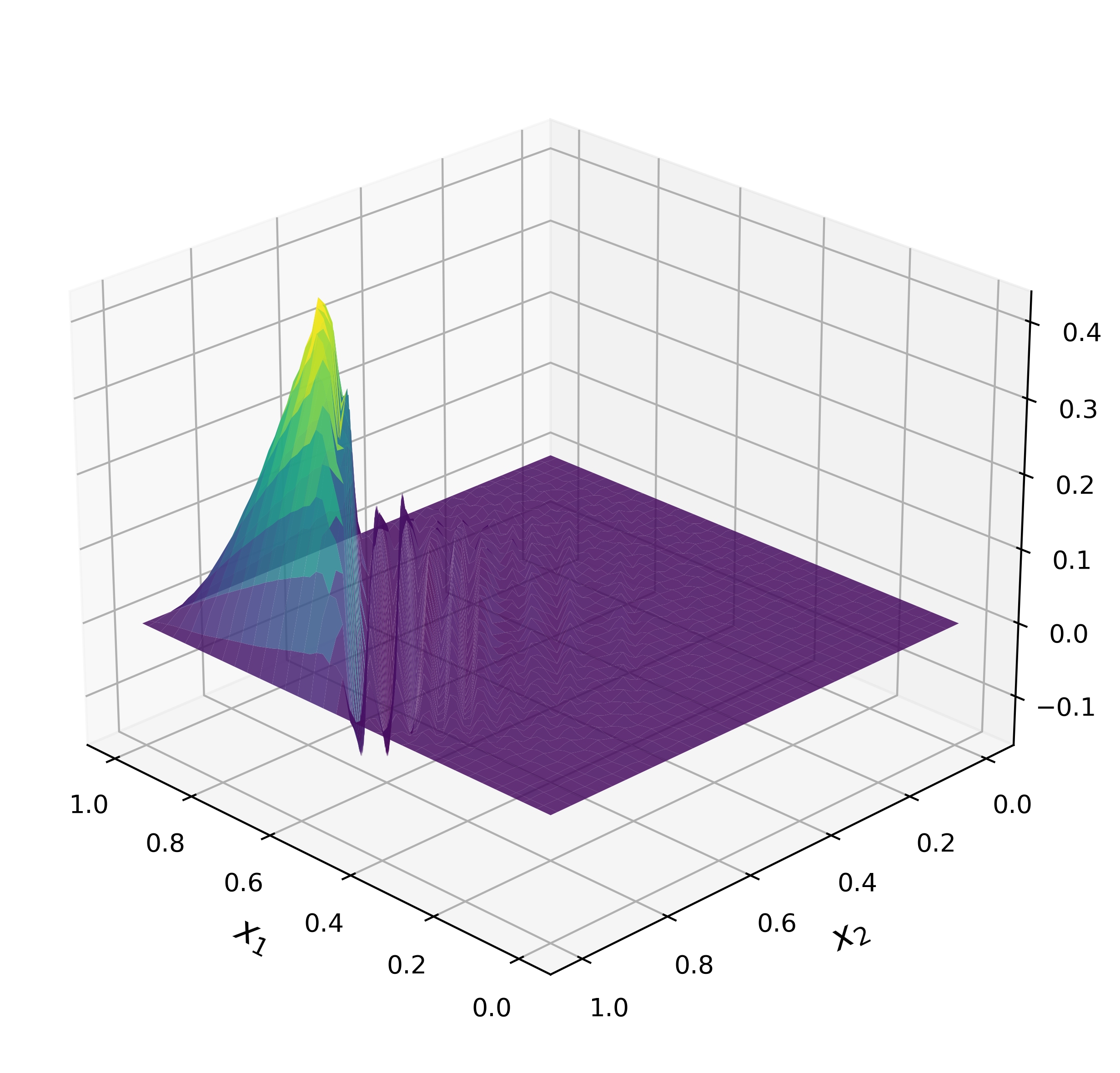}
		\caption{}
		\label{fig:fig7b}
	\end{subfigure}
	\hfill
	\begin{subfigure}[b]{0.245\linewidth}
		\includegraphics[width=\linewidth]{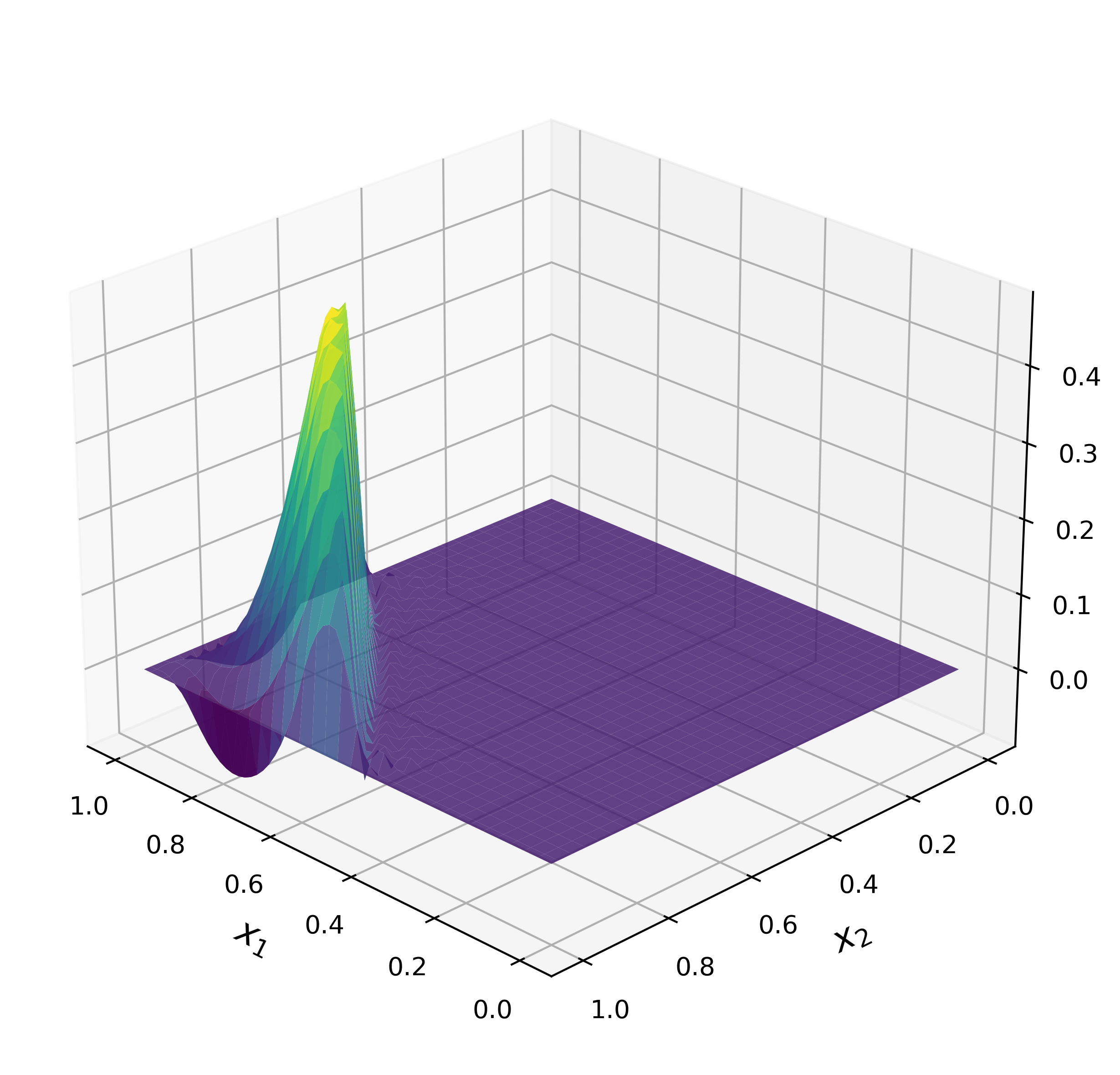}
		\caption{}
		\label{fig:fig7c}
	\end{subfigure}
	\hfill
	\begin{subfigure}[b]{0.245\linewidth}
		\includegraphics[width=\linewidth]{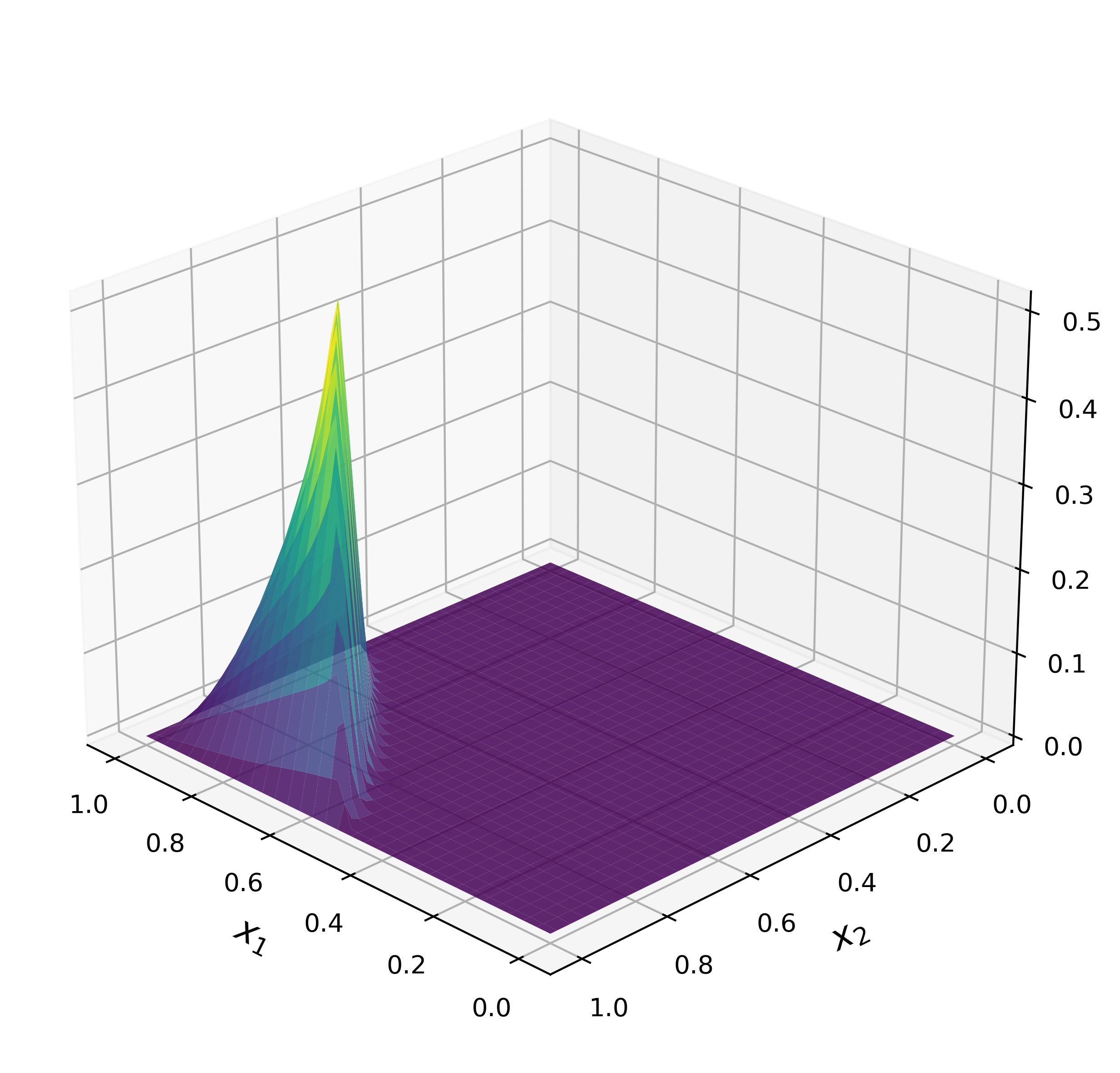}
		\caption{}
		\label{fig:fig7d}
	\end{subfigure}
	\caption{Three-dimensional surface plots of the solution 
		to \textit{Example~3} at $t_\text{f} = 1$ with 
		$\varepsilon = 10^{-8}$: 
		(\subref{fig:fig7a})~analytical solution; 
		(\subref{fig:fig7b})~SUPG finite element approximation; 
		(\subref{fig:fig7c})~SUPG-YZ$\beta$ finite element 
		approximation; 
		(\subref{fig:fig7d})~hybrid PINN correction.}
	\label{fig:7_surface}
\end{figure}

Figure~\ref{fig:7_surface} presents the three-dimensional 
solution profiles at the terminal time $t_\text{f} = 1$ for 
\textit{Example~3}. The analytical solution 
(Figure~\ref{fig:fig7a}) features a localized traveling 
wave peak resulting from the interaction of the sinusoidal 
envelope $0.5\sin(\pi x_1)\sin(\pi x_2)$ with the 
internal layer governed by the hyperbolic tangent profile 
of width $\mathcal{O}(\sqrt{\varepsilon}) = 
\mathcal{O}(10^{-4})$. This extremely thin moving front 
poses a severe challenge for all numerical methods. The SUPG solution (Figure~\ref{fig:fig7b}) fails 
dramatically: in the absence of shock-capturing, the 
unstabilized crosswind diffusion generates pronounced 
undershoots reaching approximately $-0.1$, together with 
widespread oscillatory artifacts that contaminate the 
entire domain. The addition of YZ$\beta$ shock-capturing 
(Figure~\ref{fig:fig7c}) eliminates the negative 
undershoots and suppresses the global oscillations, 
yielding a non-negative solution; however, the peak 
amplitude is visibly attenuated and the layer is 
broadened by the isotropic artificial diffusion, 
resulting in a smeared profile that underestimates 
the true solution magnitude. The hybrid PINN correction 
(Figure~\ref{fig:fig7d}) recovers both the sharp layer 
structure and the correct peak amplitude, producing a 
surface profile that closely matches the analytical 
solution. The peak is well-localized, the surrounding 
flat region remains free of spurious oscillations, and 
the transition from the elevated region to the zero 
plateau is captured with substantially higher fidelity 
than either FEM approximation. This example highlights 
the particular strength of the hybrid framework for 
problems involving extremely thin moving fronts, where 
even stabilized and shock-captured finite element 
solutions suffer from excessive numerical diffusion.

\begin{figure}[h!]
	\centering
	\begin{subfigure}[b]{0.245\linewidth}
		\includegraphics[width=\linewidth]{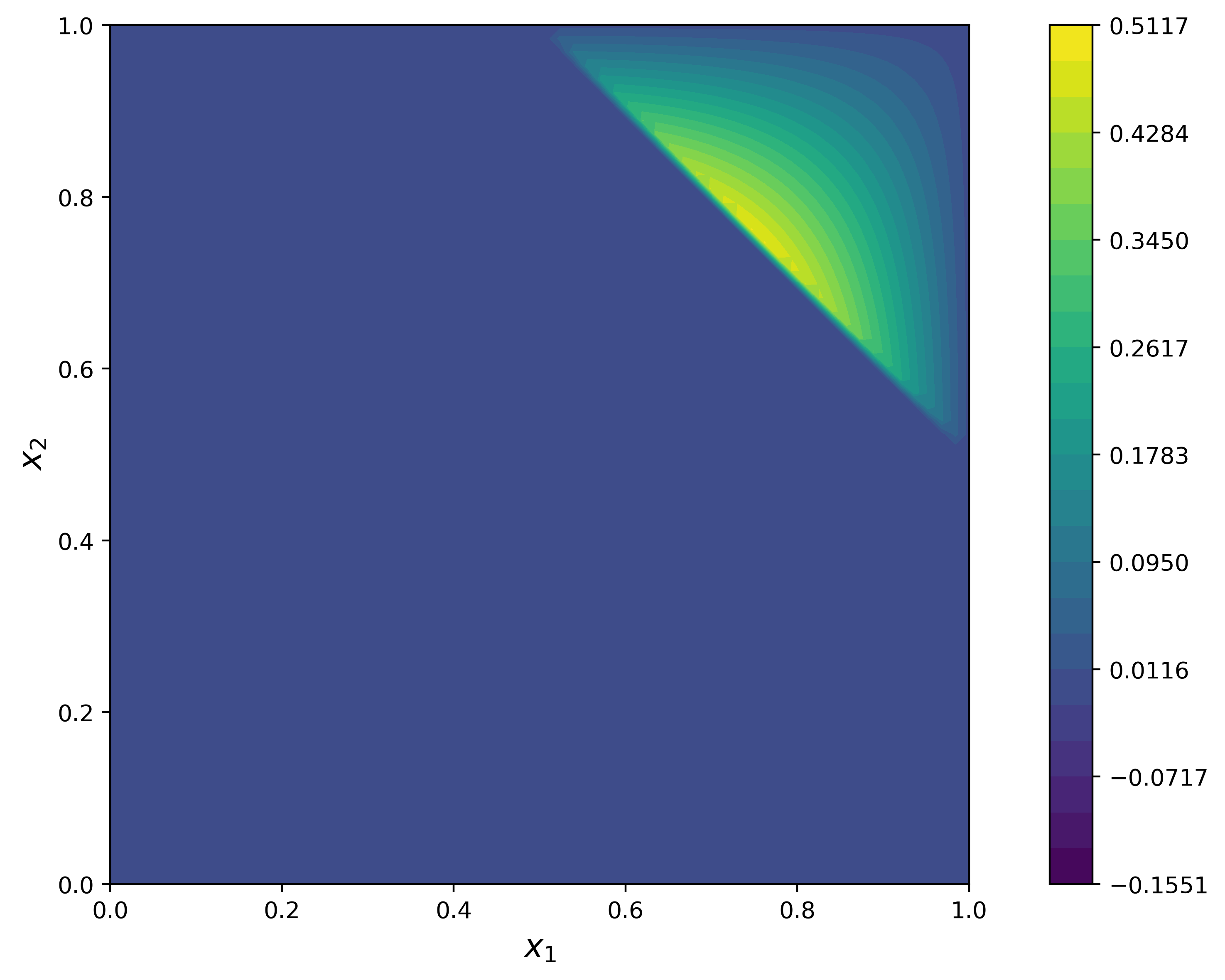}
		\caption{}
		\label{fig:fig8a}
	\end{subfigure}
	\hfill
	\begin{subfigure}[b]{0.245\linewidth}
		\includegraphics[width=\linewidth]{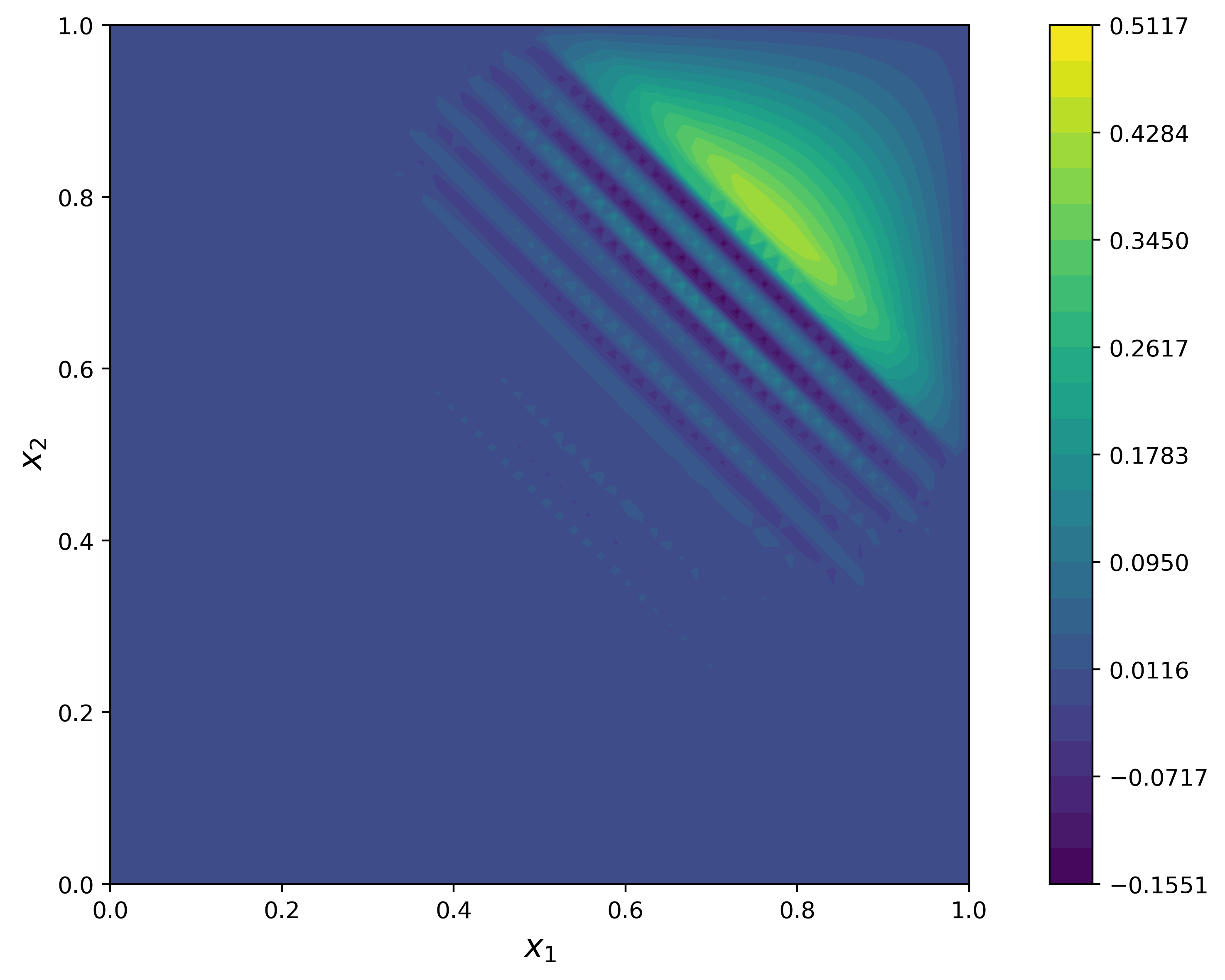}
		\caption{}
		\label{fig:fig8b}
	\end{subfigure}
	\hfill
	\begin{subfigure}[b]{0.245\linewidth}
		\includegraphics[width=\linewidth]{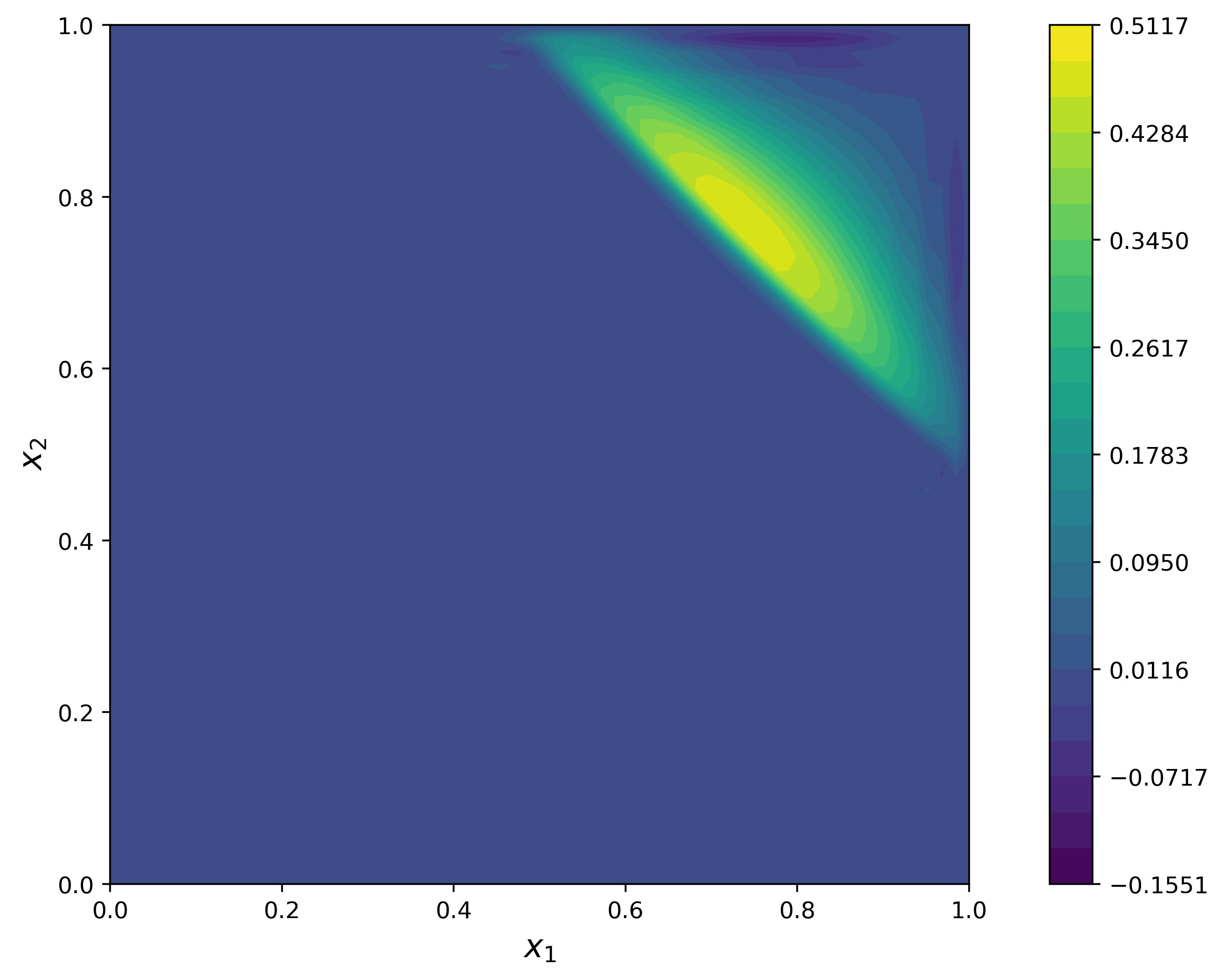}
		\caption{}
		\label{fig:fig8c}
	\end{subfigure}
	\hfill
	\begin{subfigure}[b]{0.245\linewidth}
		\includegraphics[width=\linewidth]{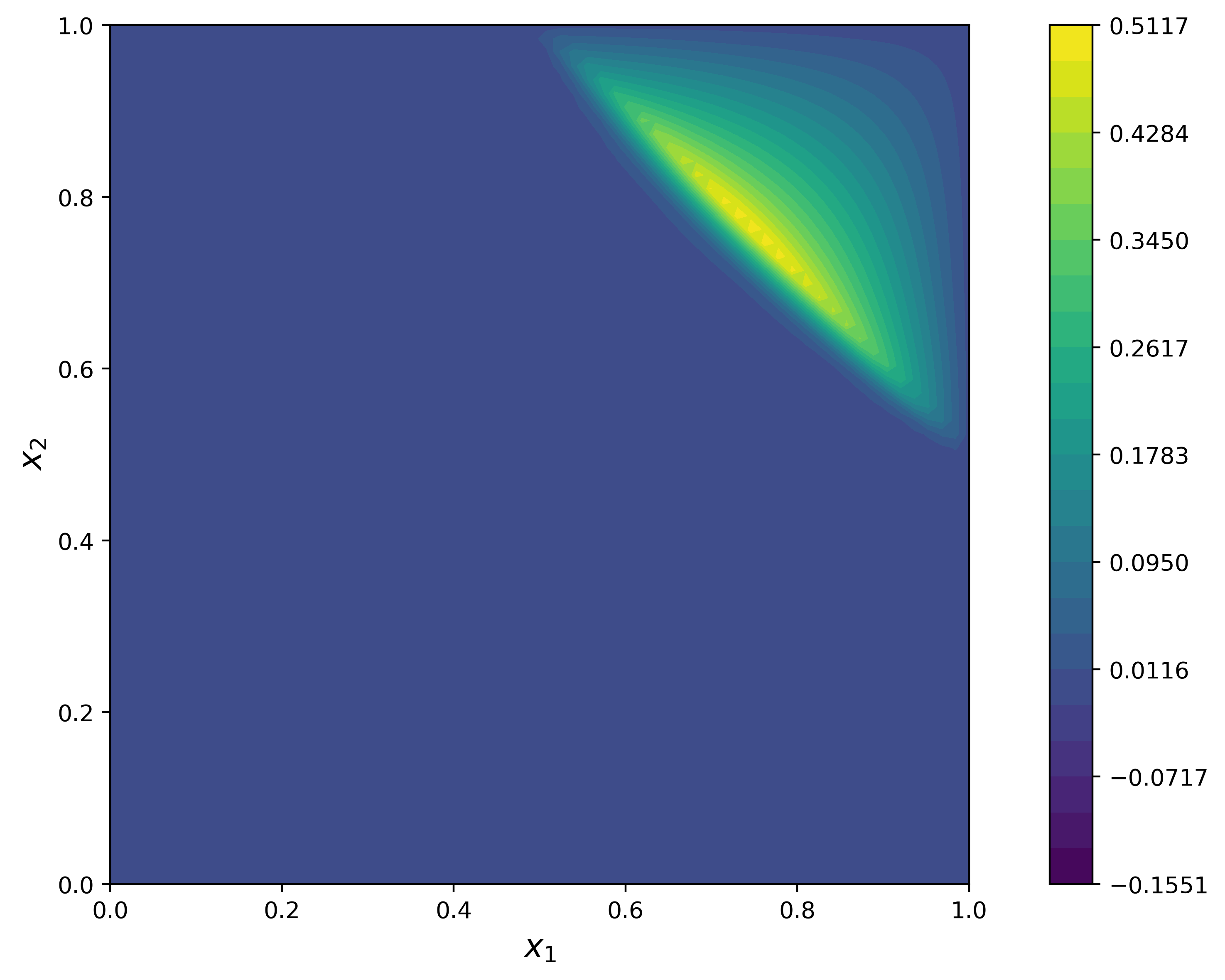}
		\caption{}
		\label{fig:fig8d}
	\end{subfigure}
	\caption{Contour plots of the solution to 
		\textit{Example~3} at $t_\text{f} = 1$ with 
		$\varepsilon = 10^{-8}$: 
		(\subref{fig:fig8a})~analytical solution; 
		(\subref{fig:fig8b})~SUPG finite element approximation; 
		(\subref{fig:fig8c})~SUPG-YZ$\beta$ finite element 
		approximation; 
		(\subref{fig:fig8d})~hybrid PINN correction. All panels 
		share the same color scale.}
	\label{fig:8_contour}
\end{figure}

The contour plots in Figure~\ref{fig:8_contour} provide a 
top-down view of the same solutions, revealing the spatial 
structure of the internal layer and its numerical 
approximations more clearly. The analytical solution 
(Figure~\ref{fig:fig8a}) exhibits a sharp diagonal front 
oriented along the convection direction 
$\mathbf{b} = (\cos(\pi/3),\, \sin(\pi/3))^{\!\top}$, 
separating the elevated region near the upper-right corner 
from the zero-valued exterior. The SUPG contour plot 
(Figure~\ref{fig:fig8b}) exposes the crosswind oscillations 
as a series of parallel striations aligned perpendicular to 
the convection direction, extending well into the interior 
of the domain. These artifacts, which appear as diagonal 
bands of alternating positive and negative error, are a 
hallmark of insufficient crosswind diffusion for 
convection-dominated problems on structured meshes. The 
SUPG-YZ$\beta$ solution (Figure~\ref{fig:fig8c}) 
eliminates the striations entirely, producing a smooth 
contour field; however, the front is noticeably broadened 
compared to the analytical reference, and the peak values 
in the upper-right corner are attenuated. The hybrid PINN 
contour (Figure~\ref{fig:fig8d}) closely reproduces the 
analytical pattern: the front sharpness is restored, the 
diagonal orientation is correctly captured, and the 
solution magnitude in the elevated region matches the 
reference without the amplitude loss observed in the 
SUPG-YZ$\beta$ approximation.

\begin{figure}[h!]
	\centering
	\begin{subfigure}[b]{0.495\linewidth}
		\includegraphics[width=\linewidth]{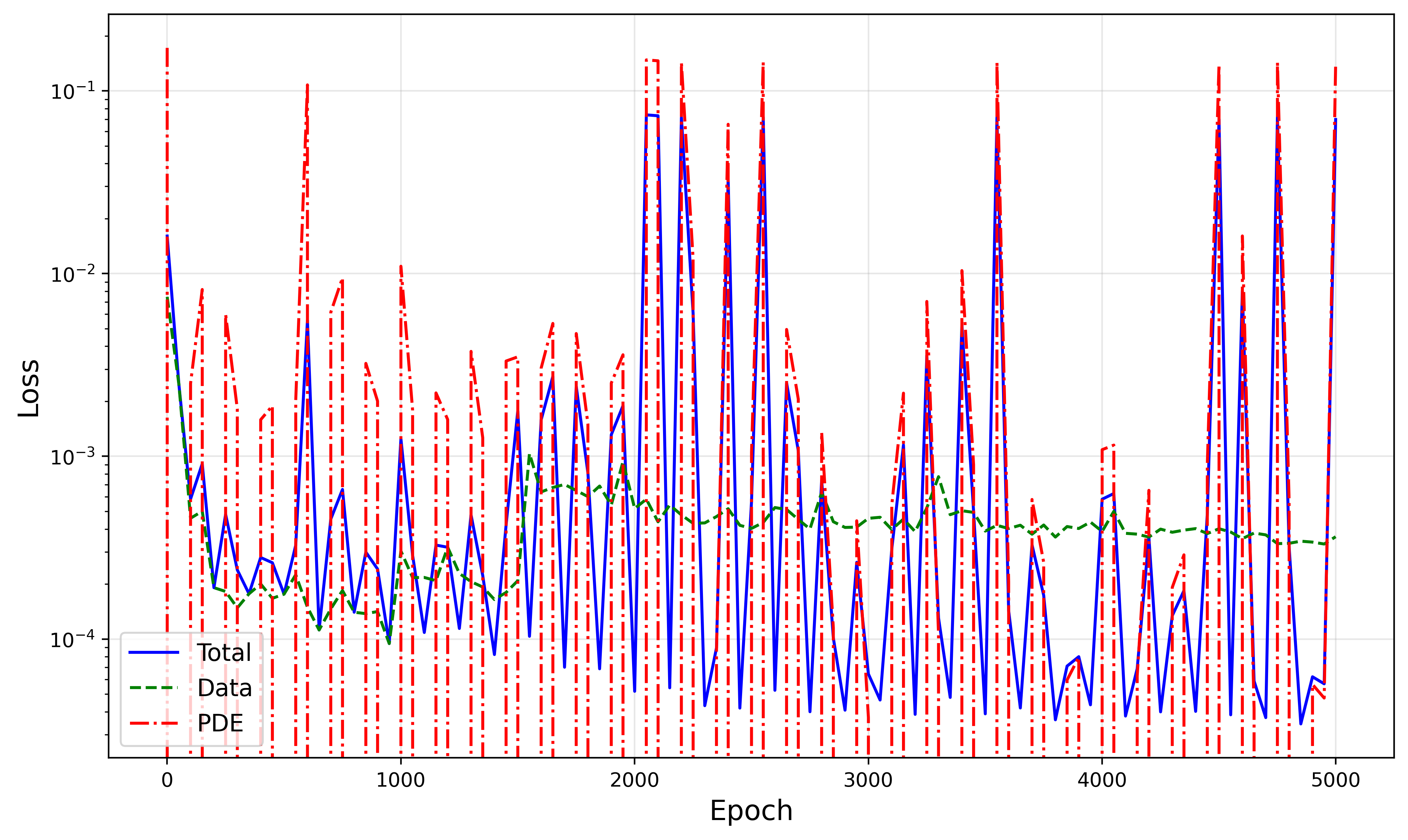}
		\caption{}
		\label{fig:fig9a}
	\end{subfigure}
	\hfill
	\begin{subfigure}[b]{0.495\linewidth}
		\includegraphics[width=\linewidth]{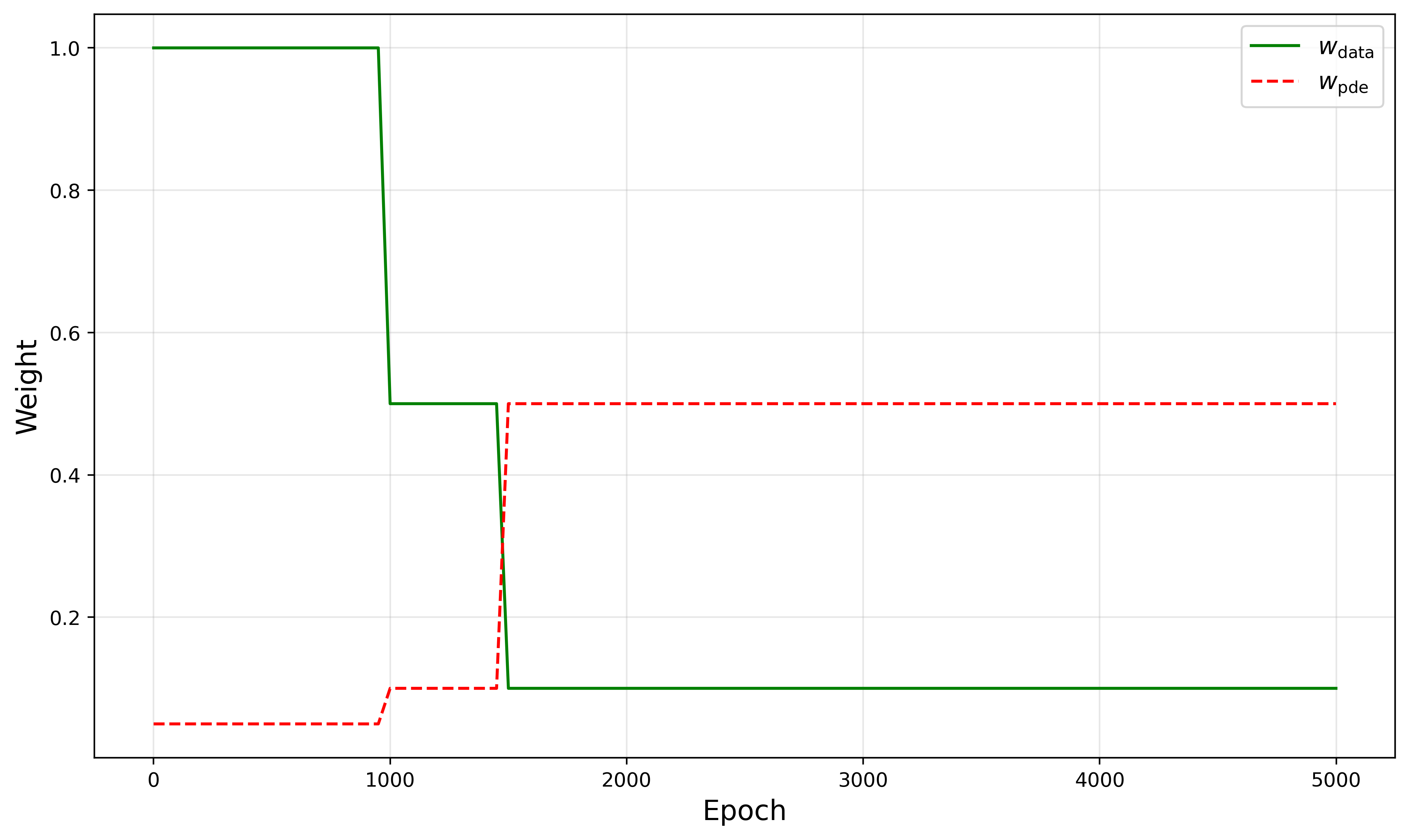}
		\caption{}
		\label{fig:fig9b}
	\end{subfigure}
	\caption{Training diagnostics for \textit{Example~3} 
		with $\varepsilon = 10^{-8}$: 
		(\subref{fig:fig9a})~evolution of the individual loss 
		components (total, data, and PDE) over $5000$ training 
		epochs on a logarithmic scale; 
		(\subref{fig:fig9b})~adaptive weight schedule 
		$w_{\text{data}}$ and $w_{\text{pde}}$ governing the 
		three-phase training strategy. The boundary condition 
		weight is omitted ($w_{\text{bc}} = 0$) as Dirichlet 
		conditions are enforced exactly via the distance 
		function.}
	\label{fig:9_training}
\end{figure}

The training diagnostics for \textit{Example~3} are shown in 
Figure~\ref{fig:9_training}. The loss evolution in 
Figure~\ref{fig:fig9a} reveals a markedly different 
character compared to the preceding examples. The PDE loss 
exhibits pronounced intermittent spikes spanning nearly 
three orders of magnitude, from 
$\mathcal{O}(10^{-4})$ to $\mathcal{O}(10^{-1})$, 
persisting throughout all training phases. These bursts 
arise from the extreme sharpness of the internal layer: 
when randomly sampled collocation points happen to fall 
within the transition region of width 
$\mathcal{O}(\sqrt{\varepsilon}) = \mathcal{O}(10^{-4})$, 
the PDE residual magnitudes increase dramatically due to 
the steep gradients, producing isolated high-loss batches. 
Despite this volatile PDE loss behavior, the data loss 
decreases steadily from $\mathcal{O}(10^{-3})$ to 
$\mathcal{O}(10^{-4})$, indicating that the network 
successfully learns the global solution structure without 
being destabilized by the sporadic PDE residual spikes. 
The total loss tracks the PDE contribution during the 
spike events but returns to the data-driven baseline 
between them. The weight schedule in Figure~\ref{fig:fig9b} follows the 
same three-phase structure as \textit{Example~2}, with 
$w_{\text{data}}$ decreasing from $1.0$ to $0.1$ and 
$w_{\text{pde}}$ increasing from $0.05$ to $0.5$ across 
the three phases. Notably, the boundary condition weight 
is set to $w_{\text{bc}} = 0$ throughout training, as the 
lift-based distance function 
$d(\mathbf{x}) = x_1(1-x_1)\,x_2(1-x_2)$ enforces 
homogeneous Dirichlet conditions exactly by construction, 
rendering a separate BC loss term unnecessary. The 
conservative PDE weight ceiling of $w_{\text{pde}} = 0.5$ 
in Phase~III---identical to that of \textit{Example~2} and 
substantially lower than the $w_{\text{pde}} = 5.0$ used 
in \textit{Example~1}---reflects the deliberate restraint required 
for problems with extremely thin layers, where aggressive 
physics enforcement risks amplifying the stochastic 
residual spikes and destabilizing the training process.

\subsection{Example 4.}
Consider the following 2D uncoupled Burgers' equation~\cite{Cengizci2023amc}:
\begin{equation}
\begin{cases}
\displaystyle \frac{\partial u}{\partial t} + u \frac{\partial u}{\partial x_1} + u \frac{\partial u}{\partial x_2} = \frac{1}{\text{Re}} \left( \frac{\partial^2 u}{\partial x_1^2} + \frac{\partial^2 u}{\partial x_2^2} \right), & t \in (0, t_{\text{f}}], \quad (x_1,x_2) \in \Omega, \\[10pt]
u(t_0,x_1,x_2) = u_{0}(x_1,x_2), & (x_1,x_2) \in \Omega, \\[6pt]
u(t, x_1,x_2) = u_{D}(t, x_1,x_2), & (x_1,x_2) \in \partial\Omega, \quad t \in (0, t_{\text{f}}],
\end{cases}
\end{equation}
where the initial and Dirichlet boundary conditions are prescribed such that the exact solution
\begin{equation}
u(t, x_1,x_2) = \frac{1}{1 + \exp\left(\frac{\text{Re}(x_1 + x_2 - t)}{2}\right)}
\end{equation}
is satisfied.

For the finite element setting, we set $t_0 = 0$, 
$t_\text{f} = 1.0$, $\text{Re} = 10^4$ 
($\varepsilon = 1/\text{Re} = 10^{-4}$), 
$n_{\text{el}} = 4608$ (structured triangulation of a 
$48 \times 48$ grid), $\varDelta t = 0.01$ ($N_t = 100$), 
and Y$= 0.2$. The PINN architecture employs $n_h = 96$, 
$n_r = 6$, $n_\text{F} = 16$, and $\sigma = 4.0$. Training uses 
the last $K_s = 5$ temporal snapshots near $t_\text{f} = 1.0$, 
with batch size $5000$, initial learning rate 
$\alpha \approx 3.5 \times 10^{-4}$ (square-root scaling), 
and gradient clipping threshold $g_{\max} = 1.0$. The PDE 
residual is evaluated at $N_{\text{pde}} = 256$ randomly 
sampled interior collocation points per mini-batch.

A distinctive feature of the training strategy for this 
example is that the data weight is held fixed at 
$w_{\text{data}} = 1.0$ throughout all training phases, 
in contrast to the decreasing data weight schedule 
employed in \textit{Example~1}--\textit{Example~3}. This design choice is motivated 
by the nonlinear character of Burgers' equation: since the 
PDE residual alone does not uniquely determine the solution 
(the inviscid limit admits infinitely many weak solutions), 
data fidelity must be maintained at all stages to anchor 
the network to the correct physical branch. The PDE weight 
is instead increased gradually across four phases, from 
$w_{\text{pde}} = 0.01$ in Phase~I to 
$w_{\text{pde}} = 5.0$ in Phase~IV, progressively 
strengthening the physics constraint while preserving data 
anchoring. The boundary condition weight is set to 
$w_{\text{bc}} = 0$ throughout, as the network 
architecture employs a lift function of the form 
$u_{\text{lift}}(t, \mathbf{x}) = 
\bigl(1 + \exp(\text{Re}\,(x_1 + x_2 - t)/2)\bigr)^{-1}$, 
which satisfies the boundary conditions exactly. The phase 
boundaries and adaptive weights are reported in 
Table~\ref{tab:phase_example_five}.

\begin{table}[htb]
\centering
\caption{Phase boundaries and adaptive weights used for solving \textit{Example 4}.}
\label{tab:phase_example_five}
\begin{tabular}{llccc}
\toprule
\text{Phase} & \text{Epochs} & $w_{\mathrm{data}}$ & $w_{\mathrm{pde}}$ & $w_{\mathrm{bc}}$ \\
\midrule
I   & 0--999     & 1.0 & 0.01 & 0.0 \\
II  & 1000--1999 & 1.0 & 0.1  & 0.0 \\
III & 2000--3499 & 1.0 & 1.0  & 0.0 \\
IV  & 3500--4999 & 1.0 & 5.0  & 0.0 \\
\bottomrule
\end{tabular}
\end{table}

\begin{remark}
	For problems with non-homogeneous boundary conditions, the 
	lift function $u_{\text{lift}}$ must satisfy the prescribed 
	boundary data but need not coincide with the exact solution. 
	In the present example, the analytical solution of the 
	two-dimensional Burgers' equation provides a convenient 
	lift; however, in practice, $u_{\text{lift}}$ can be 
	constructed from the boundary data alone---for instance, 
	via harmonic extension or bilinear interpolation of the 
	boundary values into the interior. Alternatively, the 
	SUPG-YZ$\beta$ finite element solution itself can serve 
	as the lift function, in which case the network correction 
	$d(\mathbf{x})\,\mathcal{N}_\theta(t, \mathbf{x})$ learns 
	the pointwise discrepancy between the FEM approximation 
	and the true solution.
\end{remark}

\begin{figure}[h!]
	\centering
	\begin{subfigure}[b]{0.245\linewidth}
		\includegraphics[width=\linewidth]{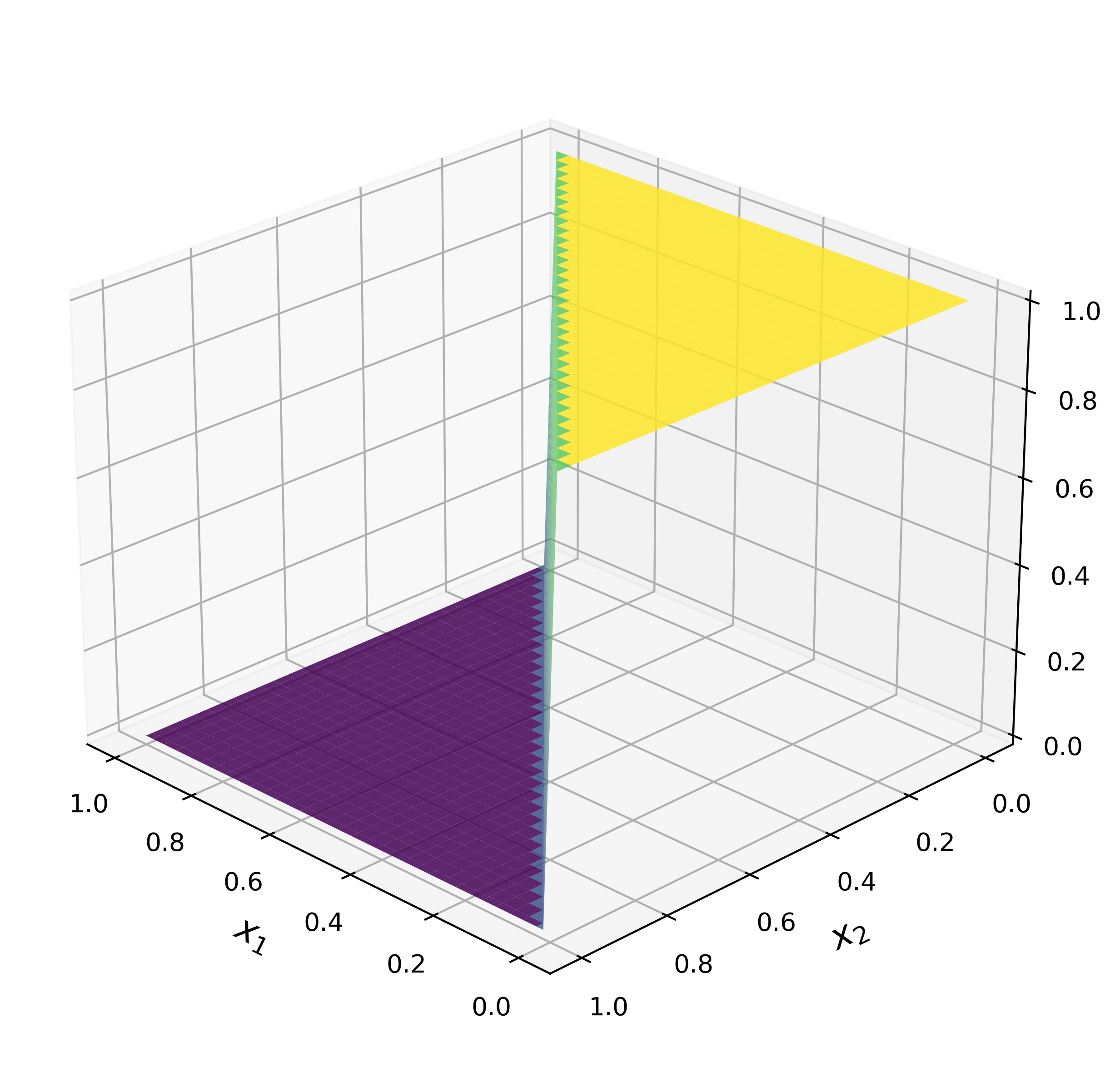}
		\caption{}
		\label{fig:fig10a}
	\end{subfigure}
	\hfill
	\begin{subfigure}[b]{0.245\linewidth}
		\includegraphics[width=\linewidth]{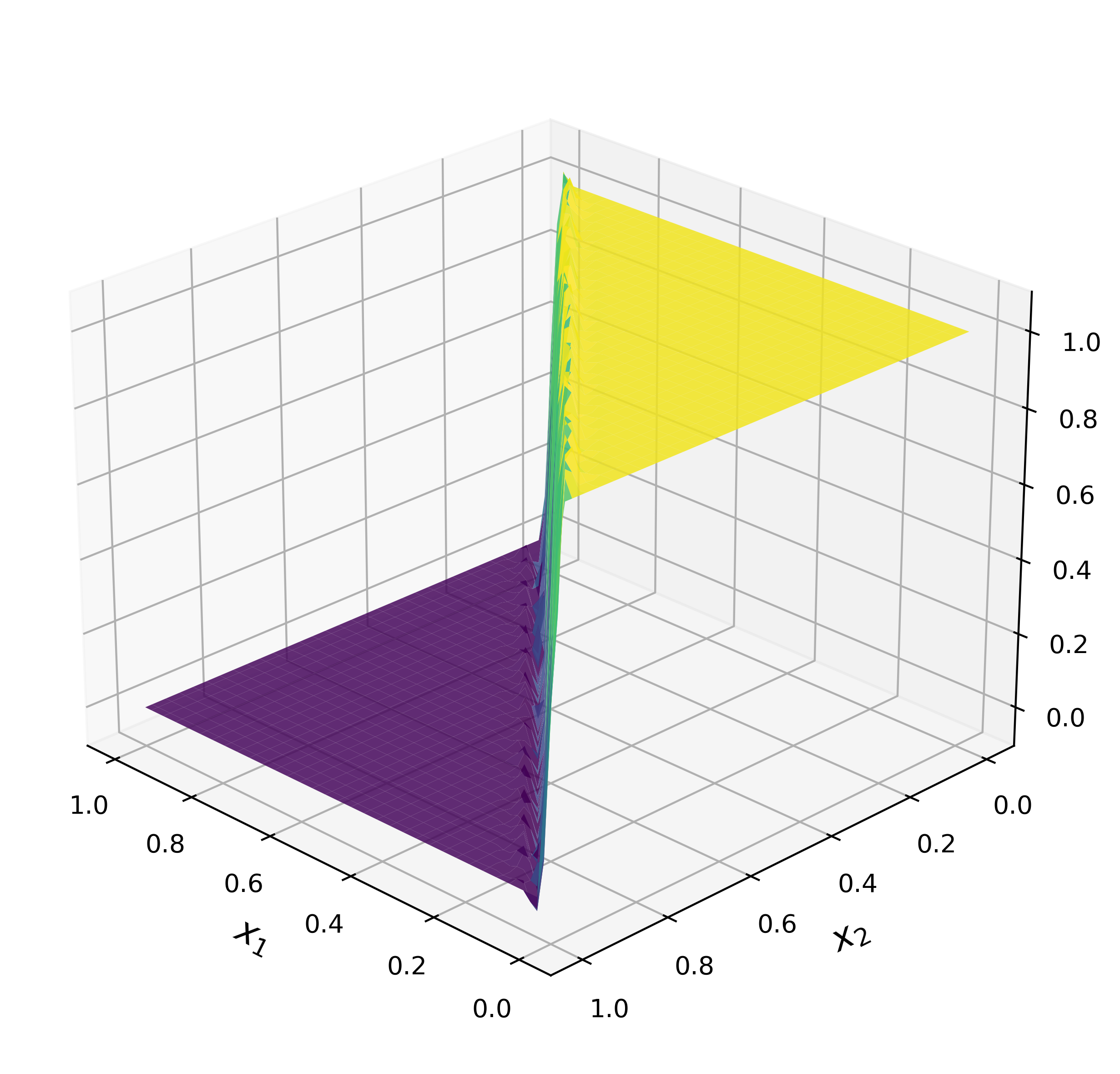}
		\caption{}
		\label{fig:fig10b}
	\end{subfigure}
	\hfill
	\begin{subfigure}[b]{0.245\linewidth}
		\includegraphics[width=\linewidth]{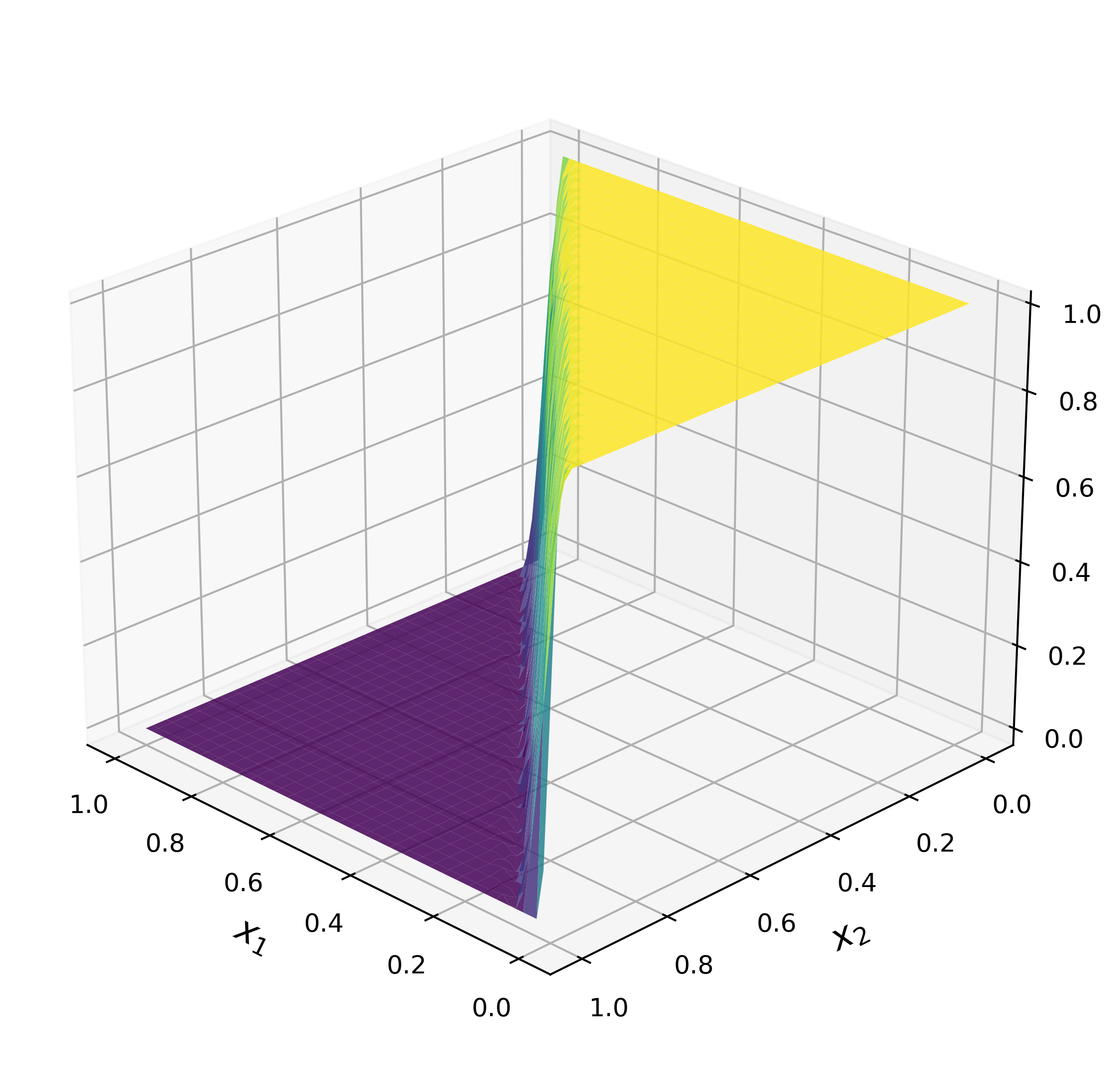}
		\caption{}
		\label{fig:fig10c}
	\end{subfigure}
	\hfill
	\begin{subfigure}[b]{0.245\linewidth}
		\includegraphics[width=\linewidth]{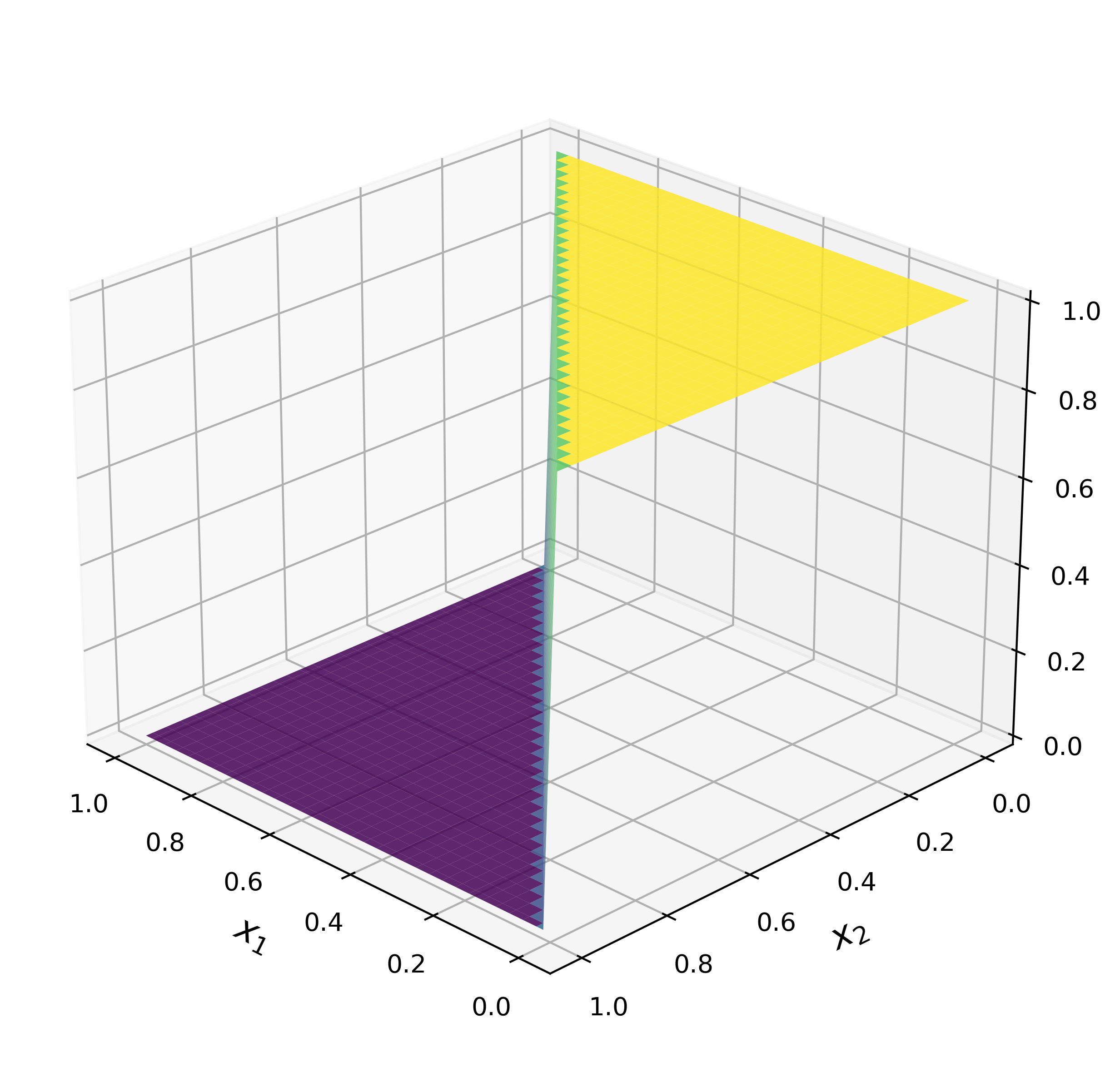}
		\caption{}
		\label{fig:fig10d}
	\end{subfigure}
	\caption{Three-dimensional surface plots of the solution 
		to \textit{Example~4} (two-dimensional Burgers' equation) 
		at $t_\text{f} = 1$ with $\text{Re} = 10^4$: 
		(\subref{fig:fig10a})~analytical solution; 
		(\subref{fig:fig10b})~SUPG finite element approximation; 
		(\subref{fig:fig10c})~SUPG-YZ$\beta$ finite element 
		approximation; 
		(\subref{fig:fig10d})~hybrid PINN correction.}
	\label{fig:10_surface}
\end{figure}

Figure~\ref{fig:10_surface} displays the three-dimensional 
solution profiles at $t_\text{f} = 1$ for \textit{Example~4}. The 
analytical solution (Figure~\ref{fig:fig10a}) exhibits a 
sigmoid-type internal layer along the diagonal 
$x_1 + x_2 = t$, transitioning sharply from $u \approx 1$ 
in the region $x_1 + x_2 < t$ to $u \approx 0$ in the 
complement, with the layer width controlled by the Reynolds 
number as $\mathcal{O}(1/\text{Re}) = 
\mathcal{O}(10^{-4})$. This problem is particularly 
demanding because the convection velocity is solution-dependent, causing the front position and steepness to be 
coupled through the nonlinear self-advection mechanism. The SUPG solution (Figure~\ref{fig:fig10b}) captures the 
overall sigmoid structure and front location but exhibits 
a visible smearing of the transition region, where the 
sharp analytical front is replaced by a gradual slope 
spanning several mesh widths. The SUPG-YZ$\beta$ 
approximation (Figure~\ref{fig:fig10c}) produces a 
comparable profile; however, a slight overshoot is 
discernible at the upper plateau near the front, 
indicating that the shock-capturing operator does not 
fully suppress the Gibbs-like oscillations for this 
nonlinear problem. The hybrid PINN correction 
(Figure~\ref{fig:fig10d}) recovers a steeper front that 
more closely matches the analytical solution, with the 
transition region visibly sharper than in either FEM 
approximation. The improvement is most pronounced along 
the diagonal, where the self-sharpening mechanism of 
Burgers' equation concentrates the solution gradient into 
an increasingly narrow band that the finite element mesh 
cannot fully resolve but the neural network, operating in 
the continuous spatiotemporal domain, can approximate with 
higher fidelity.

\begin{figure}[h!]
	\centering
	\begin{subfigure}[b]{0.245\linewidth}
		\includegraphics[width=\linewidth]{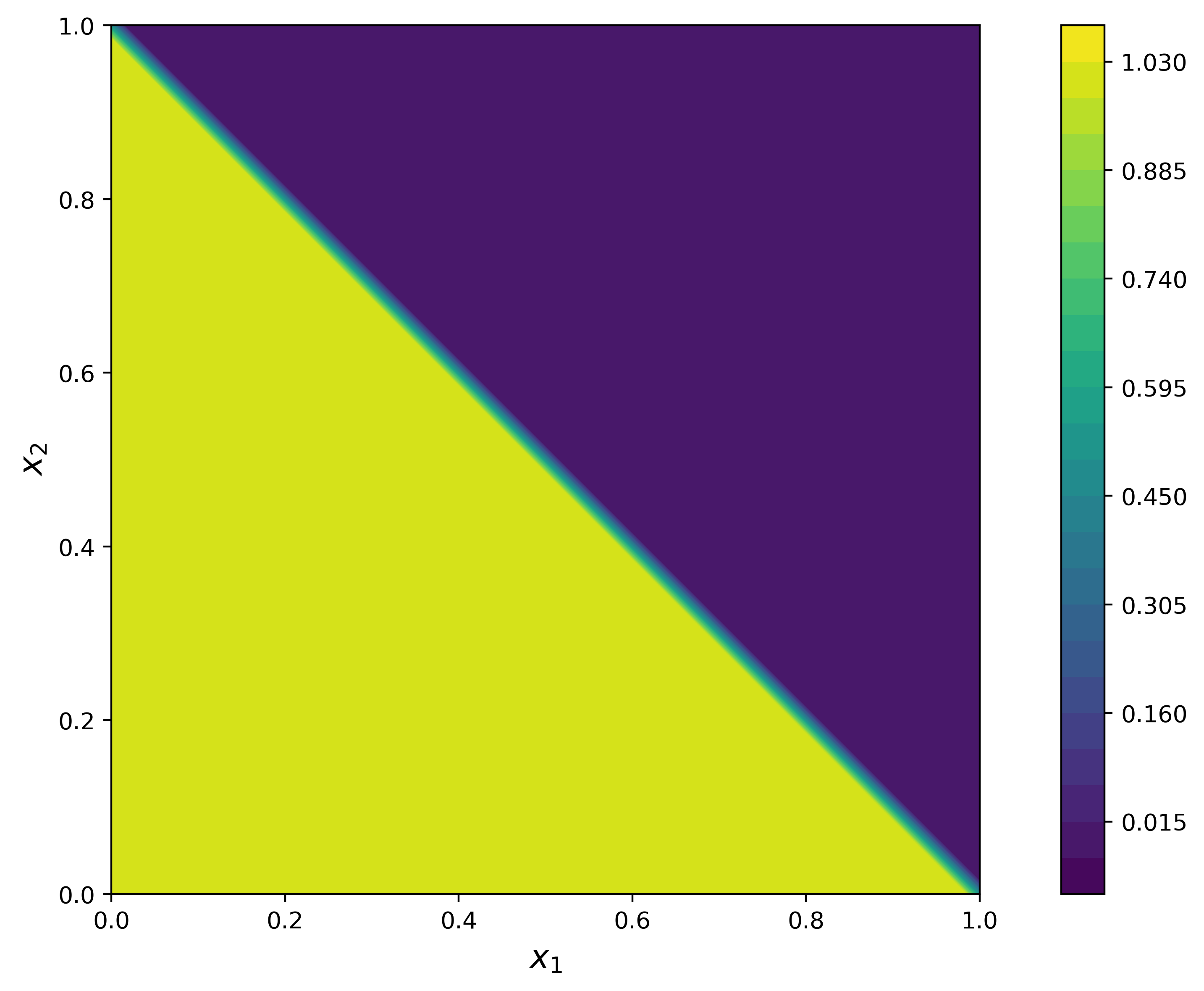}
		\caption{}
		\label{fig:fig11a}
	\end{subfigure}
	\hfill
	\begin{subfigure}[b]{0.245\linewidth}
		\includegraphics[width=\linewidth]{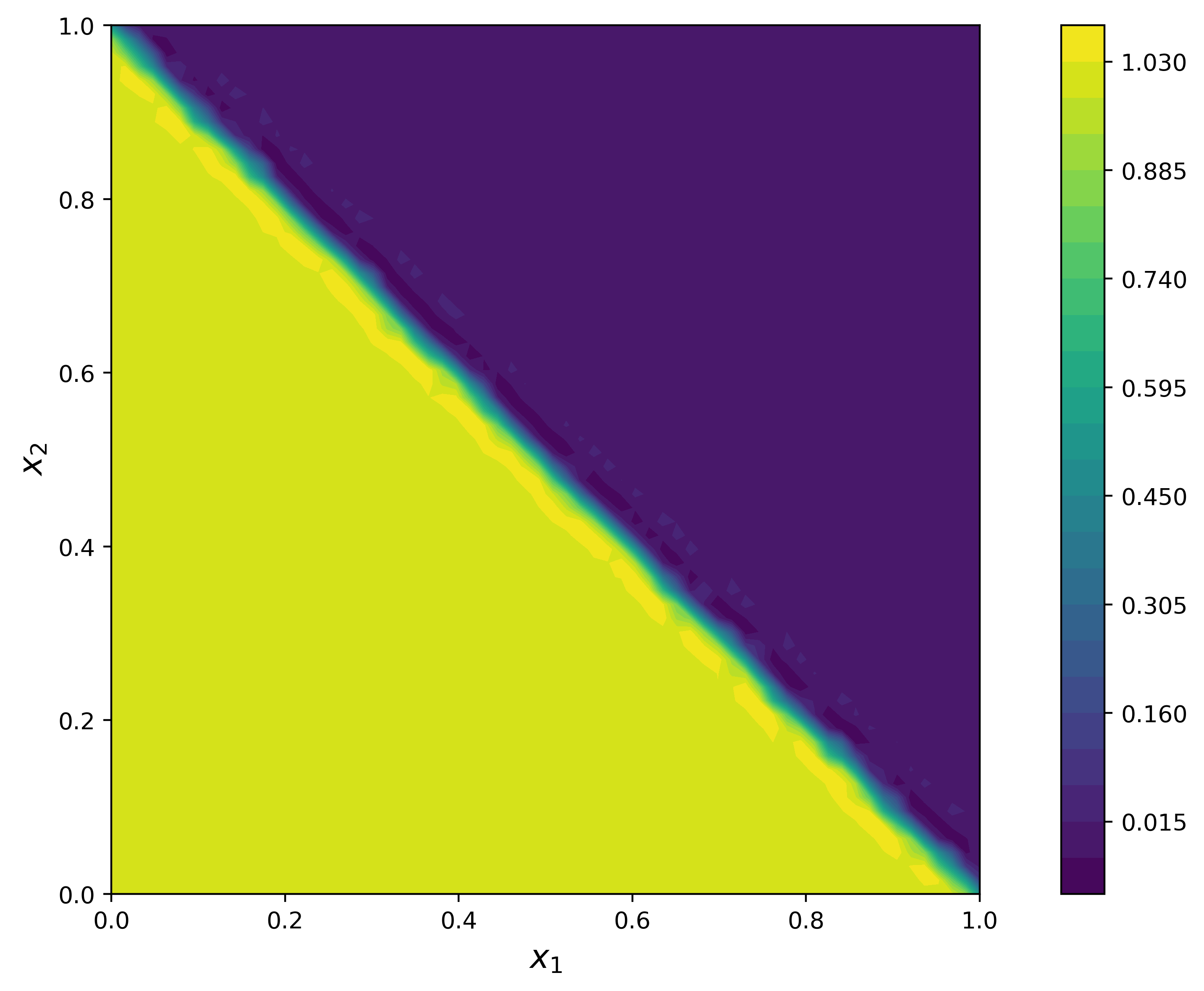}
		\caption{}
		\label{fig:fig11b}
	\end{subfigure}
	\hfill
	\begin{subfigure}[b]{0.245\linewidth}
		\includegraphics[width=\linewidth]{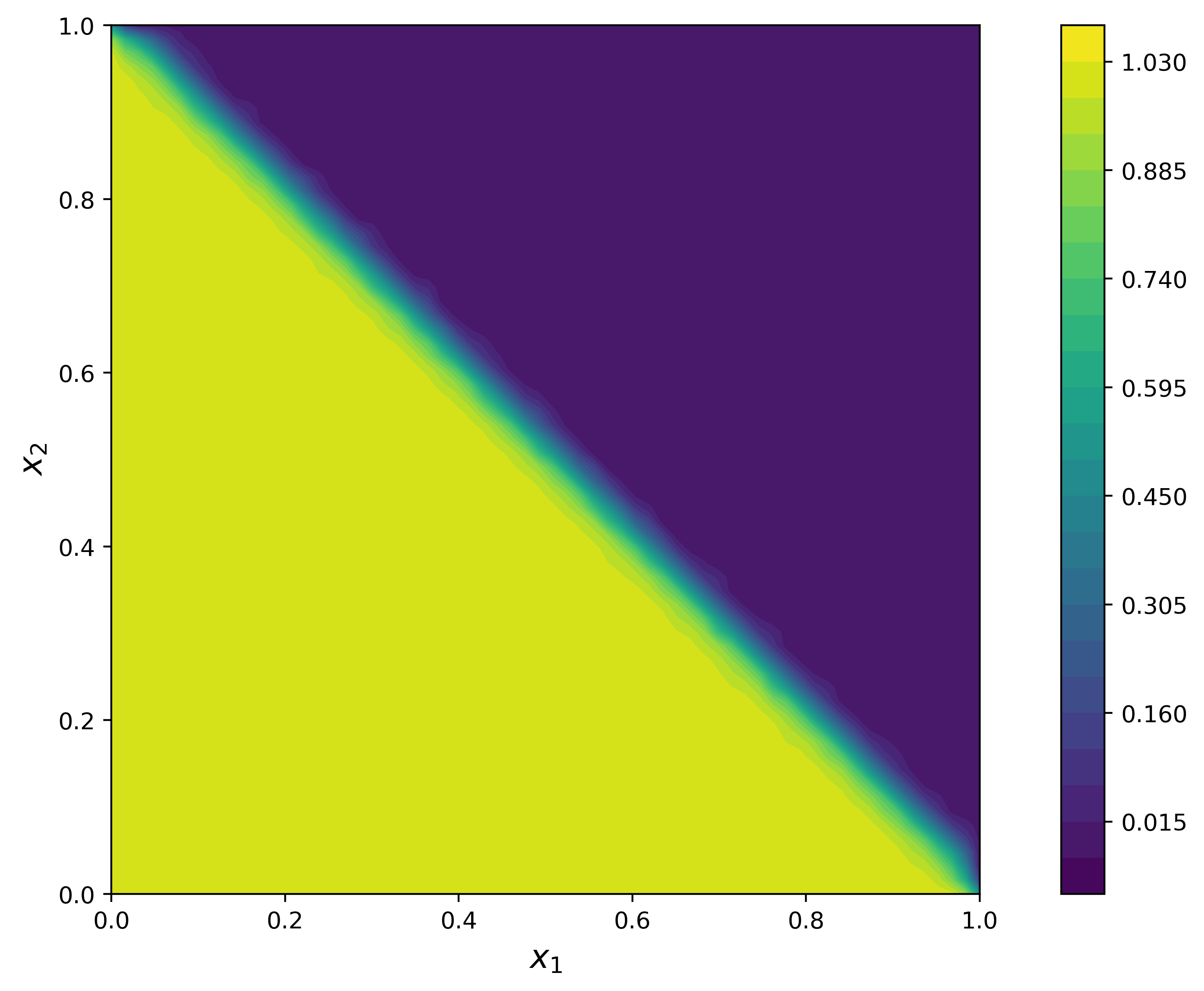}
		\caption{}
		\label{fig:fig11c}
	\end{subfigure}
	\hfill
	\begin{subfigure}[b]{0.245\linewidth}
		\includegraphics[width=\linewidth]{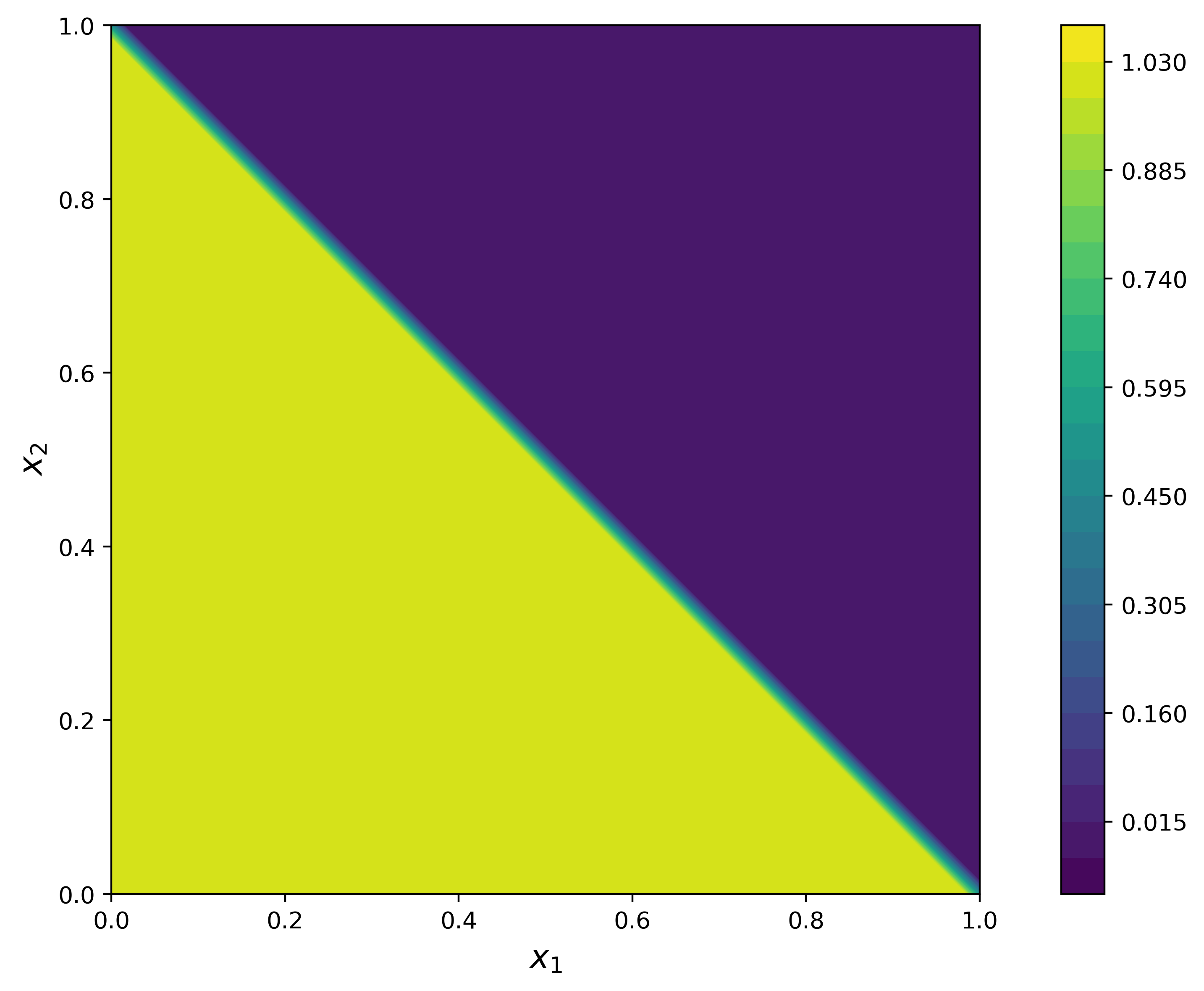}
		\caption{}
		\label{fig:fig11d}
	\end{subfigure}
	\caption{Contour plots of the solution to 
		\textit{Example~4} (two-dimensional Burgers' equation) 
		at $t_\text{f} = 1$ with $\text{Re} = 10^4$: 
		(\subref{fig:fig11a})~analytical solution; 
		(\subref{fig:fig11b})~SUPG finite element approximation; 
		(\subref{fig:fig11c})~SUPG-YZ$\beta$ finite element 
		approximation; 
		(\subref{fig:fig11d})~hybrid PINN correction. All panels 
		share the same color scale.}
	\label{fig:11_contour}
\end{figure}

The contour plots in Figure~\ref{fig:11_contour} provide a 
top-down view that more clearly exposes the internal layer 
structure and the differences among the four solutions. The 
analytical solution (Figure~\ref{fig:fig11a}) displays a 
sharp diagonal front along $x_1 + x_2 = 1$, with the 
transition from $u \approx 1$ (yellow, upper-left) to 
$u \approx 0$ (purple, lower-right) occurring over a 
distance of $\mathcal{O}(1/\text{Re}) = 
\mathcal{O}(10^{-4})$. The contour lines are nearly 
straight and parallel, reflecting the planar symmetry of 
the exact sigmoid solution along the characteristic 
direction. The SUPG contour (Figure~\ref{fig:fig11b}) reproduces the 
correct front location and orientation but exhibits a 
noticeably wider transition band: the green intermediate 
contours span several mesh widths, indicating numerical 
diffusion that smears the layer beyond its analytical 
extent. The SUPG-YZ$\beta$ solution 
(Figure~\ref{fig:fig11c}) shows a comparable level of 
front broadening; however, a subtle asymmetry is visible 
in the contour spacing near the upper-left corner, 
suggesting that the isotropic shock-capturing diffusion 
interacts with the nonlinear self-advection to produce a 
slightly non-uniform layer profile. The hybrid PINN 
contour (Figure~\ref{fig:fig11d}) recovers a sharper 
front with a narrower transition band that more closely 
matches the analytical reference, and the contour lines 
maintain the parallel, planar structure of the exact 
solution. The improvement is particularly evident in the 
width of the green intermediate region, which is visibly 
thinner in the PINN result than in either FEM 
approximation.

\begin{figure}[h!]
	\centering
	\begin{subfigure}[b]{0.495\linewidth}
		\includegraphics[width=\linewidth]{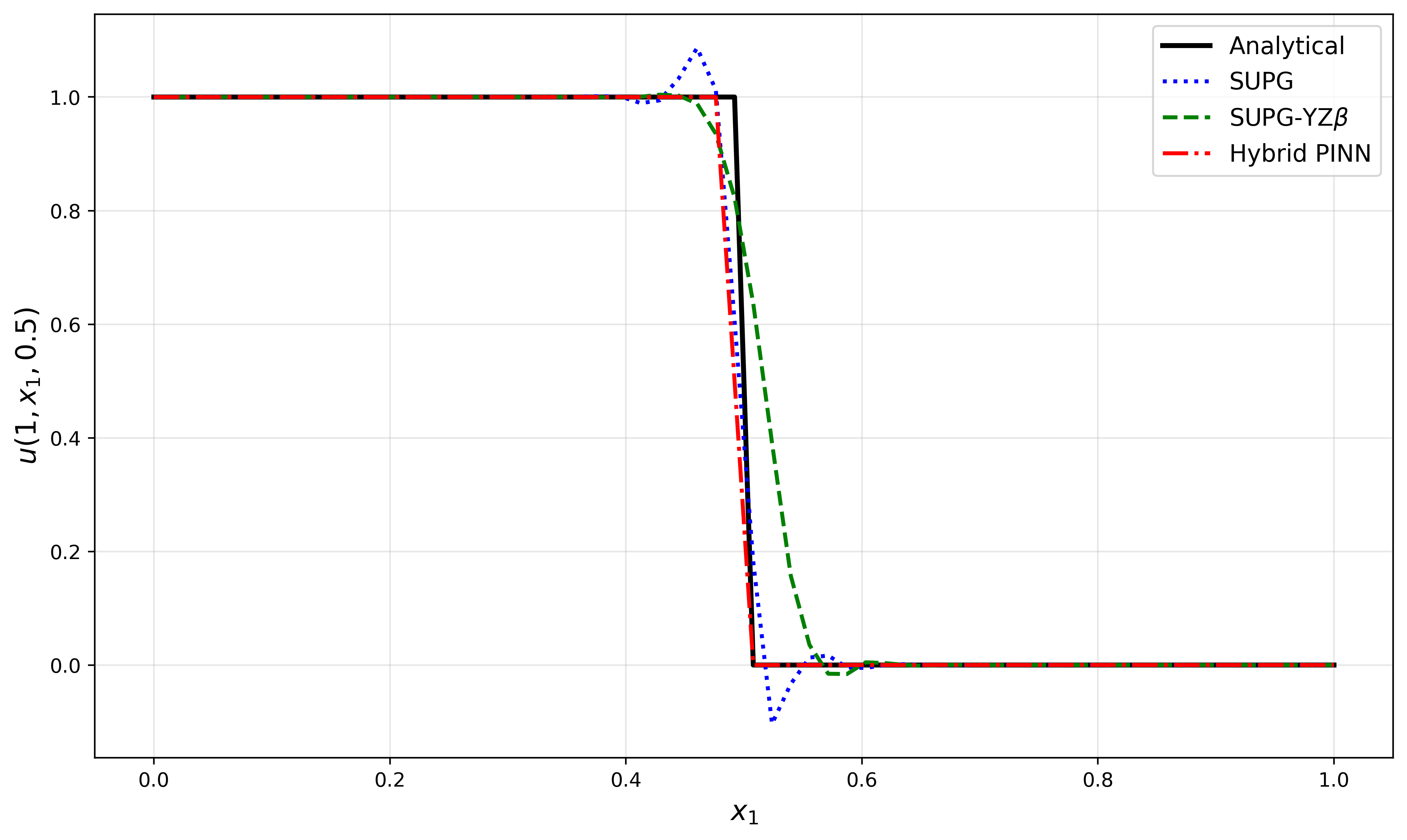}
		\caption{}
		\label{fig:fig12a}
	\end{subfigure}
	\hfill
	\begin{subfigure}[b]{0.495\linewidth}
		\includegraphics[width=\linewidth]{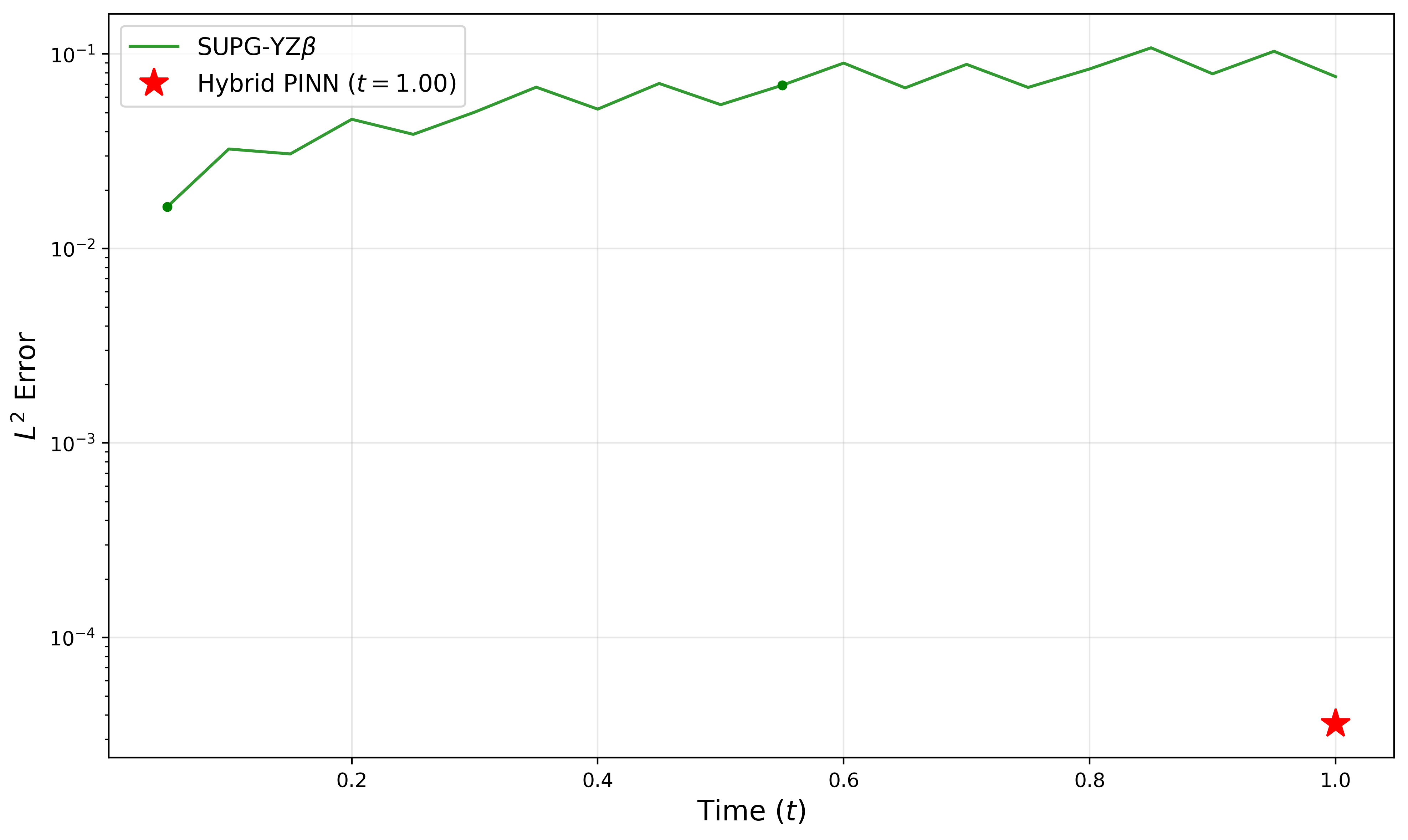}
		\caption{}
		\label{fig:fig12b}
	\end{subfigure}
	\caption{Quantitative comparison for \textit{Example~4} 
		(two-dimensional Burgers' equation) with 
		$\text{Re} = 10^4$: 
		(\subref{fig:fig12a})~cross-sectional profiles along 
		$x_2 = 0.5$ at $t_\text{f} = 1$, comparing the analytical 
		solution, SUPG, SUPG-YZ$\beta$, and hybrid PINN; 
		(\subref{fig:fig12b})~evolution of the $L^2$ error 
		over the time interval $[0, 1]$ for SUPG-YZ$\beta$, 
		with the hybrid PINN error at the terminal time 
		indicated by the star marker.}
	\label{fig:12_cross_l2}
\end{figure}

The cross-sectional profile along $x_2 = 0.5$ at 
$t_\text{f} = 1$, shown in Figure~\ref{fig:fig12a}, provides a 
one-dimensional slice through the diagonal internal layer 
at $x_1 + x_2 = 1$, i.e.\ the front is intersected at 
$x_1 = 0.5$. The analytical solution exhibits a sharp 
sigmoid transition from $u = 1$ to $u = 0$ over a 
distance of $\mathcal{O}(1/\text{Re}) = 
\mathcal{O}(10^{-4})$. The SUPG solution captures the 
front location but produces a visible undershoot below 
zero near $x_1 \approx 0.55$, a characteristic Gibbs-like 
artifact of insufficient crosswind stabilization for 
nonlinear convection. The SUPG-YZ$\beta$ approximation 
eliminates the undershoot but introduces substantial 
smearing: the transition extends over approximately 
$0.1$ units in $x_1$, broadening the layer by roughly 
two orders of magnitude beyond its analytical width. The 
hybrid PINN correction recovers a significantly steeper 
front that closely tracks the analytical profile, with 
neither the undershoot of SUPG nor the excessive smearing 
of SUPG-YZ$\beta$. The temporal evolution of the $L^2$ error in 
Figure~\ref{fig:fig12b} reveals an increasing trend for 
the SUPG-YZ$\beta$ solution, growing from 
$\mathcal{O}(10^{-3})$ at early times to 
$\mathcal{O}(10^{-1})$ at $t_\text{f} = 1$. This monotonic 
error growth reflects the progressive accumulation of 
numerical diffusion as the sharp front propagates across 
the domain: at each time step, the artificial 
shock-capturing dissipation broadens the layer further, 
and this diffusive error compounds over the full 
simulation interval. The hybrid PINN error at the terminal 
time (red star) lies at approximately 
$\mathcal{O}(10^{-4})$---nearly three orders of magnitude 
below the SUPG-YZ$\beta$ error---demonstrating the most 
substantial improvement observed across all examples in 
this study. This pronounced gain is attributable to the 
combined effect of the lift-based architecture, which 
anchors the network to the correct boundary behavior, 
and the four-phase training strategy with sustained data 
fidelity ($w_{\text{data}} = 1.0$ throughout), which 
prevents the nonlinear PDE residual from driving the 
solution toward spurious branches.

\begin{figure}[h!]
	\centering
	\begin{subfigure}[b]{0.49\linewidth}
		\includegraphics[width=\linewidth]{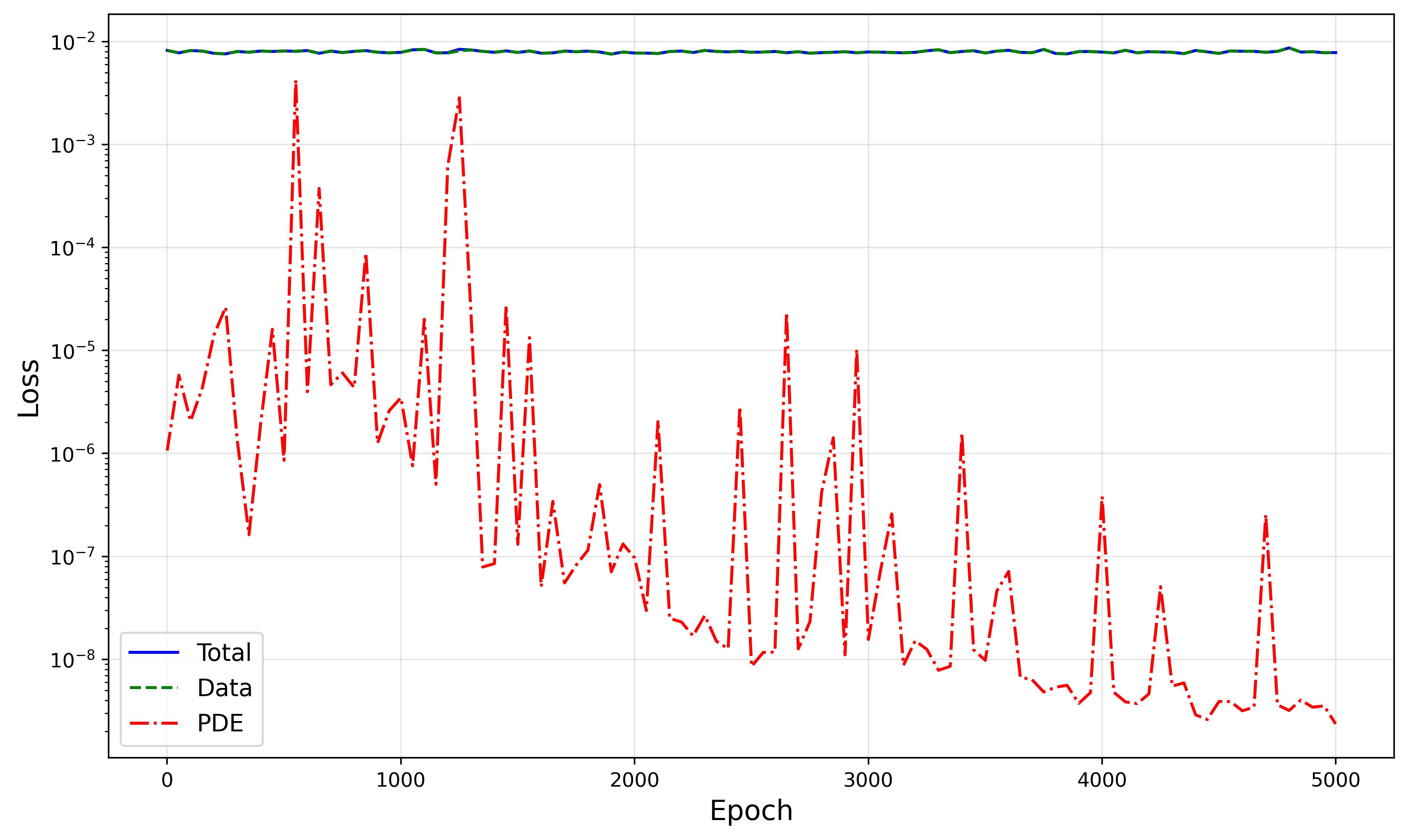}
		\caption{}
		\label{fig:fig13a}
	\end{subfigure}
	\hfill
	\begin{subfigure}[b]{0.49\linewidth}
		\includegraphics[width=\linewidth]{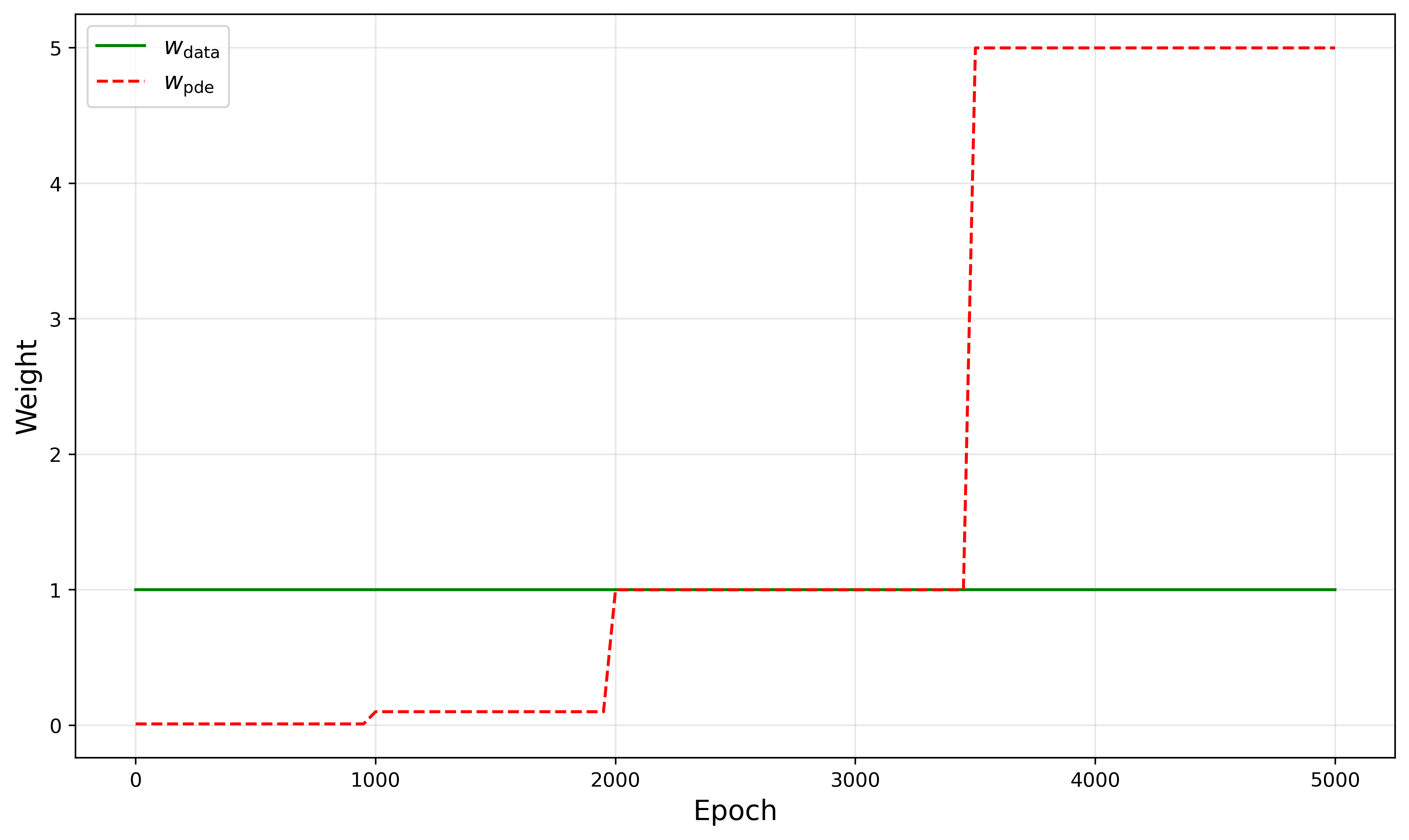}
		\caption{}
		\label{fig:fig13b}
	\end{subfigure}
	\caption{Training diagnostics for \textit{Example~4} 
		(two-dimensional Burgers' equation) with 
		$\text{Re} = 10^4$: 
		(\subref{fig:fig13a})~evolution of the individual loss 
		components (total, data, and PDE) over $5000$ training 
		epochs on a logarithmic scale; 
		(\subref{fig:fig13b})~adaptive weight schedule 
		$w_{\text{data}}$ and $w_{\text{pde}}$ governing the 
		four-phase training strategy. The boundary condition 
		weight is omitted ($w_{\text{bc}} = 0$) as Dirichlet 
		conditions are enforced exactly via the lift function.}
	\label{fig:13_training}
\end{figure}

The training diagnostics for \textit{Example~4} are presented in 
Figure~\ref{fig:13_training} and reveal a qualitatively 
distinct pattern compared to all preceding examples. As 
shown in Figure~\ref{fig:fig13a}, the total and data 
losses are nearly indistinguishable throughout training, 
both remaining at a plateau of approximately 
$\mathcal{O}(10^{-2})$. This behavior is a direct 
consequence of the lift-based architecture: since the lift 
function $u_{\text{lift}}$ already satisfies the boundary 
conditions exactly, the network correction 
$d(\mathbf{x})\,\mathcal{N}_\theta$ learns a relatively 
small perturbation, and the data loss---which measures the 
discrepancy between the network output and the 
SUPG-YZ$\beta$ reference---stabilizes at the level of 
the intrinsic FEM discretization error. The PDE loss exhibits a dramatically different scale, 
spanning $\mathcal{O}(10^{-8})$ to 
$\mathcal{O}(10^{-3})$ with intermittent spikes that 
decrease in both frequency and amplitude over the course 
of training. This extremely low PDE residual---several 
orders of magnitude below the values observed in 
\textit{Example~1}--\textit{Example~3}---reflects the smoothness of the correction 
field: the network need only learn a small adjustment to 
an already accurate lift, and the resulting PDE residual 
is correspondingly small. The spikes occur when randomly 
sampled collocation points intersect the thin internal 
layer at $x_1 + x_2 = t$, where the solution gradients 
are of $\mathcal{O}(\text{Re})$. The weight schedule in Figure~\ref{fig:fig13b} 
illustrates the four-phase strategy unique to this 
example. The data weight remains fixed at 
$w_{\text{data}} = 1.0$ throughout all phases, in 
contrast to the decreasing schedule used in \textit{Example~1}--\textit{Example~3}. 
This design choice, motivated by the non-uniqueness of 
weak solutions to Burgers' equation, ensures that the 
network remains anchored to the physically correct 
solution branch at all times. The PDE weight increases 
progressively from $w_{\text{pde}} = 0.01$ in Phase~I 
through $w_{\text{pde}} = 0.1$ and $1.0$ in Phases~II 
and~III, culminating at $w_{\text{pde}} = 5.0$ in 
Phase~IV (epochs $3500$--$4999$). Despite this aggressive 
final PDE weight, the total loss remains dominated by the 
data component because the PDE residual is inherently 
small, confirming that the physics-informed training 
refines the solution without competing against the data 
fidelity objective.

\subsection{Example 5.}
Consider the following 2D time-dependent advection--diffusion equation with an L-shaped interior layer~\cite{bazilevs_yzb_2007}:
\begin{equation}
\begin{cases}
\displaystyle \frac{\partial u}{\partial t} + \mathbf{b} \cdot \nabla u - \nabla \cdot (\varepsilon \nabla u) = 0, & (x_1,x_2) \in \Omega, \quad t \in (0,t_{\text{f}}], \\[8pt]
u(0,x_1,x_2) = u_0(x_1,x_2), & (x_1,x_2) \in \Omega, \\[6pt]
u = g(t, x_1,x_2), & (x_1,x_2) \in \Gamma_D, \quad t \in [0,t_{\text{f}}],
\end{cases}
\end{equation}
where $\Omega = (0,1)^2$ is the unit square, $t_{\text{f}} = 0.25$, and the diffusivity is $\varepsilon = 10^{-8}$. The advection velocity is
\begin{equation}
\mathbf{b} = |\mathbf{b}|(\cos\theta, \sin\theta)^\top = \left( \frac{\sqrt{2}}{2}, \frac{\sqrt{2}}{2} \right)^\top,
\end{equation}
with a magnitude $\Vert\mathbf{b}\Vert = 1$ and an advection angle $\theta = 45^\circ$ relative to the mesh. The initial condition consists of an L-shaped block:
\begin{equation}
u_0(x_1,x_2) = 
\begin{cases}
1, & \text{if } (x_1,x_2) \in \Omega_L, \\[4pt]
0, & \text{otherwise},
\end{cases}
\end{equation}
where the L-shaped region $\Omega_\text{L}$ is defined by
\begin{equation}
\Omega_L = \left\{ (x_1,x_2) \,:\, 0 \leq x_1 \leq \tfrac{1}{2}L, \; 
0 \leq x_2 \leq \tfrac{1}{4}L \right\} \cup \left\{ (x_1,x_2) \,:\, 
0 \leq x_1 \leq \tfrac{1}{4}L, \; 0 \leq x_2 \leq \tfrac{1}{2}L \right\},
\end{equation}
with $L = 1$. The Dirichlet boundary conditions are
\begin{equation}
g(t, x_1, x_2) = 
\begin{cases}
1, & \text{if } x_2 = 0 \text{ and } 0 \leq x_1 \leq \tfrac{1}{2}L, \\[4pt]
1, & \text{if } x_1 = 0 \text{ and } 0 \leq x_2 \leq \tfrac{1}{2}L, \\[4pt]
0, & \text{otherwise on } \partial\Omega.
\end{cases}
\end{equation}
For the finite element setting, we set 
$n_{\text{el}} = 8192$ (structured triangulation of a 
$64 \times 64$ grid), $\varDelta t = 0.001$ ($N_t = 250$), 
$t_\text{f} = 0.25$, and Y$= 0.25$. Since no analytical 
solution is available for this problem, the assessment 
is based on qualitative comparison of the solution 
profiles and quantitative comparison of the $L^2$ norms 
between methods. The PINN architecture employs 
$n_h = 128$, $n_r = 8$, $n_\text{F} = 24$, and $\sigma = 4.0$. 
Training uses the last $K_s = 10$ temporal snapshots near 
$t_\text{f} = 0.25$, with batch size $2048$, initial learning 
rate $\alpha = 3 \times 10^{-4}$, and gradient clipping threshold 
$g_{\max} = 1.0$. The PDE residual is evaluated at 
$N_{\text{pde}} = 4096$ randomly sampled interior 
collocation points per epoch---the largest value 
among all examples, chosen to adequately resolve the 
complex interaction of the two internal layers forming 
along the outflow boundary.

Unlike the preceding examples, the training strategy for 
this problem adopts a progressive four-phase schedule in 
which all three weights evolve simultaneously. 
This design is motivated by the absence of an analytical 
solution: the SUPG-YZ$\beta$ reference data constitutes 
the only available ground truth, and the weight schedule 
must carefully balance data fidelity with physics 
enforcement to avoid steering the network toward 
solutions that satisfy the PDE residual but deviate 
from the FEM reference without any means of independent 
verification. The data weight increases from 
$w_{\text{data}} = 0.1$ in Phase~I to 
$w_{\text{data}} = 0.8$ in Phase~IV, reflecting the 
growing importance of anchoring to the reference solution 
as physics enforcement intensifies. The PDE weight 
increases from $w_{\text{pde}} = 0.5$ to 
$w_{\text{pde}} = 10.0$, progressively strengthening 
the physics constraint. The boundary condition weight 
increases from $w_{\text{bc}} = 0.1$ to 
$w_{\text{bc}} = 1.0$, ensuring that the 
non-homogeneous piecewise Dirichlet data---which cannot 
be incorporated through a simple lift function---remains 
well-satisfied as the PDE loss becomes dominant. 
The phase boundaries and adaptive weights are reported 
in Table~\ref{tab:phase_example_six}.

\begin{table}[htb]
\centering
\caption{Phase boundaries and adaptive weights used for solving \textit{Example 5}.}
\label{tab:phase_example_six}
\begin{tabular}{llccc}
\toprule
\text{Phase} & \text{Epochs} & $w_{\mathrm{data}}$ & $w_{\mathrm{pde}}$ & $w_{\mathrm{bc}}$ \\
\midrule
I   & 0--999     & 0.1  & 0.5   & 0.1 \\
II  & 1000--1999 & 0.35 & 0.8   & 0.3 \\
III & 2000--3499 & 0.5  & 5.0   & 0.8 \\
IV  & 3500--4999 & 0.8  & 10.0  & 1.0 \\
\bottomrule
\end{tabular}
\end{table}

\begin{figure}[h!]
	\centering
	\begin{subfigure}[b]{0.325\linewidth}
		\includegraphics[width=\linewidth]{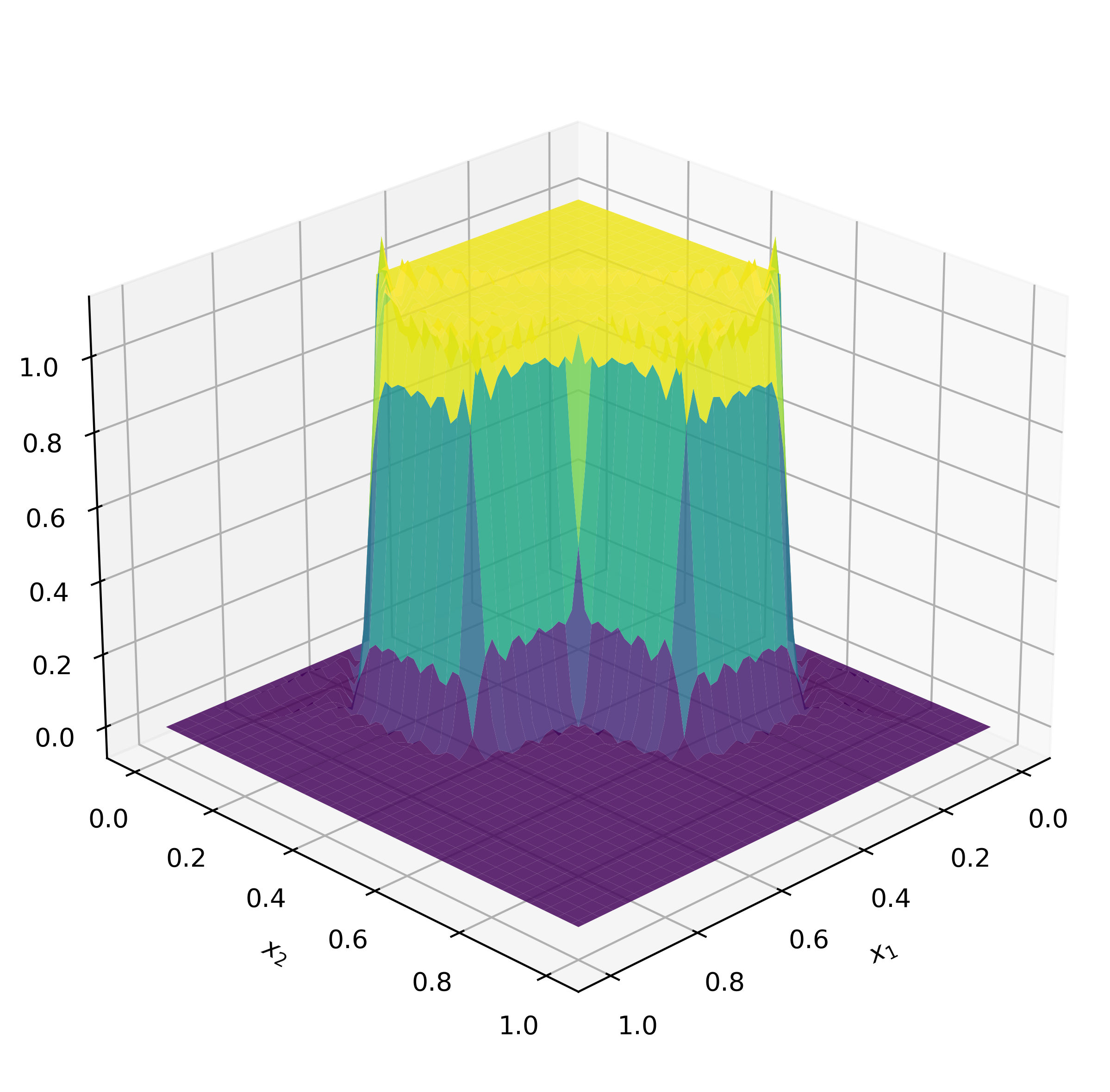}
		\caption{}
		\label{fig:fig14a}
	\end{subfigure}
	\hfill
	\begin{subfigure}[b]{0.325\linewidth}
		\includegraphics[width=\linewidth]{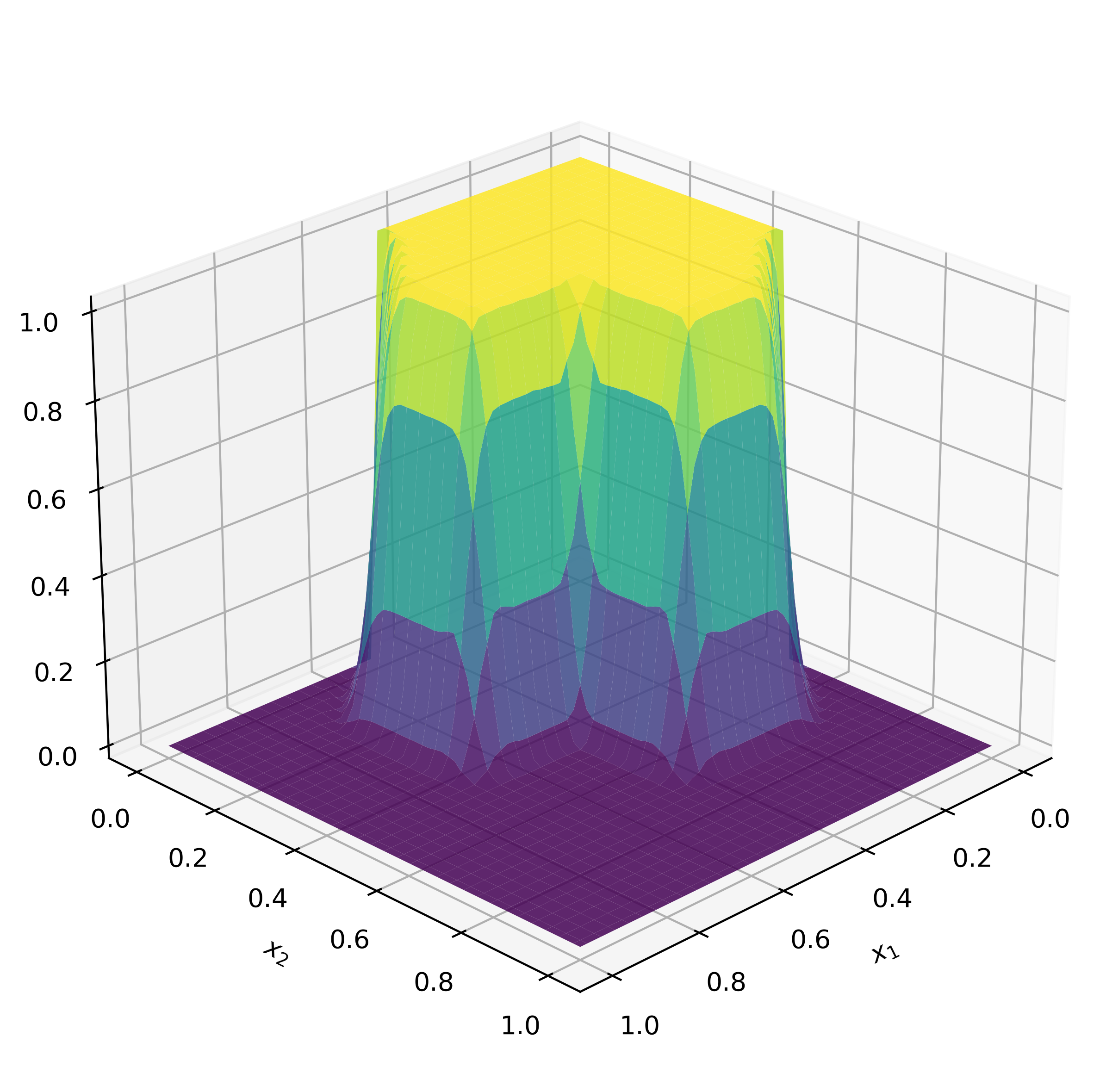}
		\caption{}
		\label{fig:fig14b}
	\end{subfigure}
	\hfill
	\begin{subfigure}[b]{0.325\linewidth}
		\includegraphics[width=\linewidth]{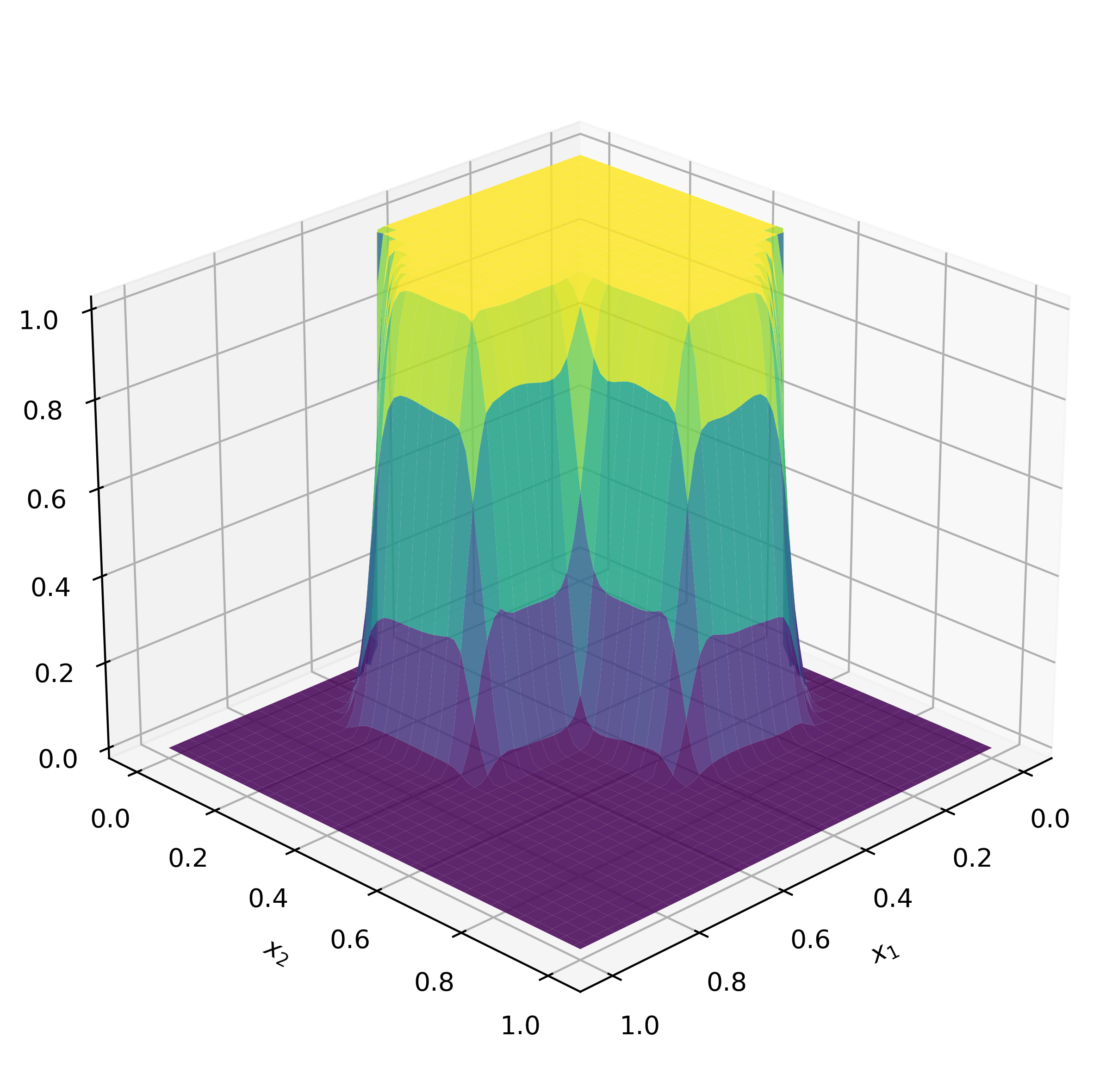}
		\caption{}
		\label{fig:fig14c}
	\end{subfigure}
	\caption{Three-dimensional surface plots of the solution 
		to \textit{Example~5} at $t_\text{f} = 0.25$ with 
		$\varepsilon = 10^{-8}$: 
		(\subref{fig:fig14a})~SUPG finite element approximation; 
		(\subref{fig:fig14b})~SUPG-YZ$\beta$ finite element 
		approximation; 
		(\subref{fig:fig14c})~hybrid PINN.}
	\label{fig:14_surface}
\end{figure}

Figure~\ref{fig:14_surface} displays the three-dimensional solution 
profiles at the terminal time $t_\text{f} = 0.25$ for \textit{Example}~5. Since 
no analytical solution is available for this problem, the assessment 
relies on qualitative comparison among the three numerical 
approximations. The SUPG solution (Figure~\ref{fig:fig14a}) 
captures the overall L-shaped block structure and its diagonal 
advection toward the upper-right corner; however, pronounced 
oscillatory artifacts are visible along the leading edges of the 
front, particularly near the corner where the horizontal and vertical 
arms of the L-shape meet. These crosswind oscillations, a 
characteristic deficiency of streamline-based stabilization in the 
presence of sharp transverse gradients, contaminate the flat plateau 
region where the solution should remain identically unity, as well as 
the surrounding zero-valued exterior. The addition of YZ$\beta$ 
shock-capturing (Figure~\ref{fig:fig14b}) substantially 
suppresses these oscillatory features, yielding a smoother and more 
monotone profile. The front edges are better resolved, the plateau 
region is flatter, and the spurious oscillations in the exterior are 
largely eliminated. However, a slight smearing of the sharp 
transitions is discernible, consistent with the isotropic artificial 
diffusion introduced by the shock-capturing operator. The hybrid PINN 
correction (Figure~\ref{fig:fig14c}) produces a surface profile 
that closely resembles the SUPG-YZ$\beta$ approximation in its 
overall structure while exhibiting two notable improvements: the solution remains 
within the physically admissible range $[0,1]$ without exhibiting the 
overshoot or undershoot artifacts present in the FEM approximations, 
and the front transitions appear marginally sharper along the leading 
edges of the advected L-shaped block. This improved boundedness is 
attributable to the combined effect of data-driven training from the 
well-behaved SUPG-YZ$\beta$ reference and the PDE residual 
regularization, which together steer the network away from 
nonphysical extrema.

\begin{figure}[h!]
	\centering
	\begin{subfigure}[b]{0.325\linewidth}
		\includegraphics[width=\linewidth]{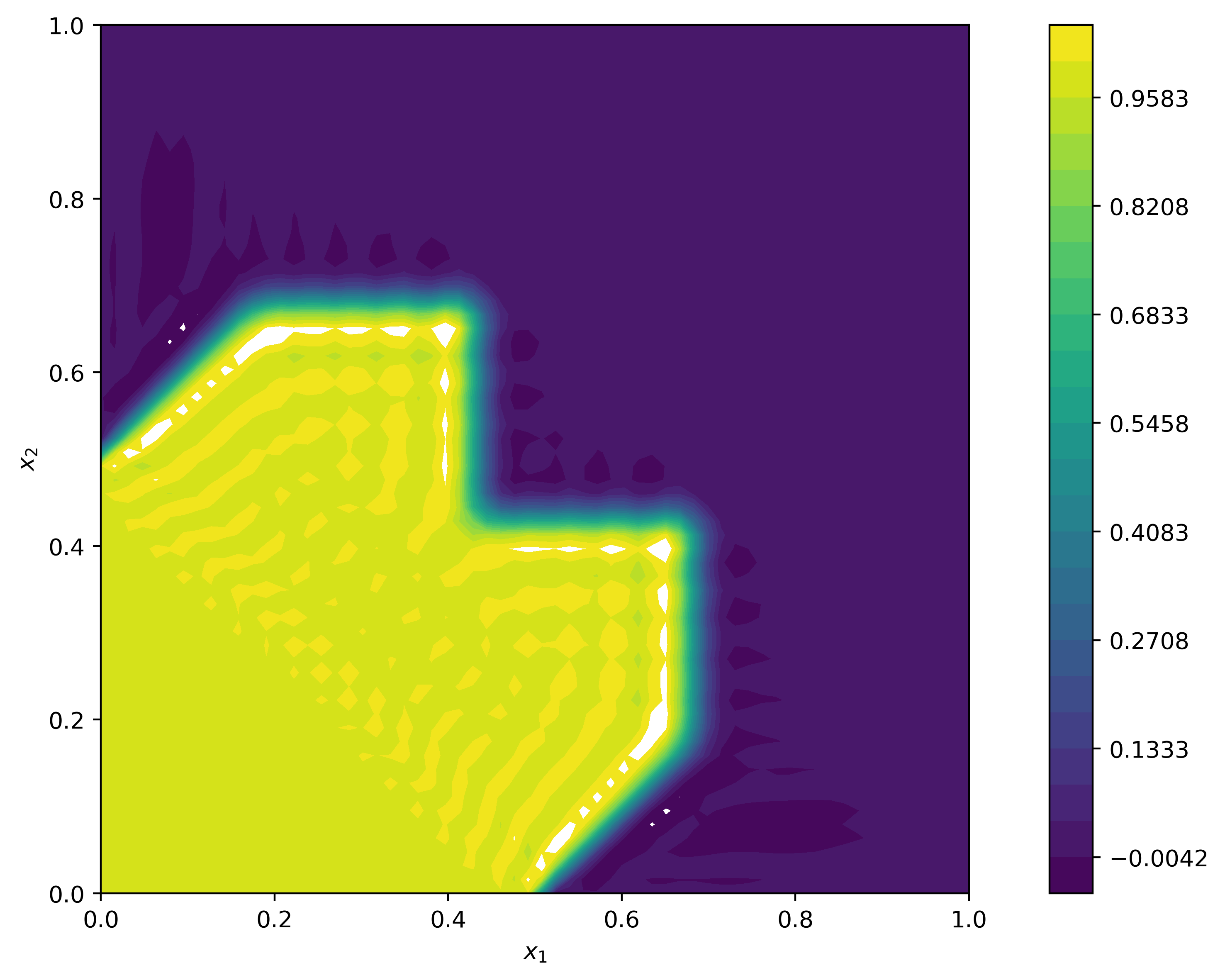}
		\caption{}
		\label{fig:fig15a}
	\end{subfigure}
	\hfill
	\begin{subfigure}[b]{0.325\linewidth}
		\includegraphics[width=\linewidth]{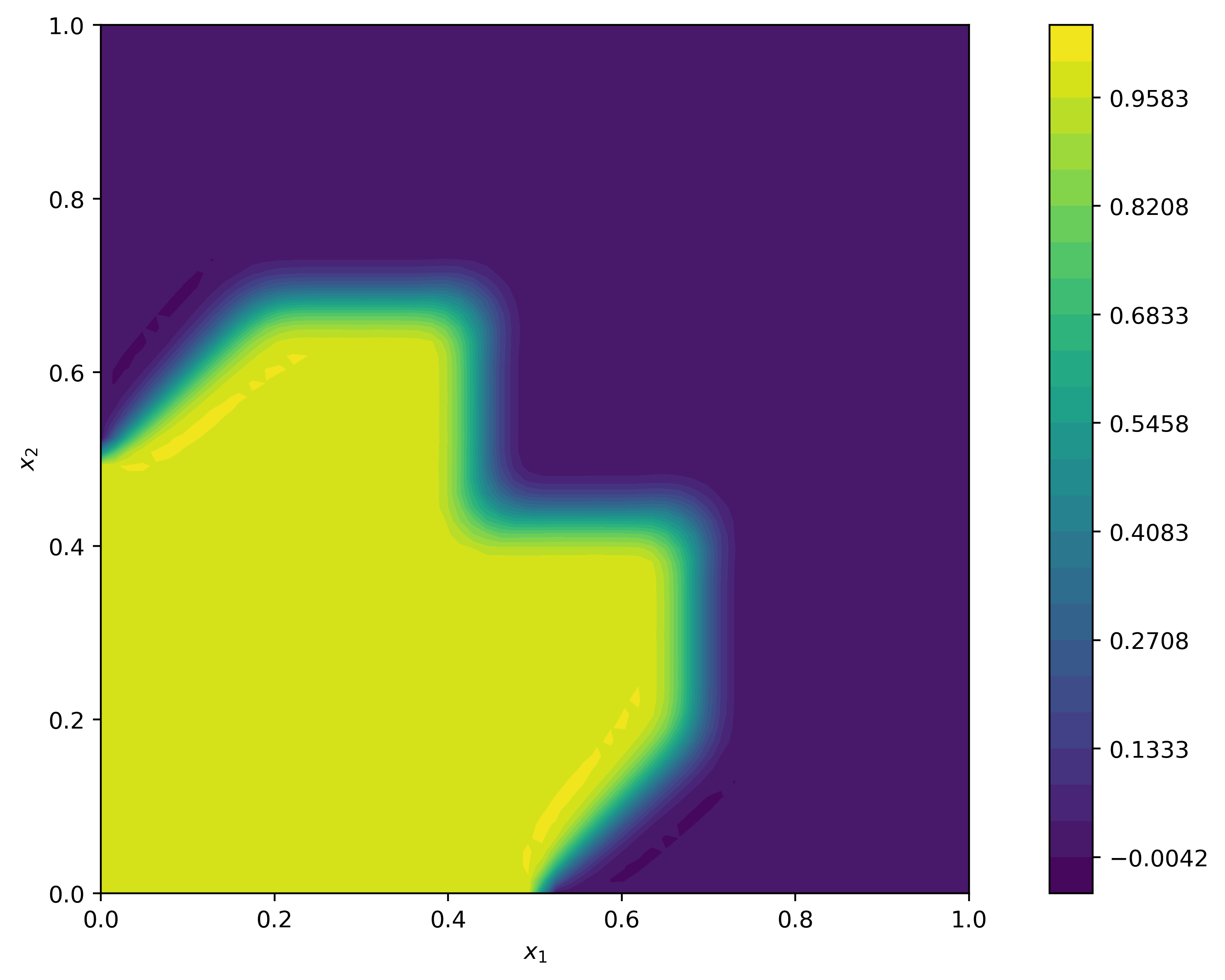}
		\caption{}
		\label{fig:fig15b}
	\end{subfigure}
	\hfill
	\begin{subfigure}[b]{0.325\linewidth}
		\includegraphics[width=\linewidth]{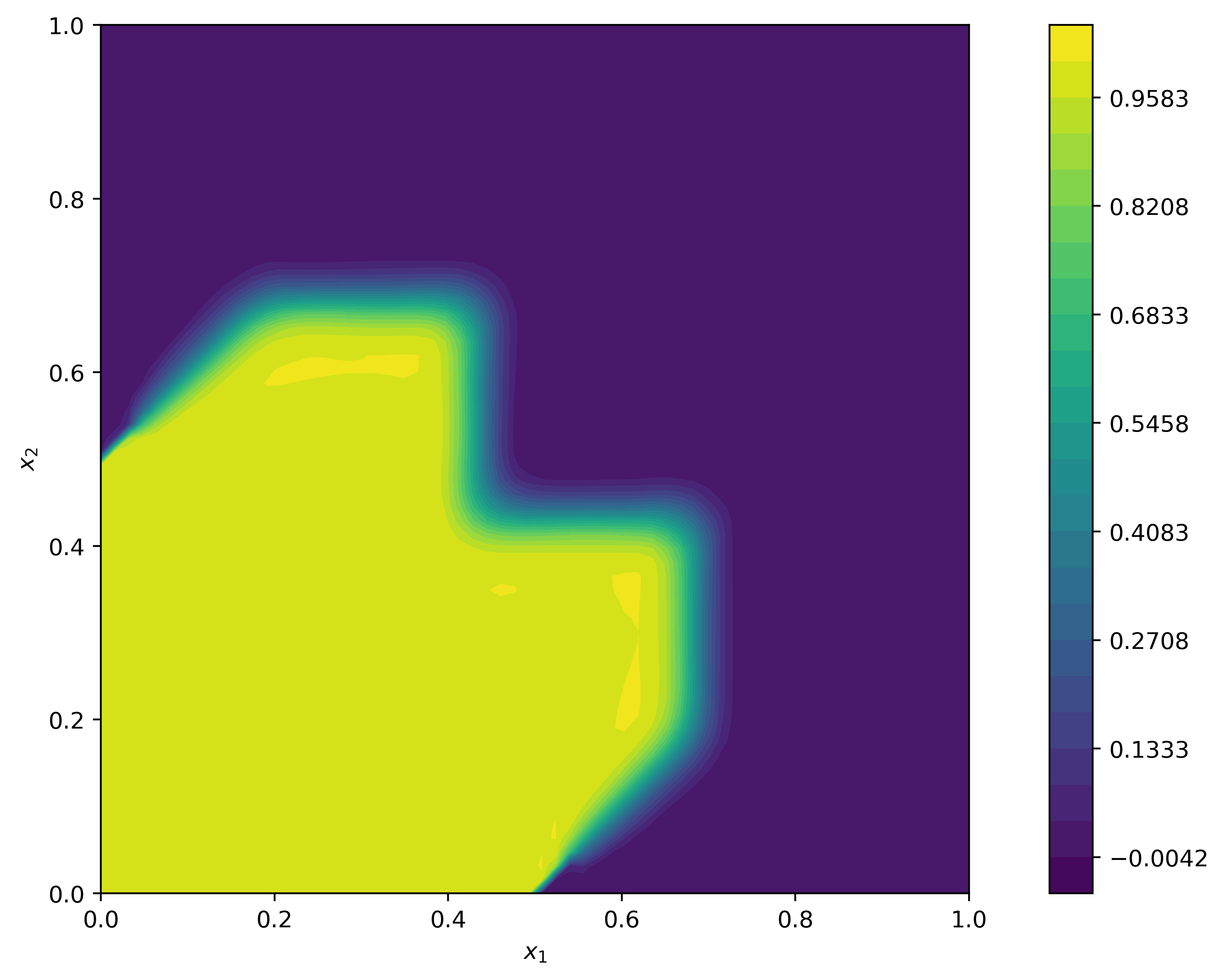}
		\caption{}
		\label{fig:fig15c}
	\end{subfigure}
	\caption{Contour plots of the solution 
		to \textit{Example~5} at $t_\text{f} = 0.25$ with 
		$\varepsilon = 10^{-8}$: 
		(\subref{fig:fig15a})~SUPG finite element approximation; 
		(\subref{fig:fig15b})~SUPG-YZ$\beta$ finite element 
		approximation; 
		(\subref{fig:fig15c})~hybrid PINN.}
	\label{fig:15_surface}
\end{figure}

The contour plots in Figure~\ref{fig:15_surface} provide a 
top-down view that more clearly exposes the spatial structure of 
the advected L-shaped front and the differences among the three 
numerical approximations. The SUPG solution 
(Figure~\ref{fig:fig15a}) reveals pronounced crosswind 
oscillations manifesting as a series of concentric ripple-like 
contour distortions along the leading edges of the front. These 
artifacts, which extend well into both the plateau region and the 
zero-valued exterior, are a hallmark of insufficient crosswind 
diffusion for convection-dominated transport of discontinuous 
initial data on structured meshes. The oscillatory pattern is 
particularly severe near the re-entrant corner of the L-shape, 
where the interaction of the horizontal and vertical fronts 
produces competing gradient directions that the streamline-based 
SUPG stabilization cannot adequately resolve. The SUPG-YZ$\beta$ solution (Figure~\ref{fig:fig15b}) 
eliminates these oscillatory features entirely, producing a smooth 
contour field with well-defined front boundaries. However, the 
front transitions are noticeably broadened compared to the initial 
sharp discontinuity, reflecting the isotropic artificial diffusion 
introduced by the shock-capturing operator. The plateau region 
behind the front maintains a nearly uniform value close to unity, 
confirming the effectiveness of the YZ$\beta$ mechanism in 
preserving the solution structure away from steep gradients. The hybrid PINN contour (Figure~\ref{fig:fig15c}) closely 
reproduces the SUPG-YZ$\beta$ pattern while exhibiting subtle 
differences along the front edges: the transition zone appears 
marginally narrower, suggesting that the PDE residual enforcement 
partially counteracts the numerical diffusion introduced by the 
shock-capturing operator. Furthermore, the hybrid PINN solution remains within $[0, 1]$ 
throughout the domain, avoiding the minor undershoot of $-0.0042$ 
visible in the colorbar of the FEM approximations---a consequence of 
the data fidelity and PDE residual losses jointly penalizing 
nonphysical solution values during training.

\begin{figure}[h!]
	\centering
	\begin{subfigure}[b]{0.245\linewidth}
		\includegraphics[width=\linewidth]{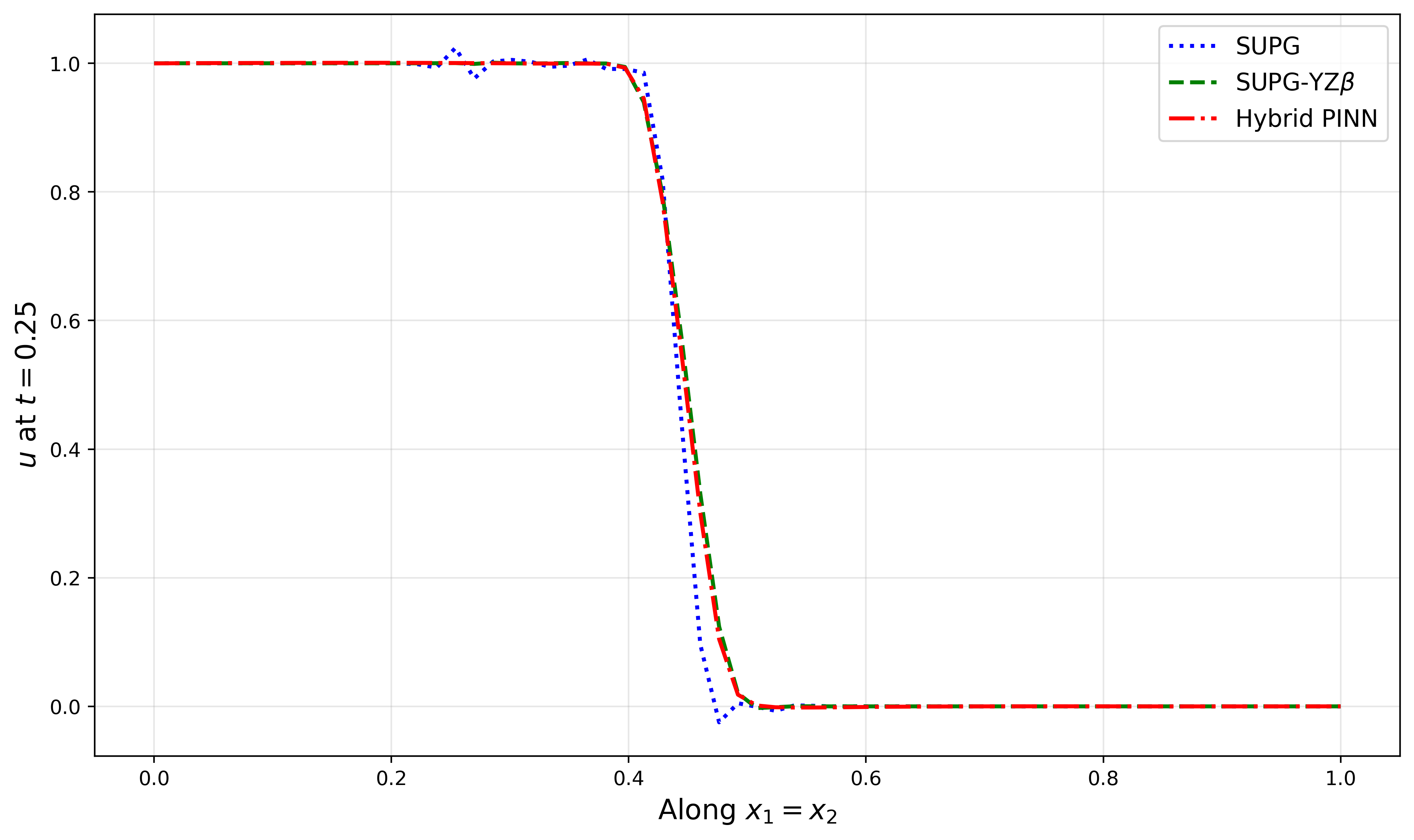}
		\caption{}
		\label{fig:fig16a}
	\end{subfigure}
	\hfill
	\begin{subfigure}[b]{0.245\linewidth}
		\includegraphics[width=\linewidth]{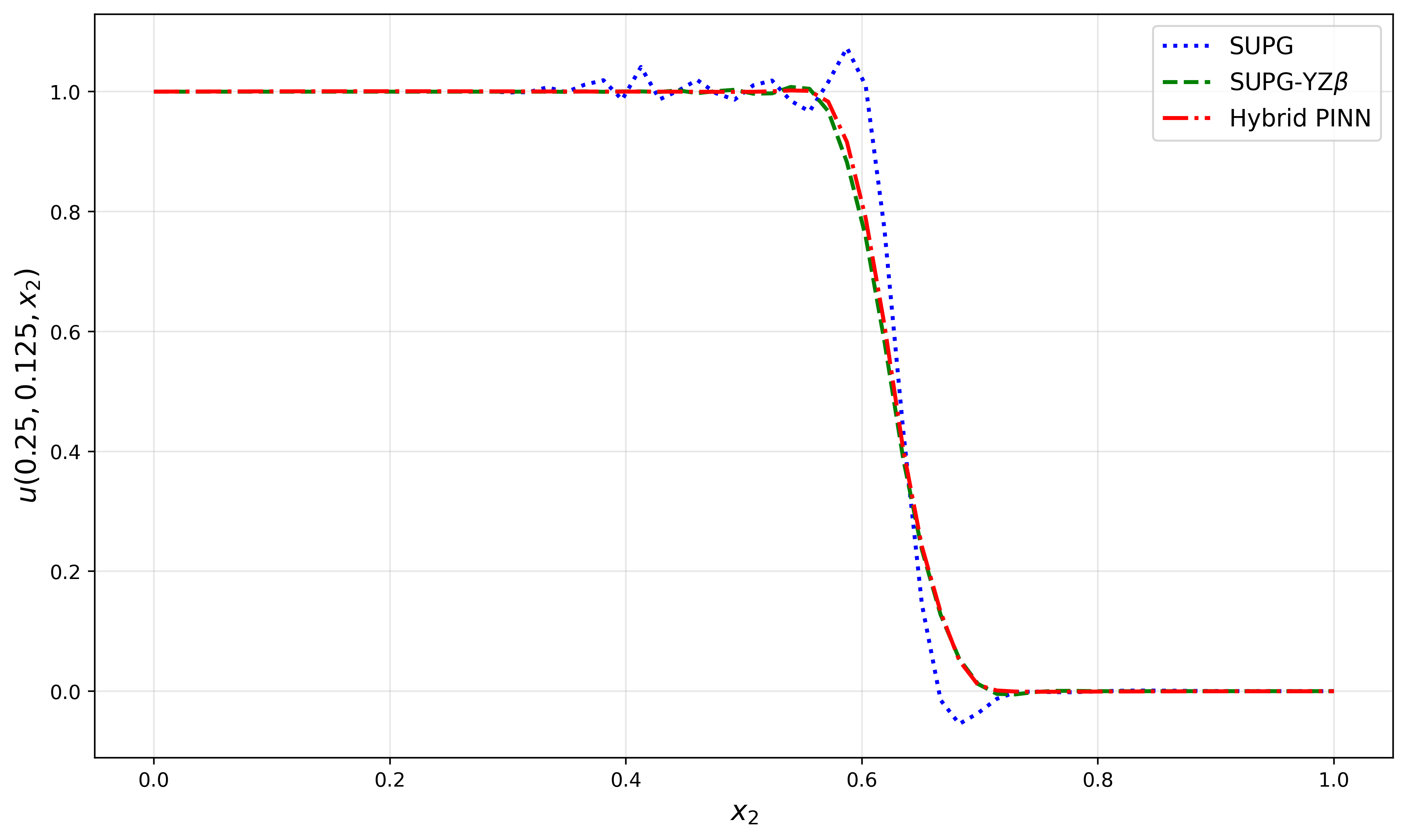}
		\caption{}
		\label{fig:fig16b}
	\end{subfigure}
	\hfill
	\begin{subfigure}[b]{0.245\linewidth}
		\includegraphics[width=\linewidth]{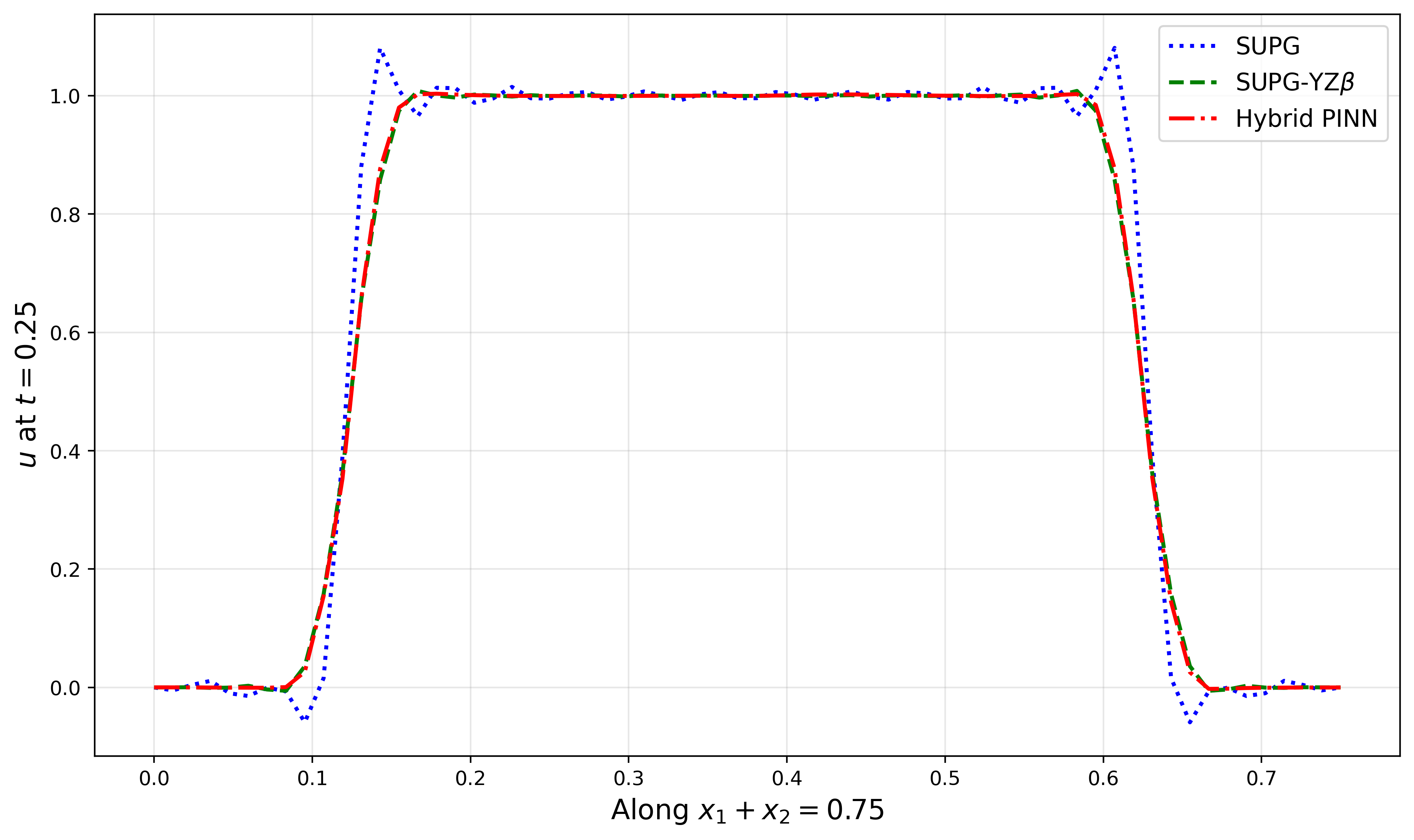}
		\caption{}
		\label{fig:fig16c}
	\end{subfigure}
	\hfill
	\begin{subfigure}[b]{0.245\linewidth}
		\includegraphics[width=\linewidth]{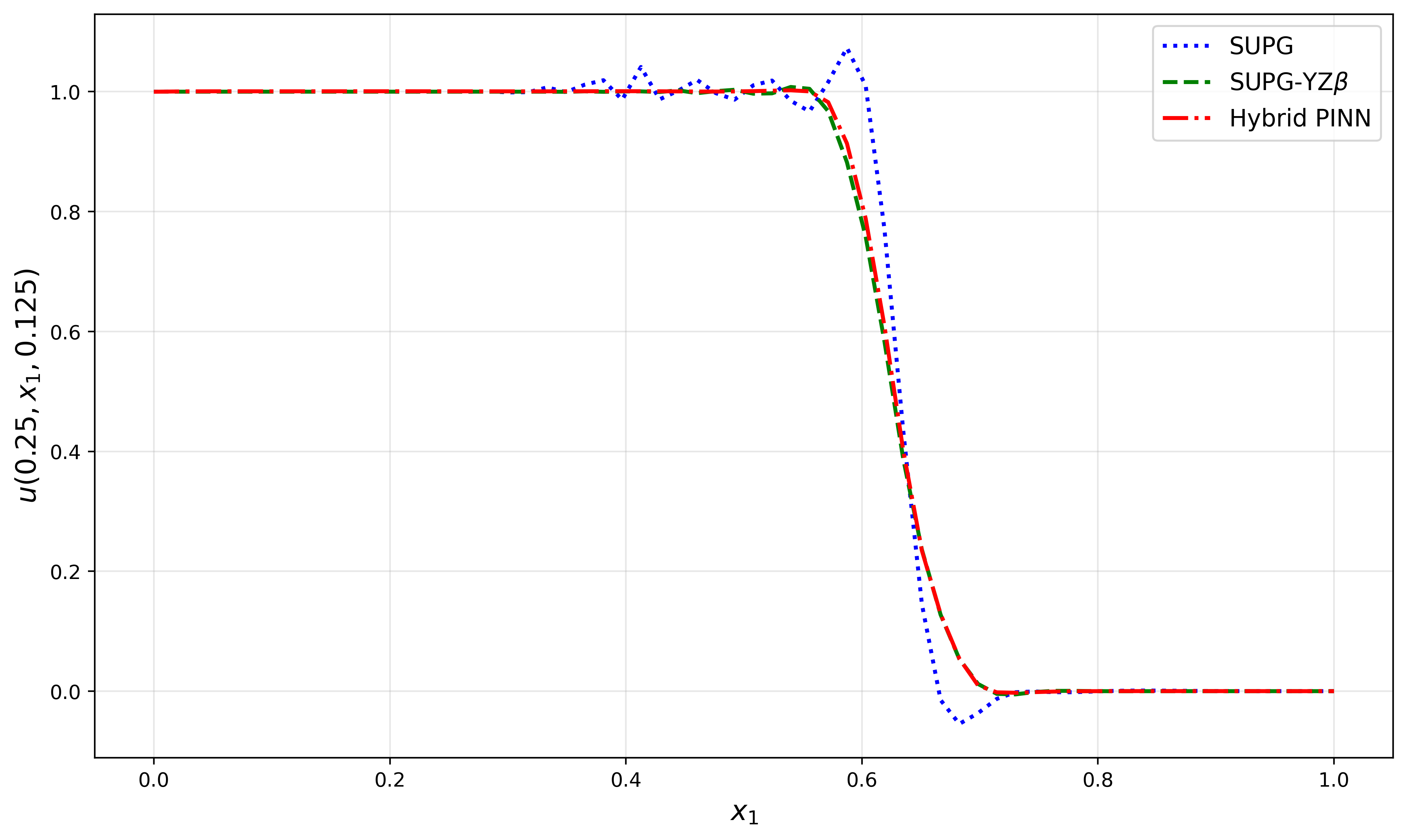}
		\caption{}
		\label{fig:fig16d}
	\end{subfigure}
	\caption{Cross-sectional profiles for \textit{Example~5} 
		at $t_\text{f} = 0.25$ with $\varepsilon = 10^{-8}$, 
		comparing SUPG, SUPG-YZ$\beta$, and hybrid PINN: 
		(\subref{fig:fig16a})~along the diagonal $x_1 = x_2$; 
		(\subref{fig:fig16b})~along $x_1 = 0.125$ (through the 
		vertical arm); 
		(\subref{fig:fig16c})~along $x_1 + x_2 = 0.75$ 
		(perpendicular to the advection direction); 
		(\subref{fig:fig16d})~along $x_2 = 0.125$ (through the 
		horizontal arm).}
	\label{fig:16_cross}
\end{figure}

Figure~\ref{fig:16_cross} presents one-dimensional cross-sectional 
profiles extracted along four representative lines at the terminal 
time $t_\text{f} = 0.25$, providing a quantitative comparison among 
the three methods. The diagonal slice along $x_1 = x_2$ 
(Figure~\ref{fig:fig16a}) intersects the advected L-shaped front 
along the convection direction $\mathbf{b} = (\sqrt{2}/2, 
\sqrt{2}/2)^\top$. All three methods capture the transition from 
$u \approx 1$ to $u \approx 0$ near $s \approx 0.45$; however, the 
SUPG solution exhibits visible oscillatory overshoots in the 
vicinity of the front, while the SUPG-YZ$\beta$ and hybrid PINN 
profiles are nearly monotone. The PINN profile tracks the 
SUPG-YZ$\beta$ solution closely, with a marginally sharper 
transition at the leading edge. The vertical slice at $x_1 = 0.125$ 
(Figure~\ref{fig:fig16b}) cuts through the vertical arm of the 
L-shape. The plateau at $u \approx 1$ extends to $x_2 \approx 0.6$, 
consistent with the advection of the initial block from 
$x_2 \in [0, 0.5]$. The SUPG solution displays pronounced 
undershoots below zero near the trailing edge of the front, while 
both the SUPG-YZ$\beta$ and hybrid PINN solutions maintain 
non-negative profiles. The hybrid PINN achieves a slightly steeper 
descent at the front compared to the SUPG-YZ$\beta$ approximation.
The slice along $x_1 + x_2 = 0.75$ 
(Figure~\ref{fig:fig16c}), oriented perpendicular to the 
advection direction, provides the most revealing comparison. This 
cross-section intersects the front transversely, where crosswind 
effects are most pronounced. The SUPG solution shows significant 
oscillatory artifacts along the entire profile, reflecting the 
well-known crosswind diffusion deficiency. The SUPG-YZ$\beta$ 
profile is substantially smoother, though slightly broadened. The hybrid PINN follows the SUPG-YZ$\beta$ profile while maintaining 
boundedness within $[0, 1]$, consistent with the regularizing effect 
of the physics-informed training. The horizontal slice at $x_2 = 0.125$ 
(Figure~\ref{fig:fig16d}) passes through the horizontal arm of 
the L-shape. The front near $x_1 \approx 0.65$ shows similar 
behavior: the SUPG solution oscillates, the SUPG-YZ$\beta$ solution 
provides a smooth but diffused transition, and the hybrid PINN 
closely tracks the stabilized FEM reference. Across all four cross-sections, the PINN correction consistently 
preserves the overall solution structure of the SUPG-YZ$\beta$ 
approximation while empirically satisfying the maximum principle---a 
property that neither FEM formulation guarantees by construction. 
This behavior, while not enforced architecturally, arises naturally 
from the combination of training on bounded reference data and PDE 
residual minimization in the interior.

\begin{figure}[h!]
	\centering
	\begin{subfigure}[b]{0.48\linewidth}
		\includegraphics[width=\linewidth]{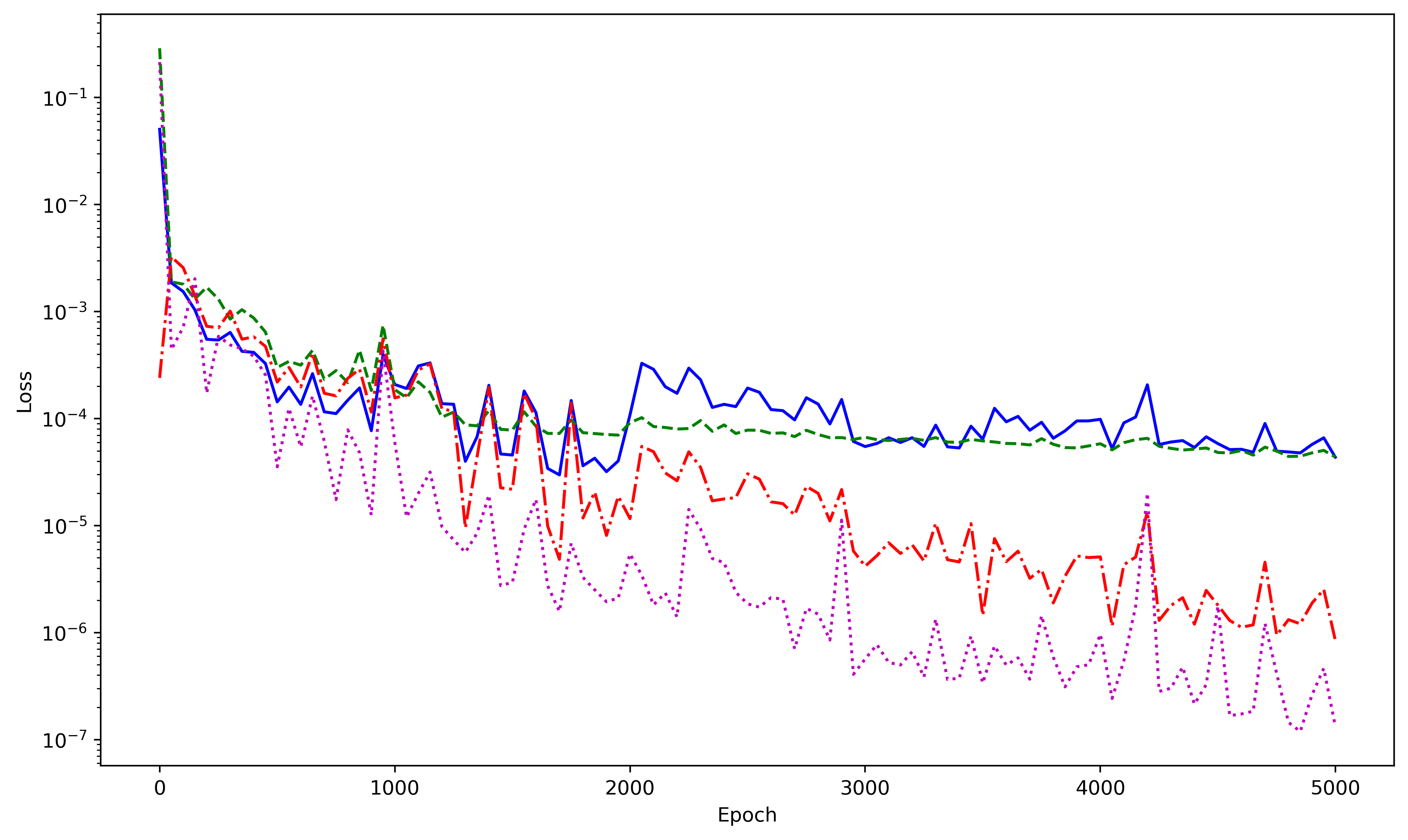}
		\caption{}
		\label{fig:fig17a}
	\end{subfigure}
	\hfill
	\begin{subfigure}[b]{0.48\linewidth}
		\includegraphics[width=\linewidth]{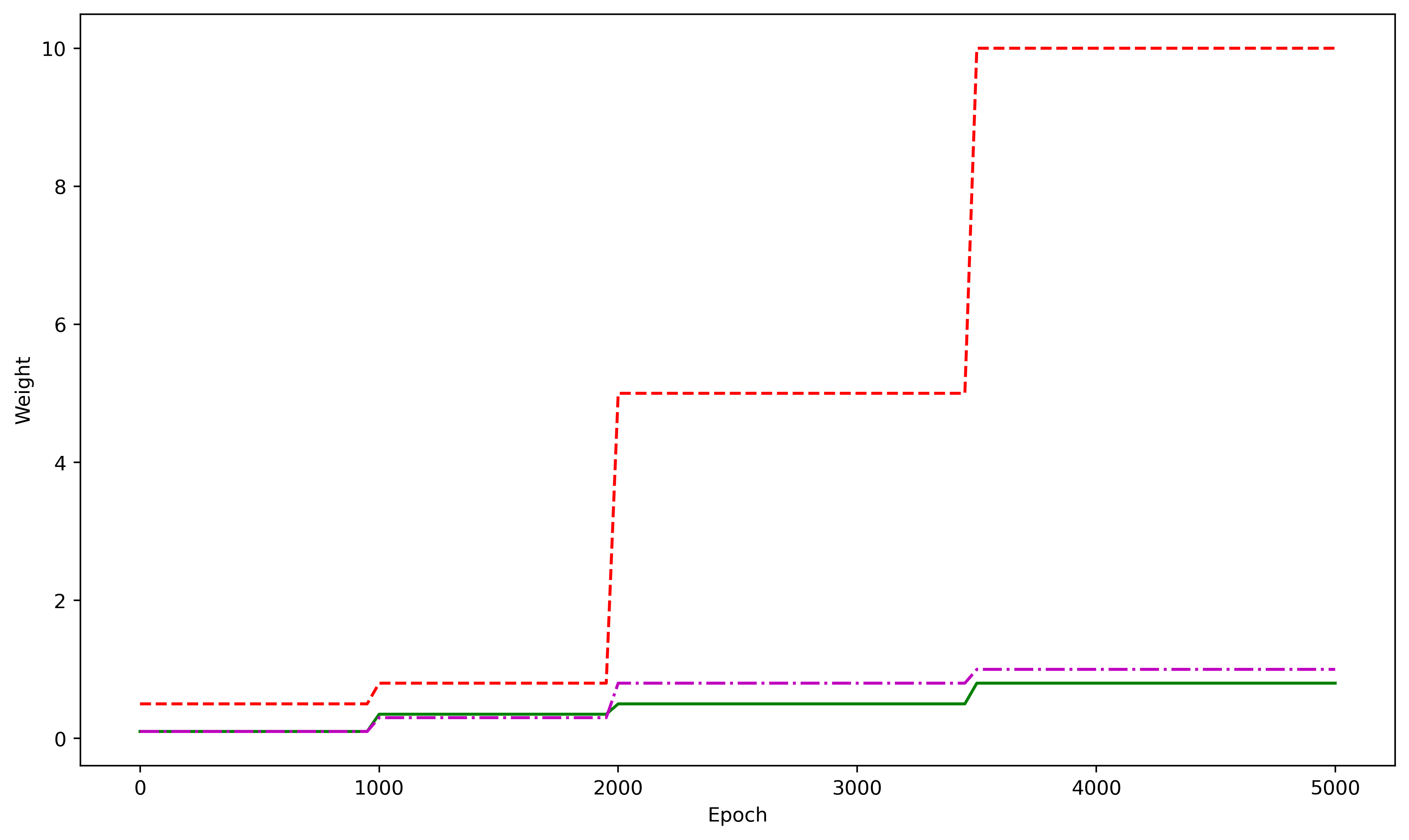}
		\caption{}
		\label{fig:fig17b}
	\end{subfigure}
	\caption{Training diagnostics for \textit{Example~5} with 
		$\varepsilon = 10^{-8}$: 
		(\subref{fig:fig17a})~evolution of the individual loss 
		components (total, data, PDE, and BC) over $5000$ 
		training epochs on a logarithmic scale; 
		(\subref{fig:fig17b})~adaptive weight schedule 
		$w_{\text{data}}$, $w_{\text{pde}}$, and 
		$w_{\text{bc}}$ governing the four-phase training 
		strategy.}
	\label{fig:17_training}
\end{figure}

The training diagnostics for Example~5 are presented in 
Figure~\ref{fig:17_training}. The loss evolution in 
Figure~\ref{fig:17_training}a reveals the distinct behavior of each 
component across the four training phases. During Phase~I (epochs 
0--999), where the PDE weight is moderate at $w_{\text{pde}} = 0.5$ 
and the data weight is low at $w_{\text{data}} = 0.1$, all loss 
components decrease rapidly from $\mathcal{O}(10^{-1})$ to 
$\mathcal{O}(10^{-4})$, indicating that the network successfully 
learns the global solution structure from the SUPG-YZ$\beta$ 
reference snapshots while simultaneously beginning to satisfy the 
PDE constraint. The data loss stabilizes near 
$\mathcal{O}(10^{-4})$ by the end of Phase~I and continues to 
decrease gradually throughout the remaining phases, confirming that 
the progressively increasing data weight 
($w_{\text{data}}: 0.1 \to 0.8$) effectively anchors the network 
to the reference solution. The PDE loss exhibits intermittent spikes 
spanning approximately two orders of magnitude, from 
$\mathcal{O}(10^{-6})$ to $\mathcal{O}(10^{-4})$, persisting 
throughout training. These bursts arise when randomly sampled 
collocation points intersect the sharp L-shaped front, where the 
solution gradients are steep and the PDE residual magnitudes 
increase accordingly. Despite this volatile behavior, the PDE loss 
trend decreases steadily over the course of training. The BC loss 
remains consistently low at $\mathcal{O}(10^{-5})$--$\mathcal{O}(10^{-7})$ 
throughout all phases, reflecting the effectiveness of the explicit 
boundary condition enforcement for the piecewise-discontinuous 
Dirichlet data. The weight schedule in Figure~\ref{fig:17_training}b illustrates the 
four-phase strategy adopted for this problem. The PDE weight 
undergoes the most dramatic evolution, increasing from 
$w_{\text{pde}} = 0.5$ in Phase~I through $w_{\text{pde}} = 0.8$ 
in Phase~II to $w_{\text{pde}} = 5.0$ in Phase~III, culminating at 
$w_{\text{pde}} = 10.0$ in Phase~IV (epochs 3500--4999). This 
aggressive final PDE weight is necessary to drive the physics-informed 
refinement beyond mere data interpolation. Concurrently, the data 
weight increases from $w_{\text{data}} = 0.1$ to 
$w_{\text{data}} = 0.8$, ensuring that the network remains anchored 
to the SUPG-YZ$\beta$ reference as the PDE constraint intensifies---a 
design choice motivated by the absence of an analytical solution 
for independent verification. The boundary condition weight increases 
gradually from $w_{\text{bc}} = 0.1$ to $w_{\text{bc}} = 1.0$, 
reflecting the growing importance of maintaining boundary fidelity 
for the non-homogeneous piecewise Dirichlet data as the training 
transitions toward physics-dominant refinement.

\section{Conclusion and Future Work} \label{sec:section6}
This work presented a hybrid computational framework that 
couples stabilized finite element methods with physics-informed 
neural networks for the numerical solution of transient 
convection-dominated CDR-type partial differential equations. 
The methodology employs the SUPG formulation augmented with 
the YZ$\beta$ shock-capturing operator to generate 
high-fidelity reference solutions, which subsequently serve 
as training data for a deep neural network equipped with 
residual blocks and random Fourier feature embeddings. Rather 
than training over the entire spatiotemporal domain, the PINN 
correction targets the last $K_s$ temporal snapshots near the 
terminal time, where the network assimilates finite element 
data while simultaneously enforcing the governing PDE in 
strong form through automatic differentiation.

A central ingredient of the proposed approach is the selective 
physics enforcement strategy, which restricts the PDE residual 
loss to interior collocation points satisfying a distance-based 
exclusion criterion. This design choice prevents the 
physics-informed training from interfering with sharp boundary 
and internal layer structures that are already well-resolved 
by the stabilized finite element discretization. The 
multiphase adaptive weight evolution---progressing from 
data-dominant through transitional to physics-dominant 
training---proved essential for balancing data fidelity with 
physical consistency and for avoiding early training 
instabilities that typically arise when physics terms dominate 
before a reasonable solution approximation has been 
established.

The framework was validated on five benchmark problems spanning 
one- and two-dimensional settings: a linear parabolic equation 
with a boundary layer (\textit{Example~1}), a 
time-dependent CDR equation with a hump changing its height 
(\textit{Example~2}), a traveling internal layer problem (\textit{Example~3}), the viscous Burgers' 
equation with a moving shock front (\textit{Example~4}), 
and an advection--diffusion equation with a discontinuous 
L-shaped interior layer (\textit{Example~5}). Across all test cases, 
the hybrid PINN correction consistently improved upon the 
standalone SUPG-YZ$\beta$ finite element solution, 
particularly in regions exhibiting steep gradients, moving 
fronts, and localized oscillatory artifacts. The results 
demonstrate that the proposed methodology inherits the 
robustness and stability of classical stabilized formulations 
while leveraging the representational flexibility of neural 
networks to achieve enhanced solution quality.

Several directions for future research emerge naturally from 
this study. First, the extension to three-dimensional 
problems represents a critical step toward practical 
engineering applications, including computational fluid 
dynamics and heat transfer in complex geometries. While the 
network architecture and training strategy are 
dimension-agnostic by construction---the spatiotemporal input 
$\mathbf{z} = (t, \mathbf{x}) \in \mathbb{R}^{n_\text{sd}+1}$ and the 
Fourier feature embedding accommodate arbitrary spatial 
dimension---the computational cost of generating 
three-dimensional SUPG-YZ$\beta$ reference data and the 
increased collocation requirements for PDE residual evaluation 
introduce additional challenges that warrant systematic 
investigation.

Second, the application to coupled systems of PDEs, such as 
the incompressible Navier--Stokes equations, 
magnetohydrodynamic flows, and multi-species reactive 
transport, constitutes a compelling extension. In such 
settings, inter-equation coupling introduces additional 
nonlinearities, and the stabilization and shock-capturing 
parameters must account for the distinct characteristic scales 
of each governing equation. The diagonal scaling matrix 
$\mathbf{Y}$ and residual vector $\mathbf{Z}$ discussed in 
Remarks~5 and~6 provide a natural starting point for 
generalizing the YZ$\beta$ mechanism to vector-valued 
problems, while the PINN component can be extended to 
multi-output architectures that simultaneously approximate 
all solution fields.

Third, the adoption of space--time finite element formulations (see, for example, \cite{Hughes1988, Mittal1992, Tezduyar06b, Tezduyar2006st, Steinbach2015})
represents a natural and methodologically significant 
extension. In the current framework, spatial and temporal 
discretizations are treated separately: the spatial domain is 
discretized with stabilized finite elements while time 
integration is performed via the backward Euler scheme. 
Space--time formulations, by contrast, treat the temporal 
dimension on an equal footing with the spatial coordinates, 
constructing finite element meshes over the full 
spatiotemporal domain 
$\mathcal{Q} = (t_0, t_f]  \times \Omega$ and solving the 
resulting variational problem in a unified manner. This 
approach offers several potential advantages for the proposed 
hybrid framework: (i)~it enables local refinement 
simultaneously in space and time, which is particularly 
beneficial for problems with moving fronts and transient 
layer structures; (ii)~SUPG stabilization and YZ$\beta$ 
shock-capturing extend naturally to the space--time setting, 
since the stabilization parameter $\tau_\text{SUPG}$ already 
incorporates both spatial and temporal scales through 
Eq.~\eqref{eq:shakib_tau2}; and (iii)~the PINN component, 
which already operates on spatiotemporal input coordinates 
$\mathbf{z} = (t, \mathbf{x})$ and computes temporal 
derivatives via automatic differentiation, is inherently 
compatible with a space--time discretization without 
architectural modification. Investigating how space--time 
SUPG-YZ$\beta$ reference solutions interact with the 
selective physics enforcement strategy---particularly 
whether the distance-based exclusion criterion should be 
generalized to a spatiotemporal metric---constitutes an 
important open question.

Fourth, and perhaps most significantly, the data-driven 
optimization of stabilization and shock-capturing parameters 
through PINNs offers a promising avenue for eliminating the 
empirical calibration that currently underpins the selection 
of $\tau_\text{SUPG}$ and $\nu_\text{SC}$. Rather than 
relying on analytical expressions, one could train auxiliary 
neural networks to predict element-wise optimal stabilization 
parameters by minimizing a loss functional that penalizes 
both PDE residuals and solution oscillations. Such an 
approach would shift the parameter selection from heuristic 
formulas to a learning-based paradigm, potentially yielding 
superior performance across a broader class of problems 
without manual intervention.

Additional avenues include the incorporation of adaptive mesh 
refinement strategies guided by PINN-based error indicators, 
the exploration of more principled dynamic loss weighting 
mechanisms such as gradient-norm balancing or 
uncertainty-based scaling to replace the current empirically 
tuned multiphase schedule, and the investigation of transfer 
learning techniques that leverage trained models from one 
problem configuration to accelerate convergence for related 
scenarios with different parameter values or geometries.

\section*{Author declarations}
\subsection*{CRediT authorship contribution statement}
\textbf{SC:} Conceptualization; Formal analysis; Funding acquisition; Investigation; Methodology; Software; Validation; Visualization;  Writing - original draft; Writing - review \& editing. 
\textbf{ÖU:} Conceptualization; Formal analysis; Investigation; Writing - review \& editing.
\textbf{SN:} Conceptualization; Formal analysis; Investigation; Writing - review \& editing.

\subsection*{Conflict of interest}
The authors have no conflicts to disclose.

\subsection*{Acknowledgment}
The first author of this study was supported by the Scientific and Technological Research Council of Türkiye (TUBITAK) under Grant Number 225M468. The authors thank TUBITAK for its support. 

\bibliographystyle{elsarticle-num-names}
\bibliography{myBiblio}
		
	\end{document}